\documentclass{article}
\usepackage{amsmath}
\usepackage{amssymb}
\usepackage{amsthm}
\usepackage{stmaryrd}
\usepackage{cite}
\usepackage{euscript}
\usepackage[arrow,curve,matrix,arc,2cell]{xy}
\UseAllTwocells
\usepackage[utf8]{inputenc}
\usepackage[unicode]{hyperref}
\DeclareFontFamily{U}{rsfs}{} \DeclareFontShape{U}{rsfs}{n}{it}{<->
rsfs10}{} \DeclareSymbolFont{mscr}{U}{rsfs}{n}{it}
\DeclareSymbolFontAlphabet{\scr}{mscr}
\def\mathscr{\scr}
\begin{document}
\def\e#1\e{\begin{equation}#1\end{equation}}
\def\ea#1\ea{\begin{align}#1\end{align}}
\def\eq#1{{\rm(\ref{#1})}}
\theoremstyle{plain}
\newtheorem{thm}{Theorem}[section]
\newtheorem{lem}[thm]{Lemma}
\newtheorem{prop}[thm]{Proposition}
\newtheorem{cor}[thm]{Corollary}
\theoremstyle{definition}
\newtheorem{dfn}[thm]{Definition}
\newtheorem{ex}[thm]{Example}
\newtheorem{rem}[thm]{Remark}
\numberwithin{figure}{section}
\numberwithin{equation}{section}
\def\dim{\mathop{\rm dim}\nolimits}
\def\codim{\mathop{\rm codim}\nolimits}
\def\depth{\mathop{\rm depth}\nolimits}
\def\ind{\mathop{\rm ind}\nolimits}
\def\Ker{\mathop{\rm Ker}}
\def\Coker{\mathop{\rm Coker}}
\def\Hol{\mathop{\rm Hol}}
\def\Re{\mathop{\rm Re}}
\def\SU{\mathop{\rm SU}}
\def\Spin{\mathop{\rm Spin}}
\def\U{\mathbin{\rm U}}
\def\Re{\mathop{\rm Re}}
\def\Hess{\mathop{\rm Hess}}
\def\Crit{\mathop{\rm Crit}}
\def\rank{\mathop{\rm rank}\nolimits}
\def\Hom{\mathop{\rm Hom}\nolimits}
\def\id{{\mathop{\rm id}\nolimits}}
\def\Id{{\mathop{\rm Id}\nolimits}}
\def\End{{\mathop{\rm End}\nolimits}}
\def\Sch{\mathop{\bf Sch}\nolimits}
\def\ev{{\rm ev}}
\def\Spec{\mathop{\rm Spec}\nolimits}
\def\Sets{{\mathop{\bf Sets}}}
\def\Man{{\mathop{\bf Man}}}
\def\Manb{{\mathop{\bf Man^b}}}
\def\Manc{{\mathop{\bf Man^c}}}
\def\Manra{{\mathop{\bf Man_{ra}}}}
\def\Manrab{{\mathop{\bf Man_{ra}^b}}}
\def\Manrac{{\mathop{\bf Man_{ra}^c}}}
\def\Manraab{{\mathop{\bf Man_{ra}^{ab}}}}
\def\Manraac{{\mathop{\bf Man_{ra}^{ac}}}}
\def\Mancin{{\mathop{\bf Man^c_{in}}}}
\def\Mancst{{\mathop{\bf Man^c_{st}}}} 
\def\cManc{{\mathop{\bf\check{M}an^c}}}
\def\cMancin{{\mathop{\bf\check{M}an^c_{in}}}} 
\def\cMancst{{\mathop{\bf\check{M}an^c_{st}}}}
\def\Manab{{\mathop{\bf Man^{ab}}}}
\def\Manac{{\mathop{\bf Man^{ac}}}}
\def\Mancac{{\mathop{\bf Man^{c,ac}}}}
\def\Manacin{{\mathop{\bf Man^{ac}_{in}}}}  
\def\Mancacin{{\mathop{\bf Man^{c,ac}_{in}}}}
\def\Manacst{{\mathop{\bf Man^{ac}_{st}}}} 
\def\Mancacst{{\mathop{\bf Man^{c,ac}_{st}}}}
\def\cManac{{\mathop{\bf\check{M}an^{ac}}}}
\def\cManacin{{\mathop{\bf\check{M}an^{ac}_{in}}}} 
\def\cManacst{{\mathop{\bf\check{M}an^{ac}_{st}}}}
\def\cMancac{{\mathop{\bf\check{M}an^{c,ac}}}}
\def\cMancacin{{\mathop{\bf\check{M}an^{c,ac}_{in}}}}
\def\cMancacst{{\mathop{\bf\check{M}an^{c,ac}_{st}}}}
\def\ul{\underline}
\def\bs{\boldsymbol}
\def\ge{\geqslant}
\def\le{\leqslant\nobreak}
\def\K{{\mathbin{\mathbb K}}}
\def\R{{\mathbin{\mathbb R}}}
\def\Z{{\mathbin{\mathbb Z}}}
\def\Q{{\mathbin{\mathbb Q}}}
\def\N{{\mathbin{\mathbb N}}}
\def\C{{\mathbin{\mathbb C}}}
\def\CP{{\mathbin{\mathbb{CP}}}}
\def\cD{{\mathbin{\cal D}}}
\def\cM{{\mathbin{\cal M}}}
\def\oM{{\mathbin{\smash{\,\,\overline{\!\!\mathcal M\!}\,}}}}
\def\cS{{\mathbin{\cal S}}}
\def\cV{{\mathbin{\cal V}}}
\def\al{\alpha}
\def\be{\beta}
\def\ga{\gamma}
\def\de{\delta}
\def\io{\iota}
\def\ep{\epsilon}
\def\la{\lambda}
\def\ka{\kappa}
\def\th{\theta}
\def\ze{\zeta}
\def\up{\upsilon}
\def\vp{\varphi}
\def\si{\sigma}
\def\om{\omega}
\def\De{\Delta}
\def\La{\Lambda}
\def\Om{\Omega}
\def\Up{\Upsilon}
\def\Ga{\Gamma}
\def\Si{\Sigma}
\def\Th{\Theta}
\def\pd{\partial}
\def\ts{\textstyle}
\def\st{\scriptstyle}
\def\sst{\scriptscriptstyle}
\def\w{\wedge}
\def\sm{\setminus}
\def\bu{\bullet}
\def\op{\oplus}
\def\ot{\otimes}
\def\ov{\overline}
\def\bigop{\bigoplus}
\def\bigot{\bigotimes}
\def\iy{\infty}
\def\es{\emptyset}
\def\ra{\rightarrow}
\def\Longra{\Longrightarrow}
\def\ab{\allowbreak}
\def\longra{\longrightarrow}
\def\hookra{\hookrightarrow}
\def\dashra{\dashrightarrow}
\def\t{\times}
\def\ci{\circ}
\def\ti{\tilde}
\def\d{{\rm d}}
\def\lb{\llbracket}
\def\rb{\rrbracket}
\def\ha{{\ts\frac{1}{2}}}
\def\md#1{\vert #1 \vert}
\def\bmd#1{\big\vert #1 \big\vert}
\def\nm#1{\Vert #1 \Vert}
\def\bnm#1{\big\Vert #1 \big\Vert}
\title{Manifolds with analytic corners}
\author{Dominic Joyce}
\date{}
\maketitle

\begin{abstract} Manifolds with boundary and with corners form categories $\Man\subset\Manb\subset\Manc$. A manifold with corners $X$ has two notions of tangent bundle: the tangent bundle $TX$, and the b-tangent bundle ${}^bTX$. The usual definition of smooth structure uses $TX$, as $f:X\ra\R$ is defined to be smooth if $\nabla^kf$ exists as a continuous section of $\bigot^kT^*X$ for all $k\ge 0$.

We define {\it manifolds with analytic corners}, or {\it manifolds with a-corners}, with a different smooth structure, in which roughly $f:X\ra\R$ is smooth if ${}^b\nabla^kf$ exists as a continuous section of $\bigot^k({}^bT^*X)$ for all $k\ge 0$. These are different from manifolds with corners even when $X=[0,\iy)$, for instance $x^\al:[0,\iy)\ra\R$ is smooth for all real $\al\ge 0$ when $[0,\iy)$ has a-corners. Manifolds with a-boundary and with a-corners form categories $\Man\subset\Manab\subset\Manac$, with well behaved differential geometry. 

Partial differential equations on manifolds with boundary may have boundary conditions of two kinds: (i) `at finite distance', e.g.\ Dirichlet or Neumann conditions, or (ii) `at infinity', prescribing the asymptotic behaviour of the solution. We argue that manifolds with corners should be used for (i), and with a-corners for (ii). We discuss many applications of manifolds with a-corners in boundary problems of type (ii), and to singular p.d.e.\ problems involving `bubbling', `neck-stretching' and `gluing'.
\end{abstract}

\setcounter{tocdepth}{2}
\tableofcontents

\section{Introduction}
\label{ac1}

We are concerned with manifolds with boundary, and with corners, which form categories $\Man\subset\Manb\subset\Manc$ extending the usual category of manifolds $\Man$, where a manifold with corners $X$ is locally modelled on $[0,\iy)^k\t\R^{m-k}$ for $0\le k\le m=\dim X$, with $k=0,1$ if $X$ is a manifold with boundary.

Partial differential equations on manifolds with boundary or corners have been extensively studied, and one usually imposes boundary conditions to make the equations well-behaved. These boundary conditions are of two kinds:
\begin{itemize}
\setlength{\itemsep}{0pt}
\setlength{\parsep}{0pt}
\item[(i)] Boundary conditions `at finite distance', e.g.\ Dirichlet or Neumann boundary conditions, Plateau's problem for minimal surfaces in $\R^3$ bounding a curve $\ga\subset\R^3$, and $J$-holomorphic curves $u:\Si\ra S$ in a symplectic manifold $S$ with $u(\pd\Si)$ in a Lagrangian $L\subset S$, as in~\cite{FOOO1, Seid}.
\item[(ii)] Boundary conditions `at infinity', where we prescribe the asymptotic behaviour of the solution. These include Asymptotically Cylindrical Riemannian manifolds or submanifolds \cite{CHNP1,CHNP2,CMR,HHN,JoSa,Kova,KoNo,Nord}, and Asymptotically Conical (sub)manifolds~\cite{BrSa,CoHe1,CoHe2,EGH,Joyc2,Joyc3,Joyc5, Joyc6,Kron1,Kron2,Lota2,Mars,Paci1}.
\end{itemize}

A manifold with boundary or corners $X$ has two notions of tangent bundle:
\begin{itemize}
\setlength{\itemsep}{0pt}
\setlength{\parsep}{0pt}
\item[(i)] The ordinary tangent bundle $TX$, with the obvious definition; and
\item[(ii)] The b-tangent bundle ${}^bTX$, as in Melrose \cite[\S 2]{Melr1}, \cite[\S 2.2]{Melr2}, \cite[\S I.10]{Melr3}.
\end{itemize}
If $(x_1,\ldots,x_m)\in[0,\iy)^k\t\R^{m-k}$ are local coordinates on an open subset $U\subset X$ then $TX\vert_U$ has basis of sections $\frac{\pd}{\pd x_1},\ldots,\frac{\pd}{\pd x_m}$, and ${}^bTX\vert_U$ basis of sections $x_1\frac{\pd}{\pd x_1},\ldots,x_k\frac{\pd}{\pd x_k},\frac{\pd}{\pd x_{k+1}},\ldots,\frac{\pd}{\pd x_m}$. It is usually appropriate to use $TX$ for boundary conditions of type (i) on $X$, and ${}^bTX$ for type~(ii).

Now in the conventional definition of manifolds with corners $X$, the smooth structure on $X$ is defined using $TX$ in (i). That is, a function $f:X\ra\R$ is smooth if the derivatives $\nabla^kf$ exist as continuous sections of $\bigot^kT^*X$ for all $k\ge 0$. The theme of this paper is that there is a second notion of manifolds with corners defined using ${}^bTX$ in (ii), so that roughly speaking $f:X\ra\R$ is smooth if the b-derivatives ${}^b\nabla^kf$ exist as continuous sections of $\bigot^k({}^bT^*X)$ for all $k\ge 0$. We call this new notion {\it manifolds with analytic corners}, or {\it manifolds with a-corners}. They form categories $\Man\subset\Manab\subset\Manac$.

We will make the case that in a boundary p.d.e.\ problem, one should use manifolds with ordinary corners for type (i) boundary conditions, and manifolds with a-corners for type (ii). Note that manifolds with corners and with a-corners are different even for the simplest example $[0,\iy)$. For instance, when $[0,\iy)$ is a manifold with a-corners, $x^\al:[0,\iy)\ra\R$ is a smooth map for all real~$\al\ge 0$.

Also, even for p.d.e.s on manifolds without boundary, manifolds with a-corners are useful for describing singular limits of solutions, e.g.\ `bubbling', `neck-stretching' and `gluing', and finite time singularities of geometric flows.

As we explain in \S\ref{ac5}--\S\ref{ac6}, there is a huge literature on boundary problems for p.d.e.s which can be rewritten in our `a-corners' language. The author believes that in many cases, manifolds with a-corners will be helpful, and lead to new insights or new methods of proof. There should be a general theory of elliptic equations on a compact manifold with a-corners $X$ (and families of elliptic equations over a proper b-fibration $f:X\ra Y$), with lots of applications. The theory is also particularly relevant to questions about smooth structures on moduli spaces near singular solutions, such as the moduli of stable $J$-holomorphic curves in symplectic geometry studied in~\cite{Fuka,FOOO1,FOOO2,FuOn,Hofe,HWZ1,HWZ2}.

Much of this paper is not wholly new material, but rather old ideas due to many authors, to which I am giving a new interpretation, or looking at from a new angle. I was particularly inspired by the work of Richard Melrose and his collaborators \cite{Loya1,Loya2,LoMe,Melr1,Melr2,Melr3,Melr4,MeMe,MeNi,MePi,Piaz} on the `b-calculus' \cite{Grie}, as explained in~\S\ref{ac55}.

Section \ref{ac2} gives background on conventional manifolds with corners $\Manc$. In \S \ref{ac3} we define manifolds with a-corners $\Manac$, plus some generalizations, and study their categorical properties. Section \ref{ac4} considers the differential geometry of manifolds with corners, including boundaries, corners, b-(co)tangent bundles, b-connections, and b-curvature.
 
Section \ref{ac5} discusses analysis on manifolds with a-corners, with a particular focus on (families of) elliptic p.d.e.s on compact manifolds with a-corners, weighted Sobolev spaces, and Fredholm properties. Section \ref{ac6} briefly describes some areas of current research which could be rewritten in the `a-corners' language, and explains how this could lead to new advances in some cases.
\medskip

\noindent{\it Acknowledgements.} I would like to thank Lino Amorim, Daniel Grieser, and Rafe Mazzeo for helpful conversations. This research was supported by EPSRC Programme Grant EP/I033343/1.

\section{Ordinary manifolds with corners}
\label{ac2}

We discuss the usual category $\Manc$ of manifolds with corners. Some references are Melrose \cite{Melr2,Melr3} and the author~\cite{Joyc10,Joyc13}.

\subsection{The definition of manifolds with corners}
\label{ac21}

\begin{dfn} Use the notation $\R^m_k=[0,\iy)^k\t\R^{m-k}$
for $0\le k\le m$, and write points of $\R^m_k$ as $u=(x_1,\ldots,x_m)$ for $x_1,\ldots,x_k\in[0,\iy)$, $x_{k+1},\ldots,x_m\in\R$. Let $U\subseteq\R^m_k$ and $V\subseteq \R^n_l$ be open, and $f=(f_1,\ldots,f_n):U\ra V$ be a continuous map, so that $f_j=f_j(x_1,\ldots,x_m)$ maps $U\ra[0,\iy)$ for $j=1,\ldots,l$ and $U\ra\R$ for $j=l+1,\ldots,n$. Then we say:
\begin{itemize}
\setlength{\itemsep}{0pt}
\setlength{\parsep}{0pt}
\item[(a)] $f$ is {\it weakly smooth\/} if all derivatives $\frac{\pd^{a_1+\cdots+a_m}}{\pd x_1^{a_1}\cdots\pd x_m^{a_m}}f_j(x_1,\ldots,x_m):U\ra\R$ exist and are continuous in for all $j=1,\ldots,m$ and $a_1,\ldots,a_m\ge 0$, including one-sided derivatives where $x_i=0$ for $i=1,\ldots,k$.
\item[(b)] $f$ is {\it smooth\/} if it is weakly smooth and every $u=(x_1,\ldots,x_m)\in U$ has an open neighbourhood $\ti U$ in $U$ such that for each $j=1,\ldots,l$, either:
\begin{itemize}
\setlength{\itemsep}{0pt}
\setlength{\parsep}{0pt}
\item[(i)] we may uniquely write $f_j(\ti x_1,\ldots,\ti x_m)=F_j(\ti x_1,\ldots,\ti x_m)\cdot\ti x_1^{a_{1,j}}\cdots\ti x_k^{a_{k,j}}$ for all $(\ti x_1,\ldots,\ti x_m)\in\ti U$, where $F_j:\ti U\ra(0,\iy)$ is weakly smooth and $a_{1,j},\ldots,a_{k,j}\in\N=\{0,1,2,\ldots\}$, with $a_{i,j}=0$ if $x_i\ne 0$; or 
\item[(ii)] $f_j\vert_{\smash{\ti U}}=0$.
\end{itemize}
\item[(c)] $f$ is {\it interior\/} if it is smooth, and case (b)(ii) does not occur.
\item[(d)] $f$ is {\it b-normal\/} if it is interior, and in case (b)(i), for each $i=1,\ldots,k$ we have $a_{i,j}>0$ for at most one $j=1,\ldots,l$.
\item[(e)] $f$ is {\it strongly smooth\/} if it is smooth, and in case (b)(i), for each $j=1,\ldots,l$ we have $a_{i,j}=1$ for at most one $i=1,\ldots,k$, and $a_{i,j}=0$ otherwise. 
\item[(f)] $f$ is a {\it diffeomorphism\/} if it is a smooth bijection with smooth inverse.
\end{itemize}
All the classes (a)--(f) include identities and are closed under composition.
\label{ac2def1}
\end{dfn}

\begin{dfn} Let $X$ be a second countable Hausdorff topological space. An {\it $m$-dimensional chart on\/} $X$ is a pair $(U,\phi)$, where
$U\subseteq\R^m_k$ is open for some $0\le k\le m$, and $\phi:U\ra X$ is a
homeomorphism with an open set~$\phi(U)\subseteq X$.

Let $(U,\phi),(V,\psi)$ be $m$-dimensional charts on $X$. We call
$(U,\phi)$ and $(V,\psi)$ {\it compatible\/} if
$\psi^{-1}\ci\phi:\phi^{-1}\bigl(\phi(U)\cap\psi(V)\bigr)\ra
\psi^{-1}\bigl(\phi(U)\cap\psi(V)\bigr)$ is a diffeomorphism between open subsets of $\R^m_k,\R^m_l$, in the sense of Definition~\ref{ac2def1}(f).

An $m$-{\it dimensional atlas\/} for $X$ is a system
$\{(U_a,\phi_a):a\in A\}$ of pairwise compatible $m$-dimensional
charts on $X$ with $X=\bigcup_{a\in A}\phi_a(U_a)$. We call such an
atlas {\it maximal\/} if it is not a proper subset of any other
atlas. Any atlas $\{(U_a,\phi_a):a\in A\}$ is contained in a unique
maximal atlas, the set of all charts $(U,\phi)$ of this type on $X$
which are compatible with $(U_a,\phi_a)$ for all~$a\in A$.

An $m$-{\it dimensional manifold with corners\/} is a second
countable Hausdorff topological space $X$ equipped with a maximal
$m$-dimensional atlas. Usually we refer to $X$ as the manifold,
leaving the atlas implicit, and by a {\it chart\/ $(U,\phi)$ on\/}
$X$, we mean an element of the maximal atlas.

Now let $X,Y$ be manifolds with corners of dimensions $m,n$, and $f:X\ra Y$ a continuous map. We call $f$ {\it weakly smooth}, or {\it smooth}, or {\it interior}, or {\it b-normal}, or {\it strongly smooth}, if whenever $(U,\phi),(V,\psi)$ are charts on $X,Y$ with $U\subseteq\R^m_k$, $V\subseteq\R^n_l$ open, then
\e
\psi^{-1}\ci f\ci\phi:(f\ci\phi)^{-1}(\psi(V))\longra V
\label{ac2eq1}
\e
is weakly smooth, or smooth, or interior, or b-normal, or strongly smooth, respectively, as maps between open subsets of $\R^m_k,\R^n_l$ in the sense of Definition~\ref{ac2def1}.

We write $\Manc$ for the category with objects manifolds with corners $X,Y,$ and morphisms smooth maps $f:X\ra Y$ in the sense above. We will also write $\Mancin,\Mancst$ for the subcategories of $\Manc$ with morphisms interior maps, and strongly smooth maps, respectively.

Write $\cManc$ for the category whose objects are disjoint unions $\coprod_{m=0}^\iy X_m$, where $X_m$ is a manifold with corners of dimension $m$, allowing $X_m=\es$, and whose morphisms are continuous maps $f:\coprod_{m=0}^\iy X_m\ra\coprod_{n=0}^\iy Y_n$, such that
$f\vert_{X_m\cap f^{-1}(Y_n)}:X_m\cap f^{-1}(Y_n)\ra
Y_n$ is a smooth map of manifolds with corners for all $m,n\ge 0$.
Objects of $\cManc$ will be called {\it manifolds with corners of
mixed dimension}. We will also write $\cMancin,\cMancst$ for the subcategories of $\cManc$ with the same objects, and morphisms interior, or strongly smooth, maps.
\label{ac2def2}
\end{dfn}

\begin{rem} There are several non-equivalent definitions of categories of manifolds with corners. Just as objects, without considering morphisms, most authors define manifolds with corners as in Definition \ref{ac2def2}. However, Melrose \cite{Melr1,Melr2,Melr3} imposes an extra condition: in \S\ref{ac22} we will define the boundary $\pd X$ of a manifold with corners $X$, with an immersion $i_X:\pd X\ra X$. Melrose requires that $i_X\vert_C:C\ra X$ should be injective for each connected component $C$ of $\pd X$ (such $X$ are sometimes called {\it manifolds with faces\/}).

There is no general agreement in the literature on how to define smooth maps, or morphisms, of manifolds with corners: 
\begin{itemize}
\setlength{\itemsep}{0pt}
\setlength{\parsep}{0pt}
\item[(i)] Our smooth maps are due to Melrose \cite[\S 1.12]{Melr3}, \cite[\S 1]{Melr1}, who calls them {\it b-maps}. Interior and b-normal maps are also due to Melrose.
\item[(ii)] The author \cite{Joyc10} defined and studied strongly smooth maps above (which were just called `smooth maps' in \cite{Joyc10}). 
\item[(iii)] Monthubert's {\it morphisms of manifolds with corners\/} \cite[Def.~2.8]{Mont} coincide with our strongly smooth b-normal maps. 
\item[(iv)] Most other authors, such as Cerf \cite[\S I.1.2]{Cerf}, define smooth maps of manifolds with corners to be weakly smooth maps, in our notation.
\end{itemize}
\label{ac2rem}
\end{rem}

\subsection{Boundaries and corners of manifolds with corners}
\label{ac22}

The material of this section broadly follows the author \cite{Joyc10,Joyc13}.

\begin{dfn} Let $U\subseteq\R^m_k$ be open. For each
$u=(x_1,\ldots,x_m)$ in $U$, define the {\it depth\/} $\depth_Uu$ of
$u$ in $U$ to be the number of $x_1,\ldots,x_k$ which are zero. That
is, $\depth_Uu$ is the number of boundary faces of $U$
containing~$u$.

Let $X$ be an $m$-manifold with corners. For $x\in X$, choose a
chart $(U,\phi)$ on the manifold $X$ with $\phi(u)=x$ for $u\in U$,
and define the {\it depth\/} $\depth_Xx$ of $x$ in $X$ by
$\depth_Xx=\depth_Uu$. This is independent of the choice of
$(U,\phi)$. For each $l=0,\ldots,m$, define the {\it depth\/ $l$
stratum\/} of $X$ to be
\begin{equation*}
S^l(X)=\bigl\{x\in X:\depth_Xx=l\bigr\}.
\end{equation*}
Then $X=\coprod_{l=0}^mS^l(X)$ and $\overline{S^l(X)}=
\bigcup_{k=l}^m S^k(X)$. The {\it interior\/} of $X$ is
$X^\ci=S^0(X)$. Each $S^l(X)$ has the structure of an
$(m-l)$-manifold without boundary.

\label{ac2def3}
\end{dfn}

\begin{dfn} Let $X$ be an $m$-manifold with corners, $x\in X$, and $k=0,1,\ldots,m$. A {\it local $k$-corner component\/ $\ga$ of\/ $X$ at\/} $x$ is a local choice of connected component of $S^k(X)$ near $x$. That is, for each small open neighbourhood $V$ of $x$ in $X$, $\ga$
gives a choice of connected component $W$ of $V\cap S^k(X)$ with
$x\in\overline W$, and any two such choices $V,W$ and $V',W'$ must
be compatible in that $x\in\overline{(W\cap W')}$. When $k=1$, we call $\ga$ a {\it local boundary component}.

As sets, define the {\it boundary\/} $\pd X$ and {\it k-corners\/} $C_k(X)$ for $k=0,1,\ldots,m$ by
\begin{align*}
\pd X&=\bigl\{(x,\be):\text{$x\in X$, $\be$ is a local boundary
component of $X$ at $x$}\bigr\},\\
C_k(X)&=\bigl\{(x,\ga):\text{$x\in X$, $\ga$ is a local $k$-corner 
component of $X$ at $x$}\bigr\}.
\end{align*}
Define $i_X:\pd X\ra X$ and $\Pi:C_k(X)\ra X$ by $i_X:(x,\be)\mapsto x$, $\Pi:(x,\ga)\mapsto x$.

If $(U,\phi)$ is a chart on $X$ with $U\subseteq\R^m_k$ open, then for each $i=1,\ldots,k$ we can define a chart $(U_i,\phi_i)$ on $\pd X$ by
\begin{align*}
&U_i=\bigl\{(x_1,\ldots,x_{m-1})\in \R^{m-1}_{k-1}:
(x_1,\ldots,x_{i-1},0,x_i,\ldots,x_{m-1})\in
U\subseteq\R^m_k\bigr\},\\
&\phi_i:(x_1,\ldots,x_{m-1})\longmapsto\bigl(\phi
(x_1,\ldots,x_{i-1}, 0,x_i,\ldots,x_{m-1}),\phi_*(\{x_i=0\})\bigr).
\end{align*}
The set of all such charts on $\pd X$ forms an atlas, making $\pd X$ into a manifold with corners of dimension $m-1$, and $i_X:\pd X\ra X$ into a smooth (but not interior) map. Similarly, we make $C_k(X)$ into an $(m-k)$-manifold with corners, and $\Pi:C_k(X)\ra X$ into a smooth map.

We call $X$ a {\it manifold without boundary\/} if $\pd X=\es$, and
a {\it manifold with boundary\/} if $\pd^2X=\es$. We write $\Man$ and $\Manb$ for the full subcategories of $\Manc$ with objects manifolds without boundary, and manifolds with boundary, so that $\Man\subset\Manb\subset\Manc$. This definition of $\Man$ is equivalent to the usual definition of the category of manifolds.
\label{ac2def4}
\end{dfn}

For $X$ a manifold with corners and $k\ge 0$, there are natural identifications 
\ea
\begin{split}
\pd^kX\cong\bigl\{(x&,\be_1,\ldots,\be_k):\text{$x\in X,$
$\be_1,\ldots,\be_k$ are distinct}\\
&\text{local boundary components for $X$ at $x$}\bigr\},
\end{split}
\label{ac2eq2}\\
\begin{split}
C_k(X)\cong\bigl\{(x&,\{\be_1,\ldots,\be_k\}):\text{$x\in X,$
$\be_1,\ldots,\be_k$ are distinct}\\
&\text{local boundary components for $X$ at $x$}\bigr\}.
\end{split}
\label{ac2eq3}
\ea
There is a natural, free, smooth action of the symmetric group $S_k$ on $\pd^kX$, by permutation of $\be_1,\ldots,\be_k$ in \eq{ac2eq2}, and \eq{ac2eq2}--\eq{ac2eq3} give a natural diffeomorphism
\e
C_k(X)\cong\pd^kX/S_k.
\label{ac2eq4}
\e

Corners commute with boundaries: there are natural isomorphisms
\e
\begin{aligned}
\pd C_k(X)\cong C_k&(\pd X)\cong \bigl\{(x,\{\be_1,\ldots,
\be_k\},\be_{k+1}):x\in X,\; \be_1,\ldots,\be_{k+1}\\
&\text{are distinct local boundary components for $X$ at
$x$}\bigr\}.
\end{aligned}
\label{ac2eq5}
\e
For products of manifolds with corners we have natural diffeomorphisms
\ea
\pd(X\t Y)&\cong (\pd X\t Y)\amalg (X\t\pd Y),
\label{ac2eq6}\\
C_k(X\t Y)&\cong \ts\coprod_{i,j\ge 0,\; i+j=k}C_i(X)\t C_j(Y).
\label{ac2eq7}
\ea

\begin{ex} The {\it teardrop\/} $T=\bigl\{(x,y)\in\R^2:x\ge 0$, $y^2\le x^2-x^4\bigr\}$, shown in Figure \ref{ac2fig}, is a manifold with corners of dimension 2. The boundary $\pd T$ is diffeomorphic to $[0,1]$, and so is connected, but $i_T:\pd T\ra T$ is not injective. Thus $T$ is not a manifold with faces, in the sense of Remark \ref{ac2rem}.
\begin{figure}[htb]
\begin{xy}
,(-1.5,0)*{}
,<6cm,-1.5cm>;<6.7cm,-1.5cm>:
,(3,.3)*{x}
,(-1.2,2)*{y}
,(-1.5,0)*{\bullet}
,(-1.5,0); (1.5,0) **\crv{(-.5,1)&(.1,1.4)&(1.5,1.2)}
?(.06)="a"
?(.12)="b"
?(.2)="c"
?(.29)="d"
?(.4)="e"
?(.5)="f"
?(.6)="g"
?(.7)="h"
?(.83)="i"
,(-1.5,0); (1.5,0) **\crv{(-.5,-1)&(.1,-1.4)&(1.5,-1.2)}
?(.06)="j"
?(.12)="k"
?(.2)="l"
?(.29)="m"
?(.4)="n"
?(.5)="o"
?(.6)="p"
?(.7)="q"
?(.83)="r"
,"a";"j"**@{.}
,"b";"k"**@{.}
,"c";"l"**@{.}
,"d";"m"**@{.}
,"e";"n"**@{.}
,"f";"o"**@{.}
,"g";"p"**@{.}
,"h";"q"**@{.}
,"i";"r"**@{.}
\ar (-1.5,0);(3,0)
\ar (-1.5,0);(-3,0)
\ar (-1.5,0);(-1.5,2)
\ar (-1.5,0);(-1.5,-2)
\end{xy}
\caption{The teardrop, a 2-manifold with corners}
\label{ac2fig}
\end{figure}
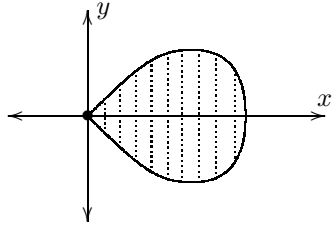

\label{ac2ex1}
\end{ex}

The following lemma is easy to prove from Definition \ref{ac2def1}(b).

\begin{lem} Let\/ $f:X\ra Y$ be a smooth map of manifolds with corners. Then $f$ \begin{bfseries}is compatible with the depth stratifications\end{bfseries} $X=\coprod_{k\ge 0}S^k(X),$ $Y=\coprod_{l\ge 0}S^l(Y)$ in Definition\/ {\rm\ref{ac2def3},} in the sense that if\/ $\es\ne W\subseteq S^k(X)$ is a connected subset for some $k\ge 0,$ then $f(W)\subseteq S^l(Y)$ for some unique $l\ge 0$.
\label{ac2lem}
\end{lem}

It is {\it not\/} true that general smooth $f:X\ra Y$ induce maps $\pd f:\pd X\ra\pd Y$ or $C_k(f):C_k(Y)\ra C_k(Y)$. For example, if $f:X\ra Y$ is the inclusion $[0,\iy)\hookra\R$ then no map $\pd f:\pd X\ra\pd Y$ exists, as $\pd X\ne\es$ and $\pd Y=\es$. So boundaries and $k$-corners do not give functors on $\Manc$. However, if we work in the enlarged category $\cManc$ of Definition \ref{ac2def2} and consider the full corners $C(X)=\coprod_{k\ge 0}C_k(X)$, we can define a functor.

\begin{dfn} Define the {\it corners\/} $C(X)$ of a manifold with corners $X$ by
\begin{align*}
&C(X)=\ts\coprod_{k=0}^{\dim X}C_k(X)\\
&=\bigl\{(x,\ga):\text{$x\in X$, $\ga$ is a local $k$-corner 
component of $X$ at $x$, $k\ge 0$}\bigr\},
\end{align*}
considered as an object of $\cManc$ in Definition \ref{ac2def2}, a manifold with corners of mixed dimension. Define $\Pi:C(X)\ra X$ by $\Pi:(x,\ga)\mapsto x$. This is smooth (i.e. a morphism in $\cManc$) as the maps $\Pi:C_k(X)\ra X$ are smooth for~$k\ge 0$.

Let $f:X\ra Y$ be a smooth map of manifolds with corners, and suppose $\ga$ is a local $k$-corner component of $X$ at $x\in X$. For each sufficiently small open neighbourhood $V$ of $x$ in $X$, $\ga$ gives a choice of connected component $W$ of $V\cap S^k(X)$ with $x\in\overline W$, so by Lemma \ref{ac2lem} $f(W)\subseteq S^l(Y)$ for some $l\ge 0$. As $f$ is continuous, $f(W)$ is connected, and $f(x)\in\ov{f(W)}$. Thus there is a unique local $l$-corner component $f_*(\ga)$ of $Y$ at $f(x)$, such that if $\ti V$ is a sufficiently small open neighbourhood of $f(x)$ in $Y$, then the connected component $\ti W$ of $\ti V\cap S^l(Y)$ given by $f_*(\ga)$ has $f(W)\cap\ti W\ne\es$. This $f_*(\ga)$ is independent of the choice of sufficiently small $V,\ti V$, so is well-defined.

Define a map $C(f):C(X)\ra C(Y)$ by $C(f):(x,\ga)\mapsto (f(x),f_*(\ga))$. Then $C(f)$ is an interior morphism in $\cManc$. If $g:Y\ra Z$ is another smooth map of manifolds with corners then $C(g\ci f)=C(g)\ci C(f):C(X)\ra C(Z)$, so $C:\Manc\ra\cMancin\subset\cManc$ is a functor, which we call the {\it corner functor}.
\label{ac2def5}
\end{dfn}

Equations \eq{ac2eq5} and \eq{ac2eq7} imply that if $X,Y$ are manifolds with corners, we have natural isomorphisms
\ea
\pd C(X)&\cong C(\pd X),
\label{ac2eq8}\\
C(X\t Y)&\cong C(X)\t C(Y).
\label{ac2eq9}
\ea
The corner functor $C$ {\it preserves products and direct products}. That is, if $f:W\ra Y,$ $g:X\ra Y,$ $h:X\ra Z$ are smooth then the following commute
\begin{equation*}
\xymatrix@C=60pt@R=20pt{ *+[r]{C(W\t X)} \ar[d]^\cong \ar[r]_{C(f\t
h)} & *+[l]{C(Y\t Z)} \ar[d]_\cong \\ *+[r]{C(W)\!\t\! C(X)}
\ar[r]^{\raisebox{8pt}{$\st C(f) \t C(h)$}} &
*+[l]{C(Y)\!\t\! C(Z),} }\;
\xymatrix@C=65pt@R=3pt{ & C(Y\t Z) \ar[dd]^\cong \\
C(X) \ar[ur]^(0.4){C((g,h))} \ar[dr]_(0.4){(C(g),C(h))} \\
& C(Y)\!\t\! C(Z), }
\end{equation*}
where the columns are the isomorphisms~\eq{ac2eq9}.

\begin{ex}{\bf(a)} Let $X=[0,\iy)$, $Y=[0,\iy)^2$, and define $f:X\ra Y$
by $f(x)=(x,x)$. We have
\begin{align*}
C_0(X)&\cong[0,\iy), \qquad\quad C_1(X)\cong\{0\}, & C_0(Y)&\cong[0,\iy)^2,\\
C_1(Y)&\cong\bigl(\{0\}\t [0,\iy)\bigr)\amalg \bigl([0,\iy)\t\{0\}\bigr),&
C_2(Y)&\cong\{(0,0)\}.
\end{align*}
Then $C(f)$ maps $C_0(X)\ra C_0(Y)$, $x\mapsto (x,x)$, and
$C_1(X)\ra C_2(Y)$, $0\mapsto(0,0)$.

\smallskip

\noindent{\bf(b)} Let $X=*$, $Y=[0,\iy)$ and define $f:X\ra Y$ by
$f(*)=0$. Then $C_0(X)\cong *$, $C_0(Y)\cong [0,\iy)$,
$C_1(Y)\cong\{0\}$, and $C(f)$ maps $C_0(X)\ra C_1(Y)$, $*\mapsto
0$.
\label{ac2ex2}
\end{ex}

Note that $C(f)$ need not map $C_k(X)\ra C_k(Y)$.

\subsection{Tangent bundles and b-tangent bundles}
\label{ac23}

Manifolds with corners $X$ have two notions of tangent bundle with
functorial properties, the ({\it ordinary\/}) {\it tangent bundle\/}
$TX$, the obvious generalization of tangent bundles of manifolds
without boundary, and the {\it b-tangent bundle\/} ${}^bTX$
introduced by Melrose \cite[\S 2]{Melr1}, \cite[\S 2.2]{Melr2}, \cite[\S I.10]{Melr3}. Taking duals gives two notions of cotangent
bundle $T^*X,{}^bT^*X$. First we discuss vector bundles:

\begin{dfn} Let $X$ be a manifold with corners. A {\it vector
bundle\/ $E\ra X$ of rank\/} $k$ is a manifold with corners $E$ and
a smooth map $\pi:E\ra X$, such that each fibre $E_x:=\pi^{-1}(x)$ for $x\in X$ is given the structure of a real vector space, and $X$ may be covered by open $U\subseteq X$ with diffeomorphisms $\pi^{-1}(U)\cong
U\t\R^k$ identifying $\pi\vert_{\pi^{-1}(U)}:\pi^{-1}(U)\ra U$ with
the projection $U\t\R^k\ra\R^k$, and the vector space structure on
$E_x$ with that on $\{x\}\t\R^k\cong\R^k$, for each $x\in U$. A {\it section\/} of $E$ is a smooth map $s:X\ra E$ with $\pi\ci
s=\id_X$.

We write $\Ga^\iy(E)$ for the vector space of smooth sections of $E$, and $C^\iy(X)$ for the $\R$-algebra of smooth functions $X\ra\R$. Then $\Ga^\iy(E)$ is a $C^\iy(X)$-module.

Morphisms of vector bundles, dual vector bundles, tensor products of
vector bundles, exterior products, and so on, all work as usual.
\label{ac2def6}
\end{dfn}

\begin{dfn} Let $X$ be an $m$-manifold with corners. The {\it tangent bundle\/} $\pi:TX\ra X$ and {\it b-tangent bundle\/} $\pi:{}^bTX\ra X$ are natural rank $m$ vector bundles on $X$, with a vector bundle morphism $I_X:{}^bTX\ra TX$. We may describe $TX,{}^bTX,I_X$ in local coordinates as follows. 

If $(U,\phi)$ is a chart on $X$, with $U\subseteq\R^m_k$ open, and $(x_1,\ldots,x_m)$ are the coordinates on $U$, then over $\phi(U)$, $TX$ is the trivial vector bundle with basis of sections $\frac{\pd}{\pd x_1},\ldots,\frac{\pd}{\pd x_m}$, and ${}^bTX$ is the trivial vector bundle with basis of sections~$x_1\frac{\pd}{\pd x_1},\ldots,x_k\frac{\pd}{\pd x_k},\frac{\pd}{\pd x_{k+1}},\ldots,\frac{\pd}{\pd x_m}$. 

We have corresponding charts $(TU,T\phi)$ on $TX$ and $({}^bTU,{}^bT\phi)$ on ${}^bTX$, where $TU={}^bTU=U\t\R^m\subseteq\R^{2m}_k$, such that $(x_1,\ldots,x_m,q_1,\ldots,q_m)$ in $TU$ represents the vector $q_1\frac{\pd}{\pd x_1}+\cdots+q_m\frac{\pd}{\pd x_m}$ over $\phi(x_1,\ldots,x_m)\in X$, and $(x_1,\ldots,x_m,r_1,\ldots,r_m)$ in ${}^bTU$ represents $r_1x_1\frac{\pd}{\pd x_1}+\cdots+r_kx_k\frac{\pd}{\pd x_k}+r_{k+1}\frac{\pd}{\pd x_{k+1}}+\cdots+r_m\frac{\pd}{\pd x_m}$ over $\phi(x_1,\ldots,x_m)$ in $X$, and $I_X$ maps $(x_1,\ldots,x_m,r_1,\ldots,r_m)$ in ${}^bTU$ to $(x_1,\ldots,x_m,r_1x_1,\ldots,r_kx_k,r_{k+1},\ldots,r_m)$ in $TU$.

Under change of coordinates $(x_1,\ldots,x_m)\rightsquigarrow(\ti x_1,\ldots,\ti x_m)$ from $(U,\phi)$ to $(\ti U,\ti\phi)$, the corresponding change $(x_1,\ldots,x_m,q_1,\ldots,q_m)\rightsquigarrow(\ti x_1,\ldots,\ti q_m)$ from $(TU,T\phi)$ to $(T\ti U,T\ti\phi)$ is determined by $\frac{\pd}{\pd x_i}=\sum_{j=1}^m\frac{\pd\ti x_j}{\pd x_i}(x_1,\ldots,x_m)\cdot\frac{\pd}{\pd\ti x_j}$, so that $\ti q_j=\sum_{i=1}^m\frac{\pd\ti x_j}{\pd x_i}(x_1,\ldots,x_m)q_i$, and similarly for $({}^bTU,{}^bT\phi),({}^bT\ti U,{}^bT\ti\phi)$.

Elements of $\Ga^\iy(TX)$ are called {\it vector fields}, and of $\Ga^\iy({}^bTX)$ are called {\it b-vector fields}. The map $(I_X)_*:\Ga^\iy({}^bTX)\ra \Ga^\iy(TX)$ is injective, and identifies $\Ga^\iy({}^bTX)$ with the vector subspace of $v\in \Ga^\iy(TX)$ such that $v\vert_{S^k(X)}$ is tangent to $S^k(X)$ for all $k=1,\ldots,\dim X$.

The {\it cotangent bundle\/} $T^*X$ and {\it b-cotangent bundle\/} ${}^bT^*X$ are the dual vector bundles of $TX,{}^bTX$. If $(U,\phi)$ is a chart on $X$, with $U\subseteq\R^m_k$ open, and $(x_1,\ldots,x_m)$ are the coordinates on $U$, then $T^*X$ has basis of sections $\d x_1,\ldots,\d x_m$ and ${}^bT^*X$ basis of sections $x_1^{-1}\d x_1,\ldots,x_k^{-1}\d x_k,\d x_{k+1},\ldots,\d x_m$ over $\phi(U)$. We have a vector bundle morphism $I_X^*:T^*X\ra{}^bT^*X$ dual to $I_X$. There is a {\it de Rham differential\/} $\d:C^\iy(X)\ra\Ga^\iy(T^*X)$, which acts in local coordinates $(x_1,\ldots,x_m)$ on $X$ by $\d:c\mapsto \frac{\pd c}{\pd x_1}\d x_1+\cdots+\frac{\pd c}{\pd x_m}\d x_m$. The {\it b-de Rham differential\/} is ${}^b\d=I_X^*\ci\d:C^\iy(X)\ra\Ga^\iy({}^bT^*X)$.

Now suppose $f:X\ra Y$ is a smooth map of manifolds with corners. Then there is a natural smooth map $Tf:TX\ra TY$ so that the following
commutes:
\begin{equation*}
\xymatrix@C=80pt@R=15pt{ *+[r]{TX} \ar[d]^\pi \ar[r]_{Tf} &
*+[l]{TY} \ar[d]_\pi \\ *+[r]{X} \ar[r]^f & *+[l]{Y.} }
\end{equation*}
Let $(U,\phi)$ and $(V,\psi)$ be coordinate charts on $X,Y$ with $U\subseteq\R^m_k$, $V\subseteq\R^n_l$, with coordinates $(x_1,\ldots,x_m)\in U$ and $(y_1,\ldots,y_n)\in V$, and let $(TU,T\phi)$, $(TV,T\psi)$ be the corresponding charts on $TX,TY$, with coordinates $(x_1,\ab\ldots,\ab x_m,\ab q_1,\ab\ldots,\ab q_m)\in TU$ and $(y_1,\ldots,y_n,r_1,\ldots,r_n)\in TV$. Equation \eq{ac2eq1} defines a map $\psi^{-1}\ci f\ci\phi$ between open subsets of $U,V$. Write $\psi^{-1}\ci f\ci\phi=(f_1,\ldots,f_n)$, for $f_j=f_j(x_1,\ldots,x_m)$. Then the corresponding $T\psi^{-1}\ci Tf\ci T\phi$ maps
\begin{align*}
&T\psi^{-1}\ci Tf\ci T\phi:(x_1,\ldots,x_m,q_1,\ldots,q_m)\longmapsto
\bigl(f_1(x_1,\ldots,x_m),\ldots,\\
&f_n(x_1,\ldots,x_m),\ts\sum_{i=1}^m\frac{\pd f_1}{\pd x_i}(x_1,\ldots,x_m)q_i,\ldots,\ts\sum_{i=1}^m\frac{\pd f_n}{\pd x_i}(x_1,\ldots,x_m)q_i\bigr).
\end{align*}
We can also regard $Tf$ as a vector bundle morphism $\d f:TX\ra f^*(TY)$ on $X$, which has dual morphism $\d f:f^*(T^*Y)\ra T^*X$.

If $g:Y\ra Z$ is smooth then $T(g\ci f)=Tg\ci Tf:TX\ra
TZ$, and $T(\id_X)=\id_{TX}:TX\ra TX$. Thus, the assignment
$X\mapsto TX$, $f\mapsto Tf$ is a functor, the {\it tangent
functor\/} $T:\Manc\ra\Manc$. It restricts to $T:\Mancin\ra\Mancin$. 

As in \cite[\S 2]{Melr1}, the analogue of the morphisms $Tf:TX\ra TY$ for b-tangent bundles works only for {\it interior\/} maps $f:X\ra Y$. So let $f:X\ra Y$ be an interior map of manifolds with corners. If $f$ is interior, there is a unique interior map ${}^bTf:{}^bTX\ra{}^bTY$ so that the following commutes:
\begin{equation*}
\xymatrix@C=30pt@R=15pt{ {{}^bTX} \ar[ddr]_(0.6)\pi \ar[dr]^{I_X}
\ar[rrr]_{{}^bTf} &&&
{{}^bTY} \ar[dr]^{I_Y} \ar[ddr]_(0.6)\pi \\
& {TX} \ar[d]^\pi \ar[rrr]^{Tf} &&&
{TY} \ar[d]^\pi \\
& {X} \ar[rrr]^f &&& {Y.\!} }
\end{equation*}
The assignment $X\mapsto {}^bTX$, $f\mapsto {}^bTf$ is a functor, the {\it b-tangent functor\/} ${}^bT:\Mancin\ra\Mancin$. The maps $I_X:{}^bTX\ra TX$ give a natural transformation $I:{}^bT\ra T$ of functors~$\Mancin\ra\Mancin$.

We can also regard ${}^bTf$ as a vector bundle morphism ${}^b\d
f:{}^bTX\ra f^*({}^bTY)$ on $X$. It has a dual morphism ${}^b\d
f:f^*({}^bT^*Y)\ra{}^bT^*X$.

Note that if $f:X\ra Y$ is a smooth map in $\Manc$ then $C(f):C(X)\ra C(Y)$ is interior, so ${}^bTC(f):{}^bTC(X)\ra{}^bTC(Y)$ is well defined, and we can use this as a substitute for ${}^bTf:{}^bTX\ra{}^bTY$ when $f$ is not interior.
\label{ac2def7}
\end{dfn}

The next definition follows Melrose \cite[\S I]{Melr1}, \cite[\S 2]{Melr2}, \cite[\S 2.4]{Melr4}.

\begin{dfn} Let $f:X\ra Y$ be an interior map of manifolds with corners. We call $f$ a {\it b-submersion\/} if ${}^b\d f:{}^bTX\ra f^*({}^bTY)$ is a surjective morphism of vector bundles on $X$. We call $f$ a {\it b-fibration\/} if $f$ is b-normal, in the sense of Definitions \ref{ac2def1}(d) and \ref{ac2def2}, and a b-submersion.
\label{ac2def8}
\end{dfn}

\begin{ex}{\bf(i)} Any projection $\pi_X:X\t Y\ra X$ for $X,Y$ manifolds with corners is b-normal, a b-submersion, and a b-fibration.
\smallskip

\noindent{\bf(ii)} Define $f:[0,\iy)^2\ra[0,\iy)$ by $f(x,y)=xy$. Then ${}^b\d f$ is given by the matrix $\bigl(\begin{smallmatrix} 1 \\ 1 \end{smallmatrix}\bigr)$ with respect to the bases $\bigl(x\frac{\pd}{\pd x},y\frac{\pd}{\pd y}\bigr)$ for ${}^bT\bigl([0,\iy)^2\bigr)$ and $z\frac{\pd}{\pd z}$ for ${}^bT\bigl([0,\iy)\bigr)$, so ${}^b\d f$ is surjective, and $f$ is a b-submersion. Also $f$ is b-normal by Definition \ref{ac2def1}(d) as $a_{1,1}=1$, $a_{2,1}=1$ and other $a_{i,j}=0$ in Definition \ref{ac2def1}(b)(i), so $f$ is a b-fibration.
\smallskip

\noindent{\bf(iii)} Define $g:[0,\iy)\t\R\ra[0,\iy)^2$ by $g(w,x)=(w,we^x)$. Then ${}^b\d g$ is given by the matrix $\bigl(\begin{smallmatrix} 1 & 0 \\ 1 & 1 \end{smallmatrix}\bigr)$ with respect to the bases $\bigl(w\frac{\pd}{\pd w},\frac{\pd}{\pd x}\bigr)$ for ${}^bT\bigl([0,\iy)\t\R\bigr)$ and $\bigl(y\frac{\pd}{\pd y},z\frac{\pd}{\pd z}\bigr)$ for ${}^bT\bigl([0,\iy)^2\bigr)$, so $g$ is a b-submersion. However, $g$ is not b-normal as $a_{1,1}=1$, $a_{1,2}=1$ in Definition \ref{ac2def1}(b)(i), so $g$ is not a b-fibration.
\smallskip

Note that $\d f:TX\ra f^*(TY)$ and $\d g:TX\ra g^*(TY)$ are not surjective in {\bf(ii)},{\bf(iii)}, so $f,g$ are not submersions in the usual sense of differential geometry.

\label{ac2ex3}
\end{ex}

The tangent bundle $TX$ of a manifold with corners $X$ is the obvious generalization of tangent bundles $T\ti X$ of usual manifolds without boundary $\ti X$. For example, if $X\subset\ti X$ with $\ti X$ a manifold without boundary, and the inclusion $X\subset\ti X$ is locally modelled on the inclusion $\R^m_k\subset\R^m$, then $TX=T\ti X\vert_X$. But for many applications, b-tangent bundles are more useful.

When $X$ is a manifold with corners, we will define the {\it normal line bundle\/} $N_{\pd X}$ and {\it b-normal line bundle\/} ${}^bN_{\pd X}$ for the immersion $i_X:\pd X\ra X$.

\begin{dfn} Let $X$ be a manifold with corners. From \S\ref{ac23},
the map $i_X:\pd X\ra X$ induces $Ti_X:T(\pd X)\ra TX$, which we may
regard as an injective morphism $\d i_X:T(\pd X)\ra i_X^*(TX)$ of vector
bundles on $\pd X$. This fits into a natural exact sequence of
vector bundles on~$\pd X$:
\e
\xymatrix@C=25pt{ 0 \ar[r] & T(\pd X) \ar[rr]^{\d i_X} && i_X^*(TX)
\ar[rr]^{\pi_N} && N_{\pd X} \ar[r] & 0, }
\label{ac2eq10}
\e
for $N_{\pd X}\ra\pd X$ the {\it normal line bundle\/} of $\pd X$ in
$X$. While $N_{\pd X}$ is not naturally trivial, it does have an orientation by `outward-pointing' normal vectors.

For b-tangent bundles, the analogue of \eq{ac2eq10} is the exact sequence 
\e
\xymatrix@C=22pt{ 0 \ar[r] & {}^bN_{\pd X} \ar[rr]^(0.4){{}^bi_T} &&
i_X^*({}^bTX) \ar[rr]^{{}^b\pi_T} && {}^bT(\pd X) \ar[r] & 0 }
\label{ac2eq11}
\e
of vector bundles on $\pd X$. Note that \eq{ac2eq11} goes the opposite way to \eq{ac2eq10}. Here ${}^bN_{\pd X}=\pd X\t\R\ra\pd X$ is just the trivial line bundle on $\pd X$, which we call the {\it b-normal bundle\/} of $\pd X$ in~$X$. 

We define the morphisms ${}^bi_T,{}^b\pi_T$ in \eq{ac2eq11} as follows: if $(U,\phi)$ is a chart on $X$, with $U\subseteq\R^m_k$ open, and $(x_1,\ldots,x_m)$ are the coordinates on $U$, and $(x,\be)\in\pd X$ with $x=\phi(\ti x_1,\ldots,\ti x_m)$ with $\ti x_1=0$, and $\be$ is the local boundary component $x_1=0$ of $X$ at $x$, then ${}^bi_T,{}^b\pi_T$ map 
\begin{align*}
{}^bi_T&:\ts\bigl((x,\be),c\bigr)\longmapsto \bigl((x,\be),c\cdot x_1\frac{\pd}{\pd x_1}\bigr),\\
{}^b\pi_T&:\ts\bigl((x,\be),\sum_{i=1}^mc_i\cdot x_i\frac{\pd}{\pd x_i}\bigr)\longmapsto\bigl((x,\be),\sum_{i=2}^mc_i\cdot x_i\frac{\pd}{\pd x_i}\bigr),
\end{align*}
for $c,c_i\in\R$, using $(x_2,\ldots,x_m)\in\R^{m-1}_{k-1}$ as the local coordinates on $\pd X$ near~$x$.

\label{ac2def9}
\end{dfn}
 
\section{Manifolds with analytic corners}
\label{ac3}

We now introduce the category $\Manac$ of {\it manifolds with analytic corners}, or {\it manifolds with a-corners}. Our focus in this section is on definitions and categorical aspects, including functors between our various categories. 

\subsection{\texorpdfstring{A-smooth functions between open sets in $\lb 0,\iy)^k\t\R^{m-k}$}{A-smooth functions between open sets in [0,∞)ᵏ×ℝᵐ⁻ᵏ}}
\label{ac31}

We will use the following notation:

\begin{dfn} As usual, we will write intervals in $\R$ using brackets $[\cdots],(\cdots)$, where $[\,,]$ indicate a closed end of an interval, and $(\,,)$ an open end of an interval. Thus for example we have
\begin{equation*}
(0,1]=\{x\in\R:0<x\le 1\}\quad\text{and}\quad [0,\iy)=\{x\in\R:x\ge 0\}.
\end{equation*}
We will also write intervals in $\R$ using brackets $\lb\cdots\rb$. As sets, these indicate a closed end of an interval, so that $\lb\,,\rb$ just mean the same as $[\,,]$, for instance
\begin{equation*}
\lb 0,1\rb=\{x\in\R:0\le x\le 1\}\quad\text{and}\quad \lb 0,\iy)=\{x\in\R:x\ge 0\}.
\end{equation*}
The difference is that $[\,,]$ mean the interval is considered as a manifold with (ordinary) corners near this end point, and $\lb\,,\rb$ mean the interval is considered as a manifold with analytic corners near this end point, in the sense below.

Sometimes we mix the two kinds, so that $[0,1]\t\lb 0,1\rb$ is a square with coordinates $(x,y)$, such that the sides $x=0$, $x=1$ are ordinary boundaries, and the sides $y=0,$ $y=1$ are analytic boundaries. In Example \ref{ac3ex2} we will also allow $-\iy,\iy$ as end points in the analytic case, so that
\begin{equation*}
\lb-\iy,\iy\rb=\{-\iy\}\amalg\R\amalg\{\iy\},
\end{equation*}
with the obvious topology homeomorphic to $[0,1]$.

For $0\le k,l\le m$ we will write
\begin{equation*}
\R^m_l=[0,\iy)^l\t\R^{m-l}, \qquad \R^{k,m}=\lb 0,\iy)^k\t\R^{m-k},
\end{equation*}
and for $0\le k,l\le k+l\le m$ we will write
\begin{equation*}
\R_l^{k,m}=\lb 0,\iy)^k\t[0,\iy)^l\t\R^{m-k-l},
\end{equation*}
so that $\R^m_l,\R^{k,m},\R_l^{k,m}$ are the local models for manifolds with corners, and with analytic corners, and with both corners and analytic corners, respectively.

\label{ac3def1}
\end{dfn}

The next two definitions are the key to the whole paper.

\begin{dfn} Let $U\subseteq\R^{k,m}$ be open and $f:U\ra\R$ be continuous. Write points of $U$ as $(x_1,\ldots,x_m)$ with $x_1,\ldots,x_k\in\lb 0,\iy)$ and $x_{k+1},\ldots,x_m\in\R$. The {\it b-derivative\/} of $f$ (if it exists) is a map ${}^b\pd f:U\ra\R^m$, written ${}^b\pd f=({}^b\pd_1f,\ldots,{}^b\pd_mf)$ for ${}^b\pd_if:U\ra\R$, where by definition
\e
{}^b\pd_if(x_1,\ldots,x_m)=\begin{cases} 0, & x_i=0, \;\> i=1,\ldots,k, \\
x_i\frac{\pd f}{\pd x_i}(x_1,\ldots,x_m), & x_i>0, \;\> i=1,\ldots,k, \\
\frac{\pd f}{\pd x_i}(x_1,\ldots,x_m), & i=k+1,\ldots,m. 
\end{cases}
\label{ac3eq1}
\e
We say that ${}^b\pd f$ {\it exists\/} if \eq{ac3eq1} is well defined, that is, if $\frac{\pd f}{\pd x_i}$ exists on $U\cap\{x_i>0\}$ if $i=1,\ldots,k$, and $\frac{\pd f}{\pd x_i}$ exists on $U$ if $i=k+1,\ldots,m$.

We can iterate b-derivatives (if they exist), to get maps ${}^b\pd^lf:U\ra\bigot^l\R^m$ for $l=0,1,\ldots,$ by taking b-derivatives of components of ${}^b\pd^jf$ for $j=0,\ldots,l-1$.
\begin{itemize}
\setlength{\itemsep}{0pt}
\setlength{\parsep}{0pt}
\item[(i)] We say that $f$ is {\it roughly differentiable}, or {\it r-differentiable}, if ${}^b\pd f$ exists and is a continuous map ${}^b\pd f:U\ra\R^m$.
\item[(ii)] We say that $f$ is {\it roughly smooth}, or {\it r-smooth}, if ${}^b\pd^lf:U\ra\bigot^l\R^m$ is r-differentiable for all $l=0,1,\ldots.$
\item[(iii)] We say that $f$ is {\it analytically differentiable}, or {\it a-differentiable}, if it is r-differentiable and for any compact subset $S\subseteq U$ and $i=1,\ldots,k$, there exist positive constants $C,\al$ such that
\e
\bmd{{}^b\pd_if(x_1,\ldots,x_m)}\le Cx_i^\al\qquad\text{for all $(x_1,\ldots,x_m)\in S$.}
\label{ac3eq2}
\e

\item[(iv)] We say that $f$ is {\it analytically smooth}, or {\it a-smooth}, or just {\it smooth}, if ${}^b\pd^lf:U\ra\bigot^l\R^m$ is a-differentiable for all $l=0,1,\ldots.$
\end{itemize}

One can show that $f$ is a-smooth if for all $a_1,\ldots,a_m\in\N$ and for any compact subset $S\subseteq U$, there exist positive constants $C,\al$ such that
\e
\begin{split}
\left\vert\frac{\pd^{a_1+\cdots+a_m}}{\pd x_1^{a_1}\cdots\pd x_m^{a_m}}f(x_1,\ldots,x_m)\right\vert\le C \prod_{i=1,\ldots,k:\, a_i>0} x_i^{\al-a_i} \\
\text{for all $(x_1,\ldots,x_m)\in S$ with $x_i>0$ if $i=1,\ldots,k$ with $a_i>0$,}
\end{split}
\label{ac3eq3}
\e
where continuous partial derivatives must exist at the required points.

If $f,g:U\ra\R$ are a-smooth (or r-smooth) and $\la,\mu\in\R$ then $\la f+\mu g
$ and $fg:U\ra\R$ are a-smooth (or r-smooth). Thus, the set $C^\iy(U)$ of a-smooth functions $f:U\ra\R$ is an $\R$-algebra, and in fact a $C^\iy$-{\it ring\/} in the sense of~\cite{Joyc11}.

If $I\subseteq\R$ is an open interval, such as $I=(0,\iy)$, we say that a map $f:U\ra I$ is {\it a-smooth}, or just {\it smooth}, if it is a-smooth as a map $f:U\ra\R$.
\label{ac3def2}
\end{dfn}

\begin{dfn} Let $U\subseteq\R^{k,m}$ and $V\subseteq \R^{l,n}$ be open, and $f=(f_1,\ldots,f_n):U\ra V$ be a continuous map, so that $f_j=f_j(x_1,\ldots,x_m)$ maps $U\ra\lb 0,\iy)$ for $j=1,\ldots,l$ and $U\ra\R$ for $j=l+1,\ldots,n$. Then we say:
\begin{itemize}
\setlength{\itemsep}{0pt}
\setlength{\parsep}{0pt}
\item[(a)] $f$ is {\it r-smooth\/} if $f_j:U\ra\R$ is r-smooth in the sense of Definition \ref{ac3def2} for $j=l+1,\ldots,n$, and every $u=(x_1,\ldots,x_m)\in U$ has an open neighbourhood $\ti U$ in $U$ such that for each $j=1,\ldots,l$, either:
\begin{itemize}
\setlength{\itemsep}{0pt}
\setlength{\parsep}{0pt}
\item[(i)] we may uniquely write $f_j(\ti x_1,\ldots,\ti x_m)=F_j(\ti x_1,\ldots,\ti x_m)\cdot\ti x_1^{a_{1,j}}\cdots\ti x_k^{a_{k,j}}$ for all $(\ti x_1,\ldots,\ti x_m)\in\ti U$, where $F_j:\ti U\ra(0,\iy)$ is r-smooth as in Definition \ref{ac3def2}, and $a_{1,j},\ldots,a_{k,j}\in[0,\iy)$, with $a_{i,j}=0$ if $x_i\ne 0$; or 
\item[(ii)] $f_j\vert_{\smash{\ti U}}=0$.
\end{itemize}
\item[(b)] $f$ is {\it a-smooth}, or just {\it smooth}, if $f_j:U\ra\R$ is a-smooth in the sense of Definition \ref{ac3def2} for $j=l+1,\ldots,n$, and every $u=(x_1,\ldots,x_m)\in U$ has an open neighbourhood $\ti U$ in $U$ such that for each $j=1,\ldots,l$, either:
\begin{itemize}
\setlength{\itemsep}{0pt}
\setlength{\parsep}{0pt}
\item[(i)] we may uniquely write $f_j(\ti x_1,\ldots,\ti x_m)=F_j(\ti x_1,\ldots,\ti x_m)\cdot\ti x_1^{a_{1,j}}\cdots\ti x_k^{a_{k,j}}$ for all $(\ti x_1,\ldots,\ti x_m)\in\ti U$, where $F_j:\ti U\ra(0,\iy)$ is a-smooth as in Definition \ref{ac3def2}, and $a_{1,j},\ldots,a_{k,j}\in[0,\iy)$, with $a_{i,j}=0$ if $x_i\ne 0$; or 
\item[(ii)] $f_j\vert_{\smash{\ti U}}=0$.
\end{itemize}
\item[(c)] $f$ is {\it interior\/} if it is a-smooth, and case (b)(ii) does not occur.
\item[(d)] $f$ is {\it b-normal\/} if it is interior, and in case (b)(i), for each $i=1,\ldots,k$ we have $a_{i,j}>0$ for at most one $j=1,\ldots,l$.
\item[(e)] $f$ is {\it strongly a-smooth}, or just {\it strongly smooth}, if it is a-smooth, and in case (b)(i), for each $j=1,\ldots,l$ we have $a_{i,j}>0$ for at most one $i=1,\ldots,k$. 
\item[(f)] $f$ is an {\it a-diffeomorphism}, or just {\it diffeomorphism}, if it is an a-smooth bijection with a-smooth inverse.
\end{itemize}
All the classes (a)--(f) include identities and are closed under composition.
\label{ac3def3}
\end{dfn}

\begin{rem}{\bf(i)} Readers are advised to study Definitions \ref{ac3def2}, \ref{ac3def3} and Example \ref{ac3ex1} before going on, to get a feeling for r-smooth and a-smooth functions.
\smallskip

\noindent{\bf(ii)} Oversimplifying a little bit, given $f:U\ra\R$ for open $U\subseteq\R^{k,m}$ as in Definition \ref{ac3def2}, the b-derivative ${}^b\pd f:U\ra\R^m$ is
\e
{}^b\pd f=\ts\bigl(x_1\frac{\pd f}{\pd x_1},\ldots,x_k\frac{\pd f}{\pd x_k},\frac{\pd f}{\pd x_{k+1}},\ldots,\frac{\pd f}{\pd x_m}\bigr),
\label{ac3eq4}
\e
a section of the b-cotangent bundle ${}^bT^*U$ as in Definition \ref{ac2def7}, and then $f$ is r-smooth if ${}^b\pd^lf:U\ra\bigot^k\R^m$ exists for all $l\ge 0$, as a section of~$\bigot^l{}^bT^*U$.

Note that for the components $x_i\frac{\pd f}{\pd x_i}$ in \eq{ac3eq4}, we do not need $\frac{\pd f}{\pd x_i}$ to exist at points $(x_1,\ldots,x_m)$ with $x_i=0$, we just write $x_i\frac{\pd f}{\pd x_i}=0$ when $x_i=0$ whether $\frac{\pd f}{\pd x_i}$ exists or not. But we do want the resulting function $x_i\frac{\pd f}{\pd x_i}$ to be continuous, and to have further derivatives.
\smallskip

\noindent{\bf(iii)} From now on, by an abuse of notation we will usually write `$x_i\frac{\pd f}{\pd x_i}$' to mean the function which is 0 when $x_i=0$, and equals $x_i\frac{\pd f}{\pd x_i}$ when $x_i>0$. In particular, {\it we do not require the derivative $\frac{\pd f}{\pd x_i}$ to exist when\/} $x_i=0$, we just {\it define\/} $x_i\frac{\pd f}{\pd x_i}=0$ when $x_i=0$, whether $\frac{\pd f}{\pd x_i}$ exists or not.
\smallskip

\noindent{\bf(iv)} A-smooth functions are r-smooth functions satisfying decay conditions on their derivatives near the boundary hyperplanes $x_i=0$ for $i=1,\ldots,k$.

Think of a-smoothness as a `strong' boundary smoothness condition, and r-smoothness as a `weak' boundary smoothness condition. There are a range of other conditions interpolating between a-smoothness and r-smoothness, which may be useful in some problems. For example, we could replace \eq{ac3eq2} by 
\e
\bmd{{}^b\pd_if(x_1,\ldots,x_m)}\le C\md{\log x_i}^{-\al}\quad\text{when $0<x_i<1$, for $C,\al>0$.}
\label{ac3eq5}
\e

We have chosen to base our theory of manifolds with a-corners primarily upon the `strong' notion of a-smoothness, rather than on the `weak' notion of r-smoothness, or any intermediate condition. This gives our theory better properties --- some results below would be false if we replaced a-smoothness by r-smoothness. But it has the disadvantage of making it less convenient for problems involving r-smooth but not a-smooth functions.
\smallskip

\noindent{\bf(v)} Note the similarity between Definition \ref{ac2def1}(b)--(f) and Definition~\ref{ac3def3}(b)--(f).

One difference is that in the local formula $f_j(\ti x_1,\ldots,\ti x_m)=F_j(\ti x_1,\ldots,\ti x_m)\cdot\ti x_1^{a_{1,j}}\cdots\ti x_k^{a_{k,j}}$ for (a-)smooth maps $U\ra[0,\iy)$ or $U\ra\lb 0,\iy)$, we take $a_{i,j}\in\N$ in Definition \ref{ac2def1}, but $a_{i,j}\in[0,\iy)$ in Definition~\ref{ac3def3}.
\label{ac3rem1}
\end{rem}

\begin{ex}{\bf(a)} Define $f:\lb 0,\iy)\ra\R$ by $f(x)=x^\al$, for $\al>0$ in $\R$. Then $f$ is a-smooth, since \eq{ac3eq3} for any $S\subseteq\lb 0,\iy)$ reduces to
\begin{equation*}
\ts\bmd{\frac{\d^l}{\d x^l}f(x)}\le C x^{\al-l}, 
\end{equation*}
which holds with $C=\al\md{\al-1}\cdots\md{\al-l+1}$ for all $l\ge 0$, and the given $\al$.
\smallskip

\noindent{\bf(b)} Define $f:\lb 0,\iy)\ra\R$ by $f(x)=x^\al\cdot\sin(\log x)$ for $x>0$ and $f(0)=0$, for $\al>0$ in $\R$. Then $x^l\frac{\d^l}{\d x^l}f(x)=O(x^\al)$ as $x\ra 0$ for all $l\ge 0$, so $f$ is a-smooth.

In both {\bf(a)},{\bf(b)} $f$ is not smooth in the usual sense, and these examples suggest that Taylor series may not work well for a-smooth functions.
\smallskip

\noindent{\bf(c)} Define $f:\lb 0,1)\ra\R$ by $f(x)=(\log x)^{-1}$ for $0<x<1$ and $f(0)=0$. Then
\begin{equation*}
{}^b\pd^lf(x)=\ts\bigl(x\frac{\d}{\d x}\bigr)^l f(x)=(-1)^ll!(\log x)^{-l-1}\;\>\text{for $x\in (0,1)$, and ${}^b\pd^lf(0)=0$.}
\end{equation*}
Thus ${}^b\pd^lf:\lb 0,1)\ra\R$ is continuous for all $l\ge 0$, and $f$ is r-smooth, and satisfies the intermediate decay condition \eq{ac3eq5}. However, $f$ is not a-smooth, since we do not have ${}^b\pd^lf=O(x^\al)$ as $x\ra 0$ for $l>0$ and any $\al>0$.
\smallskip

\noindent{\bf(d)} Define $f:\lb 0,1)\ra\R$ by $f(x)=(\log(-\log x))^{-1}$ for $0<x<1$, and $f(0)=0$. As in {\bf(c)}, $f$ is r-smooth and satisfies \eq{ac3eq5}, but is not a-smooth.
 \smallskip

\noindent{\bf(e)} Let $U\subseteq \R_k^m=\R^{k,m}$ and $V\subseteq \R_l^n=\R^{l,m}$ be open, and $f:U\ra V$ be smooth (or strongly smooth, or interior) in the sense of Definition \ref{ac2def1}. Then $f$ is also a-smooth (or strongly a-smooth, or interior, respectively) in the sense of Definition \ref{ac3def3}. The converse does not hold, as {\bf(a)},{\bf(b)} show.
\label{ac3ex1}
\end{ex}

\subsection{The definition of manifolds with a-corners}
\label{ac32}

We can now define manifolds with a-corners.

\begin{dfn} Let $X$ be a second countable Hausdorff topological space. An {\it $m$-dimensional a-chart on\/} $X$ is a pair $(U,\phi)$, where
$U\subseteq\R^{k,m}$ is open for some $0\le k\le m$, and $\phi:U\ra X$ is a homeomorphism with an open set~$\phi(U)\subseteq X$.

Let $(U,\phi),(V,\psi)$ be $m$-dimensional a-charts on $X$. We call
$(U,\phi)$ and $(V,\psi)$ {\it compatible\/} if
$\psi^{-1}\ci\phi:\phi^{-1}\bigl(\phi(U)\cap\psi(V)\bigr)\ra
\psi^{-1}\bigl(\phi(U)\cap\psi(V)\bigr)$ is an a-diffeomorphism between open subsets of $\R^{k,m},\R^{l,m}$, as in Definition~\ref{ac3def3}(f).

An $m$-{\it dimensional a-atlas\/} for $X$ is a system
$\{(U_a,\phi_a):a\in A\}$ of pairwise compatible $m$-dimensional
a-charts on $X$ with $X=\bigcup_{a\in A}\phi_a(U_a)$. We call such an
a-atlas {\it maximal\/} if it is not a proper subset of any other
a-atlas. Any a-atlas $\{(U_a,\phi_a):a\in A\}$ is contained in a unique
maximal a-atlas, the set of all a-charts $(U,\phi)$ on $X$ which are compatible with $(U_a,\phi_a)$ for all~$a\in A$.

An $m$-{\it dimensional manifold with analytic corners}, or {\it with a-corners}, is a second countable Hausdorff topological space $X$ equipped with a maximal $m$-dimensional a-atlas. Usually we refer to $X$ as the manifold with a-corners, leaving the atlas implicit, and by an {\it a-chart\/ $(U,\phi)$ on\/} $X$, we mean an element of the maximal atlas.

Now let $X,Y$ be manifolds with a-corners of dimensions $m,n$, and $f:X\ra Y$ a continuous map. We call $f$ {\it r-smooth}, or {\it a-smooth\/} (or just {\it smooth\/}), or {\it interior}, or {\it b-normal}, or {\it strongly a-smooth\/} (or just {\it strongly smooth\/}), if whenever $(U,\phi),(V,\psi)$ are a-charts on $X,Y$ with $U\subseteq\R^{k,m}$, $V\subseteq\R^{l,n}$ open, then
\e
\psi^{-1}\ci f\ci\phi:(f\ci\phi)^{-1}(\psi(V))\longra V
\label{ac3eq6}
\e
is r-smooth, or a-smooth, or interior, or b-normal, or strongly a-smooth, respectively, as maps between open subsets of $\R^{k,m},\R^{l,n}$ in the sense of Definition~\ref{ac3def3}. 

We call $f:X\ra Y$ an {\it a-diffeomorphism}, or just {\it diffeomorphism}, if $f$ is a bijection and $f:X\ra Y$, $f^{-1}:Y\ra X$ are a-smooth. 

These five classes of (a) r-smooth, (b) a-smooth, (c) interior, (d) b-normal, and (e) strongly a-smooth maps of manifolds with a-corners, all contain identities and are closed under composition, so each makes manifolds with a-corners into a category. Here r-smooth maps may not be a-smooth.

We write $\Manac$ for the category with objects manifolds with a-corners $X,Y,$ and morphisms a-smooth maps $f:X\ra Y$ in the sense above. We write $\Manacin,\Manacst$ for the subcategories of $\Manac$ with morphisms interior maps, and strongly a-smooth maps, respectively.

Write $\cManac$ for the category whose objects are disjoint unions $\coprod_{m=0}^\iy X_m$, where $X_m$ is a manifold with a-corners of dimension $m$, allowing $X_m=\es$, and whose morphisms are continuous maps $f:\coprod_{m=0}^\iy X_m\ra\coprod_{n=0}^\iy Y_n$, such that
$f\vert_{X_m\cap f^{-1}(Y_n)}:X_m\cap f^{-1}(Y_n)\ra
Y_n$ is an a-smooth map of manifolds with a-corners for all $m,n\ge 0$.
Objects of $\cManac$ will be called {\it manifolds with a-corners of
mixed dimension}. We regard $\Manac$ as a full subcategory of~$\cManac$.

Alternatively, we can regard $\cManac$ as the category defined as for $\Manac$ above, except that in defining a-atlases $\{(U_a,\phi_a):a\in A\}$ on $X$, we do not require a-charts $(U_a,\phi_a)$ in the atlas to have the same dimension~$\dim U_a=m$.

We will also write $\cManacin,\cManacst$ for the subcategories of $\cManc$ with the same objects, and morphisms interior, or strongly a-smooth, respectively.
\label{ac3def4}
\end{dfn}

\begin{dfn} We will define a functor $F_\Manc^\Manac:\Manc\ra\Manac$, for $\Manc$ as in \S\ref{ac21} and $\Manac$ as above. Let $X$ be a manifold with corners of dimension $m$, with maximal atlas $\{(U_a,\phi_a):a\in A\}$, so that $U_a\subseteq\R^m_k$ is open for $a\in A$. Identifying $\R^m_k=[0,\iy)^k\t\R^{m-k}=\R^{k,m}$, so that $U_a\subseteq\R^{k,m}$ is open, we see that $(U_a,\phi_a)$ is an $m$-dimensional a-chart on $X$ in the sense of Definition~\ref{ac3def4}. 

If $a,b\in A$ then $\phi_b^{-1}\ci\phi_a:\phi_a^{-1}(\phi_b(U_b))\ra\phi_b^{-1}(\phi_a(U_a))$ is a diffeomorphism between open subsets of $\R_k^m,\R_l^m$, as $(U_a,\phi_a),(U_b,\phi_b)$ are compatible charts. Thus $\phi_b^{-1}\ci\phi_a$ is also an a-diffeomorphism between open subsets of $\R^{k,m},\R^{l,m}$ by Example \ref{ac3ex1}(e). So $(U_a,\phi_a),(U_b,\phi_b)$ are compatible a-charts. Hence $\{(U_a,\phi_a):a\in A\}$ is an a-atlas on $X$ in the sense of Definition \ref{ac3def4}. In general it will {\it not\/} be a maximal a-atlas (though it is a maximal atlas). 

Write $\{(\ti U_a,\ti\phi_a):a\in\ti A\}$ for the unique maximal a-atlas on $X$ containing $\{(U_a,\phi_a):a\in A\}$. Then $\bigl(X,\{(\ti U_a,\ti\phi_a):a\in\ti A\}\bigr)$ is a manifold with a-corners. Set~$F_\Manc^\Manac\bigl(X,\{(U_a,\phi_a):a\in A\}\bigr)=\bigl(X,\{(\ti U_a,\ti\phi_a):a\in\ti A\}\bigr)$.

Let $X,Y$ be manifolds with corners, with maximal atlases $\{(U_a,\phi_a):a\in A\}$, $\{(V_b,\psi_b):b\in B\}$, and write $\{(\ti U_a,\ti\phi_a):a\in\ti A\}$, $\{(\ti V_b,\ti\psi_b):b\in\ti B\}$ for the corresponding maximal a-atlases. Let $f:X\ra Y$ be smooth (i.e.\ $f$ is a morphism in $\Manc$). Then for all $a\in A$ and $b\in B$, 
\begin{equation*}
\psi_b^{-1}\ci f\ci\phi_a:(f\ci\phi_a)^{-1}(\psi_b(V_b))\longra V_b
\end{equation*}
is a smooth map between open subsets of $\R^m_k,\R^n_l$. Hence it is also an a-smooth map between open subsets of $\R^{k,m},\R^{l,n}$ by Example \ref{ac3ex1}(e). As this holds for sets of a-charts $(U_a,\phi_a),(V_b,\psi_b)$ covering $X$ and $Y$, we see that $f:X\ra Y$ is an a-smooth map between the manifolds with a-corners $\bigl(X,\{(\ti U_a,\ti\phi_a):a\in \ti A\}\bigr)$, $\bigl(Y,\{(\ti V_b,\ti\psi_b):b\in\ti B\}\bigr)$, that is, $f$ is a morphism in $\Manac$. Set~$F_\Manc^\Manac(f)=f$. 

Clearly, this defines a functor $F_\Manc^\Manac:\Manc\ra\Manac$, which is faithful (injective on morphisms). Example \ref{ac3ex1}(a),(b) show that $F_\Manc^\Manac$ is not full (surjective on morphisms). Comparing Definitions \ref{ac2def1}(c),(e) and \ref{ac3def3}(c),(e) we see that $F_\Manc^\Manac$ maps $\Mancin\ra\Manacin$ and $\Mancst\ra\Manacst$. It extends to $F_\cManc^\cManac:\cManc\ra\cManac$ in the obvious way.
\label{ac3def5}
\end{dfn}

\begin{ex} Given any interval $I$ in $\R$, such as $[a,b]$, $(a,b)$, $[0,\iy),\ldots,$ we regard $I$ as a manifold with corners in the usual way, and make $I$ into a manifold with a-corners by applying $F_\Manc^\Manac$. When we do this, as in Definition \ref{ac3def1} we replace brackets `$[\,,]$' by `$\lb\,,\rb$', so that $\lb a,b\rb=F_\Manc^\Manac([a,b])$, for instance.

We extend this notation to allow $-\iy,\iy$ as closed end points of intervals for manifolds with a-corners, so that for example $\lb -\iy,\iy\rb$ is a manifold with a-corners, with topological space $\{-\iy\}\amalg\R\amalg\{\iy\}$, with the obvious topology homeomorphic to $[0,1]$. We do this such that $\bigl(\lb 0,\ep),\phi_\pm\bigr)$ are a-charts on $\lb -\iy,\iy\rb$ near $\pm\iy$ for $\ep>0$, where $\phi_\pm:\lb 0,\ep)\ra\lb -\iy,\iy\rb$ are defined by
\e
\phi_+(x)=\begin{cases} -\log x, & x>0, \\ \iy,& x=0, \end{cases}\qquad
\phi_-(x)=\begin{cases} \log x, & x>0, \\ -\iy,& x=0. \end{cases}
\label{ac3eq7}
\e

Note that the a-smooth structure on $\lb -\iy,\iy\rb$ near $\pm\iy$ is determined by the asymptotic behaviour as $x\ra 0$ of the functions $\pm\log x$ in \eq{ac3eq7}, so for instance $\phi_\pm(x)=\pm x^{-1}$ for $x>0$ would give a different a-smooth structure. We chose \eq{ac3eq7} so that in the notation of \S\ref{ac42}, the vector field $\frac{\d}{\d x}$ on $\R$ extends to an a-smooth, nonvanishing section of the b-tangent bundle ${}^bT\lb -\iy,\iy\rb$.
\label{ac3ex2}
\end{ex}

\begin{ex} Consider $\lb 0,\iy)\t\R$ as a manifold with a-corners, with coordinates $(x,y)$. Fix $\al>0$, and let $\Z$ act freely on $\lb 0,\iy)\t\R$ by
\begin{equation*}
n:(x,y)\longmapsto (x^{\al^n},y+n).
\end{equation*}
This is an a-diffeomorphism for each $n\in\Z$. Define $X=(\lb 0,\iy)\t\R)/\Z$. Then $X$ has a natural manifold with a-corners structure, such that the projection $\lb 0,\iy)\t\R\ra X$ is a-smooth, and a local a-diffeomorphism. We will show in Example \ref{ac4ex3} that if $\al\ne 1$ then $X\not\cong F_\Manc^\Manac(Y)$ in $\Manac$ for any manifold with corners $Y$. Thus $F_\Manc^\Manac:\Manc\ra\Manac$ is not essentially surjective.
\label{ac3ex3}
\end{ex}

\begin{ex} Products work as usual in $\Manac$. For manifolds with a-corners $X,Y$, there is a unique manifold with a-corner structure on $X\t Y$ with dimension $\dim X+\dim Y$, such that if $(U,\phi),(V,\psi)$ are a-charts on $X,Y$, then $(U\t V,\phi\t\psi)$ is an a-chart on $X\t Y$. If $g:X\ra Y$, $h:X\ra Z$ are a-smooth maps then the {\it direct product\/} $(g,h):X\ra Y\t Z$, $(g,h)(x)\mapsto(g(x),h(x))$, is a-smooth. If $g:W\ra Y$, $h:X\ra Z$ are a-smooth then the {\it product\/} $g\t h:W\t X\ra Y\t Z$, $(g\t h)(w,x)\mapsto(g(w),h(x))$, is a-smooth.
\label{ac3ex4}
\end{ex}

\subsection{\texorpdfstring{`Gluing profiles', and a functor $\Manacst\ra\Mancst$}{\textquoteleft Gluing profiles\textquoteright, and a functor Manᵃᶜ → Manᶜ}}
\label{ac33}

The material of this section is inspired by the method of `gluing profiles' in Hofer, Wysocki and Zehnder's theory of polyfolds \cite[\S 2]{Hofe}, \cite[\S 2.1]{HWZ1}. A parallel idea occurs in Fukaya, Oh, Ohta and Ono's theory of Kuranishi spaces \cite[\S A1.4, p.~777]{FOOO1}, \cite{FOOO2}. A very similar notion also occurs in Hassell, Mazzeo and Melrose \cite[\S 2.5]{HMM}, under the name `total logarithmic blowup', where it is regarded as way of modifying the smooth structure of a manifold with corners. We explain the connection to \cite{FOOO1,FOOO2,Hofe,HWZ1} in Remark~\ref{ac3rem2}(c).

Definition \ref{ac3def5} defined a functor $F_\Manc^\Manac:\Manc\ra\Manac$. We would like to go the other way, and define a functor $\Manac\ra\Manc$ turning manifolds with a-corners into manifolds with corners. It turns out that this works only for {\it strongly a-smooth\/} maps, giving a functor~$F_\Manacst^\Mancst:\Manacst\ra\Mancst$.

\begin{dfn} Define a homeomorphism $\vp:[0,\iy)\ra\lb 0,\iy)$ by
\e
\vp(x)=\begin{cases} 0, & x=0, \\ e^{x-x^{-1}}, & x>0. \end{cases}
\label{ac3eq8}
\e
The inverse map is $\vp^{-1}:\lb 0,\iy)\ra[0,\iy)$, given by
\e
\vp^{-1}(x)=\begin{cases} 0, & x=0, \\ \ha\bigl(\log x+\sqrt{(\log x)^2+4}\,\bigr), & x>0. \end{cases}
\label{ac3eq9}
\e
Here we think of the domain $\lb 0,\iy)$ of $\vp^{-1}$ as a manifold with a-corners, and the target $[0,\iy)$ as a manifold with corners, and $\vp^{-1}$ as our recipe for how to convert the manifold with a-corners $\lb 0,\iy)$ into a manifold with corners $[0,\iy)$. Note that for small $x>0$ we have $\vp(x)\approx e^{-1/x}$ and~$\vp^{-1}(x)\approx -(\log x)^{-1}$.

Suppose $X$ is a manifold with a-corners of dimension $m$, with maximal a-atlas $\bigl\{(U_a,\phi_a):a\in A\bigr\}$. Let $a\in A$, so that $U_a\subseteq\R^{k,m}$ is open. Define an open subset $\ti U_a\subseteq\R^m_k$ and a continuous map $\ti\phi_a:\ti U_a\ra X$ by
\begin{align*}
\ti U_a&=\!\bigl\{(\vp^{-1}(x_1),\ldots,\vp^{-1}(x_k),x_{k+1},\ldots,x_m)\!\in\!\R^m_k:(x_1,\ldots,x_m)\!\in\! U_a\bigr\},\\
\ti\phi_a&:(\ti x_1,\ldots,\ti x_m)\longmapsto \phi_a\bigl(\vp(\ti x_1),\ldots, \vp(\ti x_k),\ti x_{k+1},\ldots,\ti x_m\bigr).
\end{align*}
Then we have a commutative diagram of continuous maps
\begin{equation*}
\xymatrix@C=100pt@R=15pt{ \R_k^m \ar[d]^{\vp^k\t\id_\R^{m-k}}_\cong & \ti U_a \ar[d]_{\vp^k\t\id_\R^{m-k}}^\cong \ar[l]^\subset \ar[dr]^{\ti\phi_a} \\
\R^{k,m} & U_a \ar[l]_\subset \ar[r]^(0.4){\phi_a} & X, } 
\end{equation*}
with vertical maps homeomorphisms. Since $(U_a,\phi_a)$ is an a-chart on $X$, $(\ti U_a,\ti\phi_a)$ is a chart on $X$, in the sense of Definition~\ref{ac2def2}.

Let $a,b\in A$, with $U_a\subseteq\R^{k,m}$, $U_b\subseteq\R^{l,m}$ open, so that $(U_a,\phi_a),(U_b,\phi_b)$ are compatible a-charts on $X$. We claim that $(\ti U_a,\ti\phi_a),(\ti U_b,\ti\phi_b)$ are compatible charts on $X$. To prove this, set $\dot U_a=\phi_a^{-1}(\phi_b(U_b))$, $\dot U_b=\phi_b^{-1}(\phi_a(U_a))$, so that $\dot U_a\subseteq\R^{k,m}$, $\dot U_b\subseteq\R^{l,m}$ are open, and $\phi_b^{-1}\ci\phi_a:\dot U_a\ra\dot U_b$ is an a-diffeomorphism. Write $(x_1,\ldots,x_m)$ for the coordinates on $\dot U_a$, and $\phi_b^{-1}\ci\phi_a=(f_1,\ldots,f_m)$ for $f_j=f_j(x_1,\ldots,x_m)$, where $f_j:\dot U_a\ra\lb 0,\iy)$ for $j=1,\ldots,l$ and $f_j:\dot U_a\ra\R$ for $j=l+1,\ldots,m$ are a-smooth (and interior). By Definition \ref{ac3def3}(b), for $j=1,\ldots,l$ we may write
\begin{equation*}
f_j(x_1,\ldots,x_m)=F_j(x_1,\ldots,x_m)\cdot x_1^{a_{1,j}}\cdots x_k^{a_{k,j}}
\end{equation*}
locally in $\dot U_a$, where $F_j:\dot U_a\ra(0,\iy)$ is r-smooth, and $a_{1,j},\ldots,a_{k,j}\in[0,\iy)$. Since $(f_1,\ldots,f_m)$ is an a-diffeomorphism, and so strongly a-smooth, by Definition \ref{ac3def3}(b) for each $j=1,\ldots,l$ we have $a_{i,j}>0$ for at most one $i=1,\ldots,k$. 

Similarly, set $\dot{\ti U}_a=\ti\phi_a^{-1}(\ti\phi_b(\ti U_b))$, $\dot{\ti U}_b=\ti\phi_b^{-1}(\ti\phi_a(\ti U_a))$, write $(\ti x_1,\ldots,\ti x_m)$ for the coordinates on $\dot{\ti U}_a$, and $\ti\phi_b^{-1}\ci\ti\phi_a=(\ti f_1,\ldots,\ti f_m)$ for $\ti f_j=\ti f_j(\ti x_1,\ldots,\ti x_m)$, where $\ti f_j:\dot{\ti U}_a\ra[0,\iy)$ for $j=1,\ldots,l$ and $\ti f_j:\dot{\ti U}_a\ra\R$ for $j=l+1,\ldots,m$. From the definitions we find that for $j=1,\ldots,l$ we have
\ea
\begin{split}
\ti f_j(\ti x_1,\ldots,\ti x_m)&=\vp^{-1}\bigl[F_j(\vp(\ti x_1),\ldots,\vp(\ti x_k),\ti x_{k_1},\ldots,\ti x_m)\\
&\qquad {}\cdot \vp(\ti x_1)^{a_{1,j}}\cdots \vp(\ti x_k)^{a_{k,j}}\bigr],\quad 1\le j\le l,
\end{split}
\label{ac3eq10}\\
\ti f_j(\ti x_1,\ldots,\ti x_m)&=f_j(\vp(\ti x_1),\ldots,\vp(\ti x_k),\ti x_{k_1},\ldots,\ti x_m),\quad l<j\le m.
\label{ac3eq11}
\ea

We now claim that $\ti f_j:\dot{\ti U}_a\ra[0,\iy),\R$ in \eq{ac3eq10}--\eq{ac3eq11} are smooth maps of manifolds with corners in the sense of \S\ref{ac21}. To see this for \eq{ac3eq11}, note that if $i=1,\ldots,k$ and $j=l+1,\ldots,m$ then
\e
\begin{split}
\ts\frac{\pd\ti f_j}{\pd\ti x_i}(\ti x_1&,\ldots,\ti x_m)=\ts
\frac{\pd f_j}{\pd x_i}(\vp(\ti x_1),\ldots,\vp(\ti x_k),\ti x_{k_1},\ldots,\ti x_m)\cdot \frac{\pd\vp(\ti x_i)}{\pd\ti x_i}\\
&=\ts\bigl(x_i \frac{\pd f_j}{\pd x_i}\bigr)(\vp(\ti x_1),\ldots,\vp(\ti x_k),\ti x_{k_1},\ldots,\ti x_m)\cdot (1+\ti x_i^{-2}),
\end{split}
\label{ac3eq12}
\e
using $x_i=\vp(\ti x_i)$ and $\frac{\pd\vp(\ti x_i)}{\pd\ti x_i}=(1+\ti x_i^{-2})\cdot \vp(\ti x_i)$. By definition of a-smoothness in \S\ref{ac31}, locally in $\dot U_a$ we have
\begin{equation*}
\ts\bigl\vert x_i \frac{\pd f_j}{\pd x_i}(x_1,\ldots,x_m)\bigr\vert\le Cx_i^\al
\end{equation*}
for some $C,\al>0$. Hence by \eq{ac3eq12}, locally in $\dot{\ti U}_a$ we have
\e
\ts\bigl\vert\frac{\pd\ti f_j}{\pd\ti x_i}(\ti x_1,\ldots,\ti x_m)\bigr\vert\le
C\vp(\ti x_i)^\al\cdot (1+\ti x_i^{-2}).
\label{ac3eq13}
\e
The right hand side of \eq{ac3eq13} tends to 0 as $\ti x_i\ra 0_+$. It is obvious from the definitions that $\frac{\pd\ti f_j}{\pd\ti x_i}$ exists in $\dot{\ti U}_a$ where $\ti x_i>0$, but \eq{ac3eq13} implies that $\frac{\pd\ti f_j}{\pd\ti x_i}$ also extends continuously over $\ti x_i=0$ in $\dot{\ti U}_a$, with value~0.

A generalization of this argument, using a-smoothness of $f_j$, shows that all multiderivatives $\frac{\pd^{a_1+\cdots+a_m}}{\pd\ti x_1^{a_1}\cdots\pd\ti x_m^{a_m}}\ti f_j(\ti x_1,\ldots,\ti x_m)$ for $a_1,\ldots,a_m\ge 0$ exist and are continuous in $\dot{\ti U}_a$, and take the value 0 if $i=1,\ldots,k$ with $\ti x_i=0$ and $a_i>0$. Hence $\ti f_j:\dot{\ti U}_a\ra\R$ in \eq{ac3eq11} for $j=l+1,\ldots,m$ is smooth in the sense of~\S\ref{ac21}.

For \eq{ac3eq10}, fix $j=1,\ldots,l$. Then by strong a-smoothness we have $a_{i,j}>0$ for at most one $i=1,\ldots,k$. If there are no such $i$, the argument above easily shows that $\ti f_j:\dot{\ti U}_a\ra(0,\iy)$ is smooth. If there is one such $i$ then we may write
\begin{equation*}
\ti f_j(\ti x_1,\ldots,\ti x_m)=\Psi\bigl(F_j(\vp(\ti x_1),\ldots,\vp(\ti x_k),\ti x_{k_1},\ldots,\ti x_m),\ti x_i\bigr),
\end{equation*}
where $\Psi:(0,\iy)\t[0,\iy)\ra[0,\iy)$ is given by
\begin{equation*}
\Psi(s,t)=\vp^{-1}\bigl(s\cdot \vp(t)^{a_{i,j}}\bigr).
\end{equation*}
One can show by direct calculation that $\Psi$ is a strongly smooth map of manifolds with corners, with $\Psi(s,t)\approx \smash{a_{i,j}^{-1}t}$ for small $t\ge 0$. The argument above shows that $F_j(\vp(\ti x_1),\ldots,\vp(\ti x_k),\ti x_{k_1},\ldots,\ti x_m)$ is smooth on $\dot{\ti U}_a$, as $F_j$ is a-smooth on $\dot U_a$, and $\ti x_i:\dot{\ti U}_a\ra[0,\iy)$ is clearly strongly smooth. Hence $\ti f_j:\dot{\ti U}_a\ra [0,\iy)$ in \eq{ac3eq10} for $j=1,\ldots,l$ is strongly smooth in the sense of~\S\ref{ac21}.
 
This proves that the transition function $\ti\phi_b^{-1}\ci\ti\phi_a:\dot{\ti U}_a\ra\dot{\ti U}_b$ is strongly smooth. Similarly $\ti\phi_a^{-1}\ci\ti\phi_b:\dot{\ti U}_b\ra\dot{\ti U}_a$ is strongly smooth, so $(\ti U_a,\ti\phi_a),(\ti U_b,\ti\phi_b)$ are compatible charts on $X$. As this holds for all $a,b\in A$, $\bigl\{(\ti U_a,\ti\phi_a):a\in A\bigr\}$ is an atlas on $X$, which extends to a unique maximal atlas, making $X$ into an $m$-dimensional manifold with corners, which we write as $\ti X$. Set~$F_\Manacst^\Mancst(X)=\ti X$.

Now suppose $f:X\ra Y$ is a strongly a-smooth map of manifolds with a-corners, and write $\ti X=F_\Manacst^\Mancst(X)$, $\ti Y=F_\Manacst^\Mancst(Y)$ for the associated manifolds with corners. Then $X=\ti X$, $Y=\ti Y$ as topological spaces, though the (a-)atlases on $X,Y$ and $\ti X,\ti Y$ are different. We claim that $f:\ti X\ra\ti Y$ is a strongly smooth map of manifolds with corners. To see this, let $(U_a,\phi_a)$ and $(V_b,\psi_b)$ be a-charts on $X,Y$, and $(\ti U_a,\ti\phi_a)$ and $(\ti V_b,\ti\psi_b)$ the corresponding charts on $\ti X,\ti Y$. Write $\dot U_a=(f\ci\phi_a)^{-1}(\psi_b(V_b))$, so that $\dot U_a\subseteq U_a\subseteq\R^{k,m}$ is open. Then by definition $\psi_b^{-1}\ci f\ci\phi_a:\dot U_a\ra V_b$ is a strongly a-smooth map of open subsets $\dot U_a\subseteq\R^{k,m}$, $V_b\subseteq\R^{l,n}$. 

The argument above that $\ti\phi_b^{-1}\ci\ti\phi_a$ is strongly smooth uses only that $\phi_b^{-1}\ci\phi_a$ is strongly a-smooth, not that it is an a-diffeomorphism. Thus, exactly the same argument shows that $\ti\psi_b^{-1}\ci f\ci\ti\phi_a:\dot{\ti U}_a\ra\ti V_b$ is strongly smooth. As this holds for all charts $(\ti U_a,\ti\phi_a),(\ti V_b,\ti\psi_b)$ in atlases for $\ti X,\ti Y$, we see that $f:\ti X\ra\ti Y$ is strongly smooth. Define $F_\Manacst^\Mancst(f)=f$. Then $F_\Manacst^\Mancst:\Manacst\ra\Mancst$ is a functor. It is faithful, but not full. It maps interior morphisms to interior morphisms, and extends to $F_\cManacst^\cMancst:\cManacst\ra\cMancst$ in the obvious way.
\label{ac3def6}
\end{dfn}

\begin{ex}{\bf(a)} Define $f:\lb 0,\iy)\ra\lb 0,\iy)$ by $f(x)=x^\al$ for $\al>0$. Identify 
$F_\Manacst^\Mancst(\lb 0,\iy))$ with $[0,\iy)$ via $\vp^{-1}:\lb 0,\iy)\ra[0,\iy)$ in \eq{ac3eq9}. Then $\ti f=F_\Manacst^\Mancst(f)$ is a smooth map $[0,\iy)\ra[0,\iy)$, in the commutative diagram
\begin{equation*}
\xymatrix@C=90pt@R=13pt{ *+[r]{\lb 0,\iy)} \ar[r]_{\vp^{-1}} \ar[d]^f & *+[l]{[0,\iy)} \ar[d]_{\ti f} \\ *+[r]{\lb 0,\iy)} \ar[r]^{\vp^{-1}} & *+[l]{[0,\iy).\!}
 }
\end{equation*}
Substituting in \eq{ac3eq8}--\eq{ac3eq9} shows that
\begin{equation*}
\ti f(x)=\ha\bigl(\al(x-x^{-1})+\sqrt{\al^2(x-x^{-1})^2+4}\,\bigr),\quad x>0,\quad\text{and}\quad \ti f(0)=0.
\end{equation*}
This is smooth, with $\ti f(x)\approx \al^{-1}x$ for small $x$.
\smallskip

\noindent{\bf(b)} Define $g:\lb 0,\iy)^2\ra\lb 0,\iy)$ by $g(x,y)=xy$. It is a-smooth, but not strongly a-smooth. Identify $F_\Manacst^\Mancst(\lb 0,\iy)),F_\Manacst^\Mancst(\lb 0,\iy)^2)$ with $[0,\iy),[0,\iy)^2$ via $\vp^{-1},\vp^{-1}\t\vp^{-1}$. Define $\ti g:[0,\iy)^2\ra[0,\iy)$ by the commutative diagram
\begin{equation*}
\xymatrix@C=90pt@R=13pt{ *+[r]{\lb 0,\iy)^2} \ar[r]_{\vp^{-1}\t\vp^{-1}} \ar[d]^g & *+[l]{[0,\iy)^2} \ar[d]_{\ti g} \\ *+[r]{\lb 0,\iy)} \ar[r]^{\vp^{-1}} & *+[l]{[0,\iy).\!}
 }
\end{equation*}
Substituting in \eq{ac3eq8}--\eq{ac3eq9} shows that
\begin{equation*}
\ti g(x,y)=\begin{cases} \ha\bigl(x\!-\!x^{-1}\!+\!y\!-\!y^{-1}+\sqrt{(x\!-\!x^{-1}\!+\!y\!-\!y^{-1})^2\!+\!4}\,\bigr), & x,y>0, \\ 0, & xy=0.
\end{cases}
\end{equation*}
This satisfies $\ti g(x,y)\approx xy/(x+y)$ for small $x,y$, so $\ti g$ is not smooth at $(0,0)$. 

This example shows that in Definition \ref{ac3def6}, if $f:X\ra Y$ is a-smooth, but not strongly a-smooth, then $f:\ti X\ra\ti Y$ may not be smooth (in fact, it is never smooth), and $F_\Manacst^\Mancst$ cannot be extended to a functor~$\Manac\ra\Manc$.

This example is discussed by Hassell, Mazzeo and Melrose \cite[\S 2.5, Ex.~1]{HMM}. They also note that $\ti g(x,y)\approx xy/(x+y)$ is non-smooth after doing their `total logarithmic blowup', and say that in their framework one should perform a further (ordinary, not logarithmic) blowup $\pi:\hat X\ra\ti X$ of $\ti X$ at $(0,0)$, and then the pullback $\hat g=\pi^*(\ti g)$ is smooth on~$\hat X$.
\label{ac3ex5}
\end{ex}

\begin{rem}{\bf(a)} In \S\ref{ac31}--\S\ref{ac32} we chose to define manifolds with a-corners using a-smooth functions, including the $O(x^\al)$ decay condition \eq{ac3eq2} near boundary hypersurfaces $x=0$, rather than the weaker r-smooth functions.

If we had used r-smooth functions, the definition of $F_\Manacst^\Mancst$ above would not have worked, as \eq{ac3eq2} is needed in \eq{ac3eq13} to prove compatibility of charts.
\smallskip

\noindent{\bf(b)} Composing $\smash{F_\Manc^\Manac,F_\Manacst^\Mancst}$ gives functors
\begin{equation*}
F_\Manacst^\Mancst\ci F_\Manc^\Manac:\Mancst\ra\Mancst,\;\>
F_\Manc^\Manac\ci F_\Manacst^\Mancst:\Manacst\ra\Manacst.
\end{equation*}
These are not the identities, and both are faithful but not full. If $X$ is a manifold with corners, one can prove that there is a noncanonical diffeomorphism $X\cong F_\Manacst^\Mancst\ci F_\Manc^\Manac(X)$ in $\Manc$. But Examples \ref{ac3ex3} and \ref{ac4ex3} show that there exist manifolds with a-corners $X$ with $X\not\cong F_\Manc^\Manac\ci F_\Manacst^\Mancst(X)$ in~$\Manac$.
\smallskip

\noindent{\bf(c)} Hofer--Wysocki--Zehnder's {\it polyfolds\/} \cite{Hofe,HWZ1,HWZ2} and Fukaya--Oh--Ohta--Ono's {\it Kuranishi spaces\/} \cite{FOOO1,FOOO2,FuOn} are geometrical spaces used to describe moduli spaces $\oM$ of $J$-holomorphic curves in symplectic geometry, including singular curves. 

Neglecting issues to do with orbifolds, for Kuranishi spaces $\oM$ is locally modelled on $s^{-1}(0)$, where $V$ is a manifold with corners, $E\ra V$ a vector bundle, and $s:V\ra E$ a smooth section. Singular curves correspond to $v\in s^{-1}(0)$ lying in a boundary stratum of $V$. In the more complicated polyfold picture, roughly the same holds, but with $V,E$ infinite-dimensional.

An aspect of both theories is related to our functor $F_\Manacst^\Mancst$ above. One describes curves $\Si_x$ near a singular curve $\Si_0$ using a `neck length' $x\in[0,\iy)$, where the curve is singular when $x=0$. For quantities such as the Kuranishi section $s$, one proves estimates of the form $\frac{\pd^k}{\pd x^k}s=O(x^{\al-k})$ for some $\al>0$ and all $k=0,1,\ldots.$ Thus, if $x\in[0,\iy)$ is a coordinate normal to $\pd V$ in the obvious way, the Kuranishi section $s$ may not be smooth at $x=0$. Both groups deal with this by changing coordinates from $x$ to $\ti x$ with $x=e^{-1/\ti x}$, \cite[\S 2.1]{HWZ1}, \cite[p.~777]{FOOO1}, \cite[\S 1]{FOOO2}, so that $s$ is smooth as a function of $\ti x$. In the polyfold theory, the choice of coordinate change function $\ti x\mapsto e^{-1/\ti x}$ is called a `gluing profile'.

In \eq{ac3eq8} we have $\vp(\ti x)\approx e^{-1/\ti x}$ for $\ti x>0$ small, so our definition of $F_\Manacst^\Mancst$ essentially applies the change of coordinates in $V$ normal to $\pd V$ in~\cite{FOOO1,Hofe,HWZ1}.

We advocate the following point of view, discussed further in~\S\ref{ac63}:
\begin{itemize}
\setlength{\itemsep}{0pt}
\setlength{\parsep}{0pt}
\item[(i)] To describe `moduli spaces with corners' $\oM$ in differential geometry, such as those in \cite{FOOO1,FOOO2,Hofe,HWZ1,HWZ2}, the natural smooth structure to use is {\it manifolds with a-corners}, not manifolds with corners. So, for example, we should locally model $\oM$ as $s^{-1}(0)$ for $V$ a manifold with a-corners, $E\ra V$ a vector bundle, and $s:V\ra E$ an a-smooth section.
\item[(ii)] The method of gluing profiles in \cite{Hofe,HWZ1}, and its analogue in \cite{FOOO1,FOOO2} are equivalent to first constructing $V,E,s$ in $\Manac$ as in (i), and then applying $F_\Manacst^\Mancst$ to $V,E,s$ to get $\ti V,\ti E,\ti s$ in $\Manc$.
\item[(iii)] There may be advantages to working with manifolds with a-corners, rather than manifolds with corners, in moduli problems of this kind.

One advantage is that the a-smooth structure is more canonical. Another is that $F_\Manacst^\Mancst$ only works for strongly a-smooth maps. There are important morphisms between moduli spaces of $J$-holomorphic curves (e.g.\ forgetting a marked point) which may be locally modelled on a-smooth, but not strongly a-smooth, maps. Such morphisms will be a-smooth when written using manifolds with a-corners, but will be {\it non-smooth\/} when described using manifolds with corners and gluing profiles as in~\cite{FOOO1,FOOO2,Hofe,HWZ1}.
\end{itemize}
\label{ac3rem2}
\end{rem}

\subsection{Real analytic manifolds with a-corners}
\label{ac34}

We recall the definition of {\it real analytic\/} for ordinary manifolds with corners:

\begin{dfn} Let $U\subseteq\R^m_k$ be open, and $f:U\ra\R$ be a smooth function. We call $f$ {\it real analytic\/} if either of the following two equivalent conditions hold:
\begin{itemize}
\setlength{\itemsep}{0pt}
\setlength{\parsep}{0pt}
\item[(i)] There exists an open neighbourhood $V$ of $U$ in $\C^m$ and a holomorphic function $g:U\ra\C$ with $g\vert_U=f$, where~$U\subseteq\R^m_k\subseteq\R^m\subseteq\C^m$.
\item[(ii)] For each $x\in U$, the Taylor series of $f$ at $x$ converges absolutely to $f$ in a neighbourhood of $x$ in~$U$. 
\end{itemize}

Let $U\subseteq\R^m_k$ and $V\subseteq\R^n_l$ be open, and $f:U\ra V$ be a smooth function, as in Definition \ref{ac2def1}. We say that $f$ is {\it real analytic\/} if in Definition \ref{ac2def1}(b)(i), the local functions $F_j:U\ra(0,\iy)$ are real analytic as functions~$U\ra\R$.

We can now define the category $\Manrac$ of {\it real analytic manifolds with corners}, with full subcategories $\Manra\subset\Manrab\subset\Manrac$ of {\it real analytic manifolds\/} and {\it real analytic manifolds with boundary}, by following Definition \ref{ac2def2} but using real analytic maps $f:U\ra V$ rather than smooth maps throughout, for open $U\subseteq\R^m_k$, $V\subseteq\R^n_l$.
\label{ac3def7}
\end{dfn}

We generalize Definition \ref{ac3def7} to manifolds with a-corners:

\begin{dfn} Let $\R^{k,m}=\lb 0,\iy)^k\t\R^{m-k}$ have coordinates $(x_1,\ldots,x_m)$. We consider the complexification of $\R^{k,m}$ to be $\R^{k,m}\t\R^m=\R^{k,2m}$, with coordinates $(x_1,\ldots,x_m,y_1,\ldots,y_m)$, where we embed $\R^{k,m}\hookra\R^{k,m}\t\R^m$ by $(x_1,\ldots,x_m)\mapsto(x_1,\ldots,x_m,0,\ldots,0)$. Let $V\subseteq\R^{k,m}\t\R^m$ be open, and $g:V\ra\C$ be a-smooth. We call $g$ {\it a-holomorphic\/} if
\e
\begin{aligned}
\ts x_j\frac{\pd g}{\pd x_j}+i\frac{\pd g}{\pd y_j}&=0, & j&=1,\ldots,k,\quad\text{and} \\
\ts\frac{\pd g}{\pd x_j}+i\frac{\pd g}{\pd y_j}&=0, & j&=k+1,\ldots,m.
\end{aligned}
\label{ac3eq14}
\e

Let $U\subseteq\R^{k,m}$ be open, and $f:U\ra\R$ be a-smooth. We call $f$ {\it real analytic\/} if there exists an open neighbourhood $V$ of $U$ in $\R^{k,m}\t\R^m$ and an a-holomorphic function $g:V\ra\C$ with $g\vert_U=f$.

Let $U\subseteq\R^{k,m}$ and $V\subseteq\R^{l,n}$ be open, and $f:U\ra V$ be an a-smooth function, as in Definition \ref{ac3def3}. We say that $f$ is {\it real analytic\/} if in Definition \ref{ac3def3}(b)(i), the local a-smooth functions $F_j:U\ra(0,\iy)$ are real analytic as functions~$U\ra\R$.

We can now define the category $\Manraac$ of {\it real analytic manifolds with a-corners}, with full subcategories $\Manra\subset\Manraab\subset\Manraac$ of {\it real analytic manifolds\/} (as in Definition \ref{ac3def7}) and {\it real analytic manifolds with a-boundary}, by following Definition \ref{ac3def4} but using real analytic maps $f:U\ra V$ rather than a-smooth maps throughout, for open $U\subseteq\R^{k,m}$, $V\subseteq\R^{l,n}$.
\label{ac3def8}
\end{dfn}

\begin{ex} Let $\al>0$ and $\be\in\R$. Define $g:\lb 0,1)\t\R\ra\C$ by 
\begin{equation*}
g(x,y)=\begin{cases} x^\al \cdot \exp(i\al y)\cdot(\log x+iy)^\be, & x>0, \\ 0, & x=0. \end{cases}
\end{equation*}
Then $g$ is a-smooth and satisfies \eq{ac3eq14}, so $g$ is a-holomorphic. Therefore $f:\lb 0,1)\ra\R$ given by
\begin{equation*}
f(x)=\begin{cases} x^\al(\log x)^\be, & x>0, \\ 0, & x=0, \end{cases}
\end{equation*}
is a-smooth and real analytic.
\label{ac3ex6}
\end{ex}

\begin{rem}{\bf(a)} For manifolds with corners, as in Definition \ref{ac3def7}(i),(ii) we have two equivalent definitions of real analytic. For manifolds with a-corners, we have generalized (i). We do not know an equivalent definition (ii), involving `Taylor series' for a-smooth or a-holomorphic functions at corner points. 

As in \S\ref{ac55}, in Melrose's theory of analysis on manifolds with corners there is a notion of {\it polyhomogeneous conormal function}, which is like a kind of Taylor series at boundary and corner points, but our real analytic functions are not equivalent to functions which are locally the limits of a polyhomogeneous conormal expansion, as Example \ref{ac3ex6} with $\be\notin\N$ shows.
\smallskip

\noindent{\bf(b)} There are many situations in which solutions of elliptic p.d.e.s are known to be real analytic. For example, if $(X,g)$ is an Einstein Riemannian manifold, then $X$ has a unique real analytic structure with respect to which $g$ is real analytic. The author expects that in the same way, solutions of elliptic p.d.e.s on manifolds with a-corners will often be real analytic in the sense above.
\label{ac3rem3}
\end{rem}

\subsection{Manifolds with both corners and a-corners}
\label{ac35}

So far we have studied manifolds with corners $\Manc$, locally modelled on $\R^m_l=[0,\iy)^l\t\R^{m-l}$, and manifolds with a-corners $\Manac$, locally modelled on $\R^{k,m}=\lb 0,\iy)^k\t\R^{m-k}$. We will now combine both into a category $\Mancac$ of {\it manifolds with corners and a-corners}, locally modelled on $\R^{k,m}_l=\ab\lb 0,\iy)^k\t[0,\iy)^l\t\R^{m-k-l}$, containing $\Manc,\Manac$ as full subcategories.

Here are the analogues of Definitions \ref{ac2def1}, \ref{ac3def2} and~\ref{ac3def3}.

\begin{dfn} For $0\le k,l\le k+l\le m$, write $\R^{k,m}_l=\lb 0,\iy)^k\t[0,\iy)^l\t\R^{m-k-l}$, with coordinates $(x_1,\ldots,x_m)$ with $x_1,\ldots,x_k$ in $\lb 0,\iy)$, $x_{k+1},\ldots,x_{k+l}$ in $[0,\iy)$, and $x_{k+l+1},\ldots,x_m$ in $\R$.

Let $U\subseteq \R^{k,m}_l$ be open and $f:U\ra\R$ be continuous. The {\it b-derivative\/} of $f$ (if it exists) is a map ${}^b\pd f:U\ra\R^m$ given by \eq{ac3eq1}, where now for $i=k+1,\ldots,l$ and $x_i=0$, by $\frac{\pd f}{\pd x_i}(x_1,\ldots,x_m)$ we mean the one-sided derivative at $0\in [0,\iy)$. \begin{itemize}
\setlength{\itemsep}{0pt}
\setlength{\parsep}{0pt}
\item[(i)] We call $f$ {\it r-differentiable\/} if ${}^b\pd f:U\ra\R^m$ exists and is continuous.
\item[(ii)] We call $f$ {\it r-smooth\/} if ${}^b\pd^lf:U\ra\bigot^l\R^m$ is r-differentiable for $l\ge 0$.
\item[(iii)] We call $f$ {\it differentiable\/} if it is r-differentiable and for any compact subset $S\subseteq U$ and $i=1,\ldots,k$, there exist $C,\al>0$ such that \eq{ac3eq2} holds.
\item[(iv)] We call $f$ {\it smooth\/} if ${}^b\pd^lf:U\ra\bigot^l\R^m$ is differentiable for~$l\ge 0$.
\end{itemize}
Essentially, $f$ is smooth if it is a-smooth in the $x_1,\ldots,x_k$ variables in the sense of \S\ref{ac31}, and smooth in the $x_{k+1},\ldots,x_m$ variables in the sense of~\S\ref{ac21}.
\label{ac3def9}
\end{dfn}

\begin{dfn} Let $U\subseteq\R^{k,m}_l$ and $V\subseteq \R^{k',n}_{l'}$ be open, and $f=(f_1,\ldots,f_n):U\ra V$ be a continuous map, so that $f_j=f_j(x_1,\ldots,x_m)$ maps $U\ra\lb 0,\iy)$ for $j=1,\ldots,k'$ and $U\ra[0,\iy)$ for $j=k'+1,\ldots,k'+l'$ and $U\ra\R$ for $j=k'+l'+1,\ldots,n$. Then we say:
\begin{itemize}
\setlength{\itemsep}{0pt}
\setlength{\parsep}{0pt}
\item[(a)] $f$ is {\it r-smooth\/} if $f_j:U\ra\R$ is r-smooth in the sense of Definition \ref{ac3def9} for $j=k'+l'+1,\ldots,n$, and every $u=(x_1,\ldots,x_m)\in U$ has an open neighbourhood $\ti U$ in $U$ such that for each $j=1,\ldots,k'$, either:
\begin{itemize}
\setlength{\itemsep}{0pt}
\setlength{\parsep}{0pt}
\item[(i)] we may uniquely write
\begin{equation*}
f_j(\ti x_1,\ldots,\ti x_m)=F_j(\ti x_1,\ldots,\ti x_m)\cdot\ti x_1^{a_{1,j}}\cdots\ti x_{k+l}^{a_{k+l,j}}
\end{equation*}
for all $(\ti x_1,\ldots,\ti x_m)\in\ti U$, with $F_j:\ti U\ra(0,\iy)$ r-smooth as in Definition \ref{ac3def9}, and $a_{1,j},\ldots,a_{k+l,j}\in[0,\iy)$, with $a_{i,j}=0$ if $x_i\ne 0$;~or 
\item[(ii)] $f_j\vert_{\smash{\ti U}}=0$.
\end{itemize}
and for each $j=k'+1,\ldots,k'+l'$, either:
\begin{itemize}
\setlength{\itemsep}{0pt}
\setlength{\parsep}{0pt}
\item[(iii)] we may uniquely write
\begin{equation*}
f_j(\ti x_1,\ldots,\ti x_m)=F_j(\ti x_1,\ldots,\ti x_m)\cdot\ti x_{k+1}^{b_{k+1,j}}\cdots\ti x_{k+l}^{b_{k+l,j}}
\end{equation*}
for all $(\ti x_1,\ldots,\ti x_m)\in\ti U$, with $F_j:\ti U\ra(0,\iy)$ r-smooth as in Definition \ref{ac3def9}, and $b_{k+1,j},\ldots,b_{k+l,j}\in\N$, with $b_{i,j}=0$ if $x_i\ne 0$;~or 
\item[(iv)] $f_j\vert_{\smash{\ti U}}=0$.
\end{itemize}
\item[(b)] $f$ is {\it smooth\/} if (a) holds, but requiring $f_j,F_j$ to be smooth rather than r-smooth, in the sense of Definition \ref{ac3def9}.
\item[(c)] $f$ is {\it interior\/} if it is smooth, and cases (a)(ii),(iv) do not occur.
\item[(d)] $f$ is {\it b-normal\/} if it is interior, and in cases (a)(i),(iii) for each $i=1,\ldots,k$ we have $a_{i,j}>0$ for at most one $j=1,\ldots,k'$, and for each $i=k+1,\ldots,k+l$ we have $a_{i,j}>0$ or $b_{i,j}>0$ for at most one $j=1,\ldots,k'+l'$.
\item[(e)] $f$ is {\it strongly smooth\/} if it is smooth, and in case (a)(i) for each $j=1,\ldots,k'$ we have $a_{i,j}>0$ for at most one $i=1,\ldots,k+l$, and in case (a)(iii) for each $j=k'+1,\ldots,k'+l'$ we have $b_{i,j}=1$ for at most one $i=k+1,\ldots,k+l$ and $b_{i,j}=0$ otherwise. 
\item[(f)] $f$ is a {\it diffeomorphism\/} if it is a smooth bijection with smooth inverse.
\end{itemize}
All the classes (a)--(f) include identities and are closed under composition.
\label{ac3def10}
\end{dfn}

We can now follow Definition \ref{ac3def4} with obvious changes, to define a category $\Mancac$ of {\it manifolds with corners and a-corners}. An object of $\Mancac$ is a second countable Hausdorff topological space $X$ equipped with a maximal atlas of charts $(U_a,\phi_a)$ for $U_a\subseteq\R^{k,m}_l$ open and $\phi_a:U_a\ra X$ a homeomorphism with an open set $\phi_a(U_a)\subseteq X$, where $m=\dim X$. The transition functions $\phi_b^{-1}\ci\phi_a$ between charts should be diffeomorphisms between open subsets of $\R^{k,m}_l,\R^{k',m}_{l'}$ in the sense of Definition~\ref{ac3def10}(f). 

Let $X,Y$ be objects in $\Mancac$, and $f:X\ra Y$ a continuous map. We say that $f$ is {\it r-smooth}, or {\it smooth}, or {\it interior}, or {\it b-normal}, or {\it strongly smooth}, if $\psi_b^{-1}\ci f\ci\phi_a$ is an r-smooth, \ldots, strongly smooth map between open subsets of $\R^{k,m}_l,\R^{k',n}_{l'}$ in the sense of Definition \ref{ac3def10}(a)--(e) for all charts $(U_a,\phi_a)$ on $X$ and $(U_b,\phi_b)$ on $Y$. Morphisms in $\Mancac$ are defined to be smooth maps.

We write $\Mancacin,\Mancacst$ for the subcategories of $\Mancac$ with morphisms interior maps, and strongly a-smooth maps, respectively. As for $\cManc,\ab\cManac$, we define $\cMancac$ to have objects disjoint unions $\coprod_{m=0}^\iy X_m$, where $X_m$ is a manifold with corners and a-corners of dimension $m$, and morphisms smooth maps, and we define subcategories $\cMancacin,\cMancacst\subset\cMancac$ with morphisms interior, and strongly smooth, maps.

Manifolds with corners $\Manc$, and manifolds with a-corners $\Manac$, are both full subcategories of $\Mancac$, of objects $X$ covered by charts $(U_a,\phi_a)$ with $U_a\subseteq\R^{k,m}_l$ open where $k=0$ for $\Manc$ and $l=0$ for~$\Manac$.

The functor $F_\Manc^\Manac:\Manc\ra\Manac$ in Definition \ref{ac3def5} extends to a functor $F_\Mancac^\Manac:\Mancac\ra\Manac$, which restricts to $F_\Manc^\Manac$ on $\Manc\subset\Mancac$ and to the identity on $\Manac\subset\Mancac$. It replaces charts $(U_a,\phi_a)$ on $X$ with $U_a\subseteq\R^{k,m}_l$ open by $(U_a,\phi_a)$ regarded as an a-chart with $U_a\subseteq\R^{k+l,m}$ open.

The functor $F_\Manacst^\Mancst:\Manacst\ra\Mancst$ in Definition \ref{ac3def6} extends to a functor $F_\Mancacst^\Mancst:\Mancacst\ra\Mancst$, which restricts to $F_\Manacst^\Mancst$ on $\Manacst\subset\Mancacst$ and to the identity on $\Mancst\subset\Mancacst$. It replaces charts $(U_a,\phi_a)$ on $X$ with $U_a\subseteq\R^{k,m}_l$ open by charts $(\ti U_a,\ti\phi_a)$ on $X$ with $\ti U_a\subseteq\R^m_{k+l}$ open, where only the variables $x_1,\ldots,x_k$ on $U_a$ are modified using~$\vp$.

\begin{rem}{\bf(a)} Much of what we say about $\Manac$ in this paper has an extension to $\Mancac$, which is usually fairly obvious, and left to the reader.
\smallskip

\noindent{\bf(b)} Manifolds with both corners and a-corners will be important in geometric problems involving boundary conditions of both kinds, as we explain in~\S\ref{ac6}.
\smallskip

\noindent{\bf(c)} By combining Definitions \ref{ac3def7} and \ref{ac3def8}, we can define {\it real analytic\/} manifolds with corners and a-corners, in an obvious way.\label{ac3rem4}
\end{rem}

\section[Differential geometry of manifolds with a-corners]{Differential geometry of manifolds with \\ a-corners}
\label{ac4}

\subsection{Boundaries and corners of manifolds with a-corners}
\label{ac41}

This section follows \S\ref{ac22} for ordinary corners essentially without change. Proofs of results below either follow those for ordinary corners in \cite{Joyc10,Joyc13}, or are obvious.

\begin{dfn} Let $X$ be a manifold with a-corners, of dimension $m$. There is a natural {\it depth stratification\/} $X=\coprod_{l=0}^mS^l(X)$, such that if $(U,\phi)$ is an a-chart on $X$ with $U\subseteq\R^{k,m}$ open and $(x_1,\ldots,x_m)\in U$ then $\phi(x_1,\ldots,x_m)\in S^l(X)$ if $l$ out of $x_1,\ldots,x_k$ are zero. Then $\overline{S^l(X)}=
\bigcup_{k=l}^m S^k(X)$. The {\it interior\/} of $X$ is
$X^\ci=S^0(X)$. Each $S^l(X)$ is an $(m-l)$-manifold without boundary.

A {\it local $k$-corner component\/ $\ga$ of\/ $X$ at\/} $x$ is a local choice of connected component of $S^k(X)$ near $x$. When $k=1$, we call $\ga$ a {\it local boundary component}. As sets, define the {\it boundary\/} $\pd X$ and {\it k-corners\/} $C_k(X)$ for $k=0,1,\ldots,m$ by
\begin{align*}
\pd X&=\bigl\{(x,\be):\text{$x\in X$, $\be$ is a local boundary
component of $X$ at $x$}\bigr\},\\
C_k(X)&=\bigl\{(x,\ga):\text{$x\in X$, $\ga$ is a local $k$-corner 
component of $X$ at $x$}\bigr\}.
\end{align*}
Define $i_X:\pd X\ra X$ and $\Pi:C_k(X)\ra X$ by $i_X:(x,\be)\mapsto x$, $\Pi:(x,\ga)\mapsto x$.

If $(U,\phi)$ is an a-chart on $X$ with $U\subseteq\R^{k,m}$ open, then for each $i=1,\ldots,k$ we can define an a-chart $(U_i,\phi_i)$ on $\pd X$ by
\begin{align*}
&U_i\!=\!\bigl\{(x_1,\ldots,x_{m-1})\!\in\! \R^{k-1,m-1}:
(x_1,\ldots,x_{i-1},0,x_i,\ldots,x_{m-1})\!\in\!
U\subseteq\R^{k,m}\bigr\},\\
&\phi_i:(x_1,\ldots,x_{m-1})\longmapsto\bigl(\phi
(x_1,\ldots,x_{i-1}, 0,x_i,\ldots,x_{m-1}),\phi_*(\{x_i=0\})\bigr).
\end{align*}
The set of all such a-charts on $\pd X$ forms an a-atlas, making $\pd X$ into a manifold with a-corners of dimension $m-1$, and $i_X:\pd X\ra X$ into an a-smooth (but not interior) map. Similarly, we make $C_k(X)$ into an $(m-k)$-manifold with a-corners, and $\Pi:C_k(X)\ra X$ into an a-smooth map.

We call $X$ a {\it manifold without boundary\/} if $\pd X=\es$, and
a {\it manifold with a-boundary\/} if $\pd^2X=\es$. We write $\Man$ and $\Manab$ for the full subcategories of $\Manac$ with objects manifolds without boundary, and manifolds with a-boundary, so that $\Man\subset\Manab\subset\Manac$. This definition of $\Man$ is equivalent to the usual definition of the category of manifolds.
\label{ac4def1}
\end{dfn}

For $X\in\Manac$ and $k\ge 0$, as in \eq{ac2eq2}--\eq{ac2eq5} there are natural identifications 
\ea
\begin{split}
\pd^kX&\cong\bigl\{(x,\be_1,\ldots,\be_k):\text{$x\in X,$
$\be_1,\ldots,\be_k$ are distinct}\\
&\qquad\qquad\text{local boundary components for $X$ at $x$}\bigr\},
\end{split}
\label{ac4eq1}\\
\begin{split}
C_k(X)&\cong\bigl\{(x,\{\be_1,\ldots,\be_k\}):\text{$x\in X,$
$\be_1,\ldots,\be_k$ are distinct}\\
&\qquad\qquad\text{local boundary components for $X$ at $x$}\bigr\},
\end{split}
\label{ac4eq2}\\
C_k(X)&\cong\pd^kX/S_k,\qquad\qquad \pd C_k(X)\cong C_k(\pd X), 
\label{ac4eq3}
\ea
where $S_k$ acts freely and smoothly on $\pd^kX$ by permuting $\be_1,\ldots,\be_k$ in~\eq{ac4eq1}.

For products of manifolds with a-corners we have natural a-diffeomorphisms
\ea
\pd(X\t Y)&\cong (\pd X\t Y)\amalg (X\t\pd Y),
\label{ac4eq4}\\
C_k(X\t Y)&\cong \ts\coprod_{i,j\ge 0,\; i+j=k}C_i(X)\t C_j(Y).
\label{ac4eq5}
\ea
The analogue of Lemma \ref{ac2lem} holds for manifolds with a-corners.

\begin{dfn} Define the {\it corners\/} $C(X)$ of a manifold with a-corners $X$ by
\begin{align*}
C(X)&=\ts\coprod_{k=0}^{\dim X}C_k(X)\\
&=\bigl\{(x,\ga):\text{$x\in X$, $\ga$ is a local $k$-corner 
component of $X$ at $x$, $k\ge 0$}\bigr\},
\end{align*}
considered as an object of $\cManac$ in Definition \ref{ac3def4}. Define $\Pi:C(X)\ra X$ by $\Pi:(x,\ga)\mapsto x$. This is a-smooth, that is, a morphism in~$\cManac$.

Let $f:X\ra Y$ be a morphism in $\Manac$, and suppose $\ga$ is a local $k$-corner component of $X$ at $x\in X$. Then there is a unique local $l$-corner component $f_*(\ga)$ of $Y$ at $f(x)$ for some $l\ge 0$, such that if $V,\ti V$ are sufficiently small open neighbourhoods of $x,f(x)$ in $X,Y$, so that $\ga,f_*(\ga)$ determine connected components $W,\ti W$ of $V\cap S^k(X)$ and $\ti V\cap S^l(Y)$, then~$f(W)\cap\ti W\ne\es$. 

Define a map $C(f):C(X)\ra C(Y)$ by $C(f):(x,\ga)\mapsto (f(x),f_*(\ga))$. Then $C(f)$ is an interior morphism in $\cManac$, and $C:\Manac\ra\cManacin\subset\cManac$ is a functor, which we call the {\it corner functor}.
\label{ac4def2}
\end{dfn}

Equations \eq{ac4eq3} and \eq{ac4eq5} imply that if $X,Y$ are manifolds with a-corners, we have natural isomorphisms
\ea
\pd C(X)&\cong C(\pd X),
\label{ac4eq6}\\
C(X\t Y)&\cong C(X)\t C(Y).
\label{ac4eq7}
\ea
The corner functor $C$ {\it preserves products and direct products}. That is, if $f:W\ra Y,$ $g:X\ra Y,$ $h:X\ra Z$ are a-smooth then the following commute:
\begin{equation*}
\xymatrix@C=60pt@R=20pt{ *+[r]{C(W\t X)} \ar[d]^\cong \ar[r]_{C(f\t
h)} & *+[l]{C(Y\t Z)} \ar[d]_\cong \\ *+[r]{C(W)\!\t\! C(X)}
\ar[r]^{\raisebox{8pt}{$\st C(f) \t C(h)$}} &
*+[l]{C(Y)\!\t\! C(Z),} }\;
\xymatrix@C=65pt@R=3pt{ & C(Y\t Z) \ar[dd]^\cong \\
C(X) \ar[ur]^(0.4){C((g,h))} \ar[dr]_(0.4){(C(g),C(h))} \\
& C(Y)\!\t\! C(Z), }
\end{equation*}
where the columns are the isomorphisms~\eq{ac4eq7}.

The functor $F_\Manc^\Manac:\Manc\ra\Manac$ in Definition \ref{ac3def5} preserves the depth stratification, and local boundary and corner components. For $X\in\Manc$ there are canonical a-diffeomorphisms
\begin{equation*}
\pd F_\Manc^\Manac(X)\cong F_\Manc^\Manac(\pd X),\qquad C_k\bigl(F_\Manc^\Manac(X)\bigr)\cong F_\Manc^\Manac\bigl(C_k(X)\bigr).
\end{equation*}
The following diagram of functors commutes:
\begin{equation*}
\xymatrix@C=120pt@R=15pt{ *+[r]{\Manc} \ar[r]_(0.4){F_\Manc^\Manac} \ar[d]^C & *+[l]{\Manac} \ar[d]_C \\
*+[r]{\cManc} \ar[r]^(0.6){F_\cManc^\cManac} & *+[l]{\cManac.\!} }
\end{equation*}
The same holds for the functor $F_\Manacst^\Mancst:\Manacst\ra\Mancst$ of~\S\ref{ac33}.

\begin{rem} The material of this section extends to the category $\Mancac$ of \S\ref{ac35}. Manifolds with corners and a-corners $X$ have natural depth stratifications $X=\coprod_{l=0}^{\dim X}S^l(X)$, defined in the obvious way. We follow Definitions \ref{ac4def1}--\ref{ac4def2} essentially without change to define $\pd X,C_k(X),C(X)$ and the corner functor $C:\Mancac\ra\cMancac$.

For $\Mancac$, the boundary $\pd X$ splits as $\pd X=\pd^{\rm c} X\amalg\pd^{\rm ac}X$, where $\pd^{\rm c} X$ is the `ordinary boundary' and $\pd^{\rm ac}X$ the `a-boundary' of $X$. Similarly $C_k(X)$ has a natural decomposition $C_k(X)=\coprod_{i+j=k}C_{i,j}(X)$, where $C_{i,j}(X)$ is locally the intersection of $i$ ordinary boundaries and $j$ a-boundaries, so that locally $X\simeq C_{i,j}(X)\t[0,\iy)^i\t\lb0,\iy)^j$. We can also define functors $C^{\rm c},C^{\rm ac}:\Mancac\ra\cMancac$ with $C^{\rm c}(X)=\coprod_{i\ge 0}C_{i,0}(X)$ and $C^{\rm ac}(X)=\coprod_{j\ge 0}C_{0,j}(X)$, so that $C^{\rm c}$ takes only ordinary corners, and $C^{\rm ac}$ only a-corners.
\label{ac4rem1}
\end{rem}

\subsection{B-(co)tangent bundles of manifolds with a-corners}
\label{ac42}

As in \S\ref{ac23}, for manifolds with ordinary corners $X$, we can define the tangent bundle $TX$, and the b-tangent bundle ${}^bTX$. We explain in Remark \ref{ac4rem2}(a) that tangent bundles $TX$ do not make sense for manifolds with a-corners. But we will now define b-tangent bundles ${}^bTX$ for manifolds with a-corners.

\begin{dfn} Let $X$ be a manifold with a-corners. A {\it vector
bundle\/ $E\ra X$ of rank\/} $k$ is a manifold with a-corners $E$ and
an a-smooth map $\pi:E\ra X$, such that each fibre $E_x:=\pi^{-1}(x)$ for $x\in X$ is given the structure of a real vector space, and $X$ may be covered by open $U\subseteq X$ with a-diffeomorphisms $\pi^{-1}(U)\cong
U\t\R^k$ identifying $\pi\vert_{\pi^{-1}(U)}:\pi^{-1}(U)\ra U$ with
the projection $U\t\R^k\ra\R^k$, and the vector space structure on
$E_x$ with that on $\{x\}\t\R^k\cong\R^k$, for each $x\in U$. A {\it section\/} of $E$ is an a-smooth map $s:X\ra E$ with $\pi\ci
s=\id_X$.

We write $\Ga^\iy(E)$ for the vector space of a-smooth sections of $E$, and $C^\iy(X)$ for the $\R$-algebra of a-smooth functions $X\ra\R$. Then $\Ga^\iy(E)$ is a $C^\iy(X)$-module. Morphisms of vector bundles, dual vector bundles, direct sums, tensor products and exterior products of vector bundles, and so on, all work as usual.

We can also define {\it r-smooth sections\/} of $E$, that is, an r-smooth map $s:X\ra E$ with $\pi\ci s=\id_X$. Note that r-smooth sections need not be a-smooth.

\label{ac4def3}
\end{dfn}

\begin{dfn} Let $X$ be an $m$-manifold with a-corners. The {\it b-tangent bundle\/} $\pi:{}^bTX\ra X$ is a natural rank $m$ vector bundle on $X$. We may describe ${}^bTX$ in local coordinates as follows. 

If $(U,\phi)$ is an a-chart on $X$, with $U\subseteq\R^{k,m}$ open, and $(x_1,\ldots,x_m)$ are the coordinates on $U$, then over $\phi(U)$, ${}^bTX$ is the trivial vector bundle with basis of sections $x_1\frac{\pd}{\pd x_1},\ldots,x_k\frac{\pd}{\pd x_k},\frac{\pd}{\pd x_{k+1}},\ldots,\frac{\pd}{\pd x_m}$. We have a corresponding a-chart $({}^bTU,{}^bT\phi)$ on ${}^bTX$, where ${}^bTU=U\t\R^m\subseteq\R^{k,2m}$, such that $(x_1,\ldots,x_m,q_1,\ab\ldots,\ab q_m)$ in ${}^bTU$ represents the vector 
\begin{equation*}
\ts q_1x_1\frac{\pd}{\pd x_1}+\cdots+q_kx_k\frac{\pd}{\pd x_k}+q_{k+1}\frac{\pd}{\pd x_{k+1}}+\cdots+q_m\frac{\pd}{\pd x_m}
\end{equation*}
over $\phi(x_1,\ldots,x_m)$ in $X$.

Let $(\ti U,\ti\phi)$ be another a-chart on $X$, with $\ti U\subseteq\R^{l,m}$ open, and $(\ti x_1,\ldots,\ti x_m)$ be the coordinates on $\ti U$. There is a change of coordinates $(x_1,\ab\ldots,\ab x_m)\ab\rightsquigarrow(\ti x_1,\ldots,\ti x_m)$ from $\dot U=\phi^{-1}(\ti\phi(\ti U))\subseteq U\subseteq\R^{k,m}$ to $\ddot U=\ti\phi^{-1}(\phi(U))\subseteq\ti U\subseteq\R^{l,m}$. That is, we may regard $(\ti x_1,\ldots,\ti x_m):\dot U\ra\ddot U$ as an a-diffeomorphism, with $\ti x_j=\ti x_j(x_1,\ldots,x_m)$, where $\ti x_j:\dot U\ra\lb 0,\iy)$ is a-smooth for $j=1,\ldots,l$, and $\ti x_j:\dot U\ra\R$ is a-smooth for $j=l+1,\ldots,m$.

We have another a-chart $({}^bT\ti U,{}^bT\ti \phi)$ on ${}^bTX$ corresponding to $(\ti U,\ti\phi)$, with coordinates $(\ti x_1,\ldots,\ti x_m,\ti q_1,\ab\ldots,\ab\ti q_m)$. Thus there should be a change of coordinates $(x_1,\ab\ldots,\ab x_m,\ab q_1,\ab\ldots,\ab q_m)\ab\rightsquigarrow(\ti x_1,\ldots,\ti x_m,\ti q_1,\ab\ldots,\ab\ti q_m)$ from $\dot U\t\R^m\subseteq{}^bTU$ to $\ddot U\t\R^m\subseteq{}^bT\ti U$. We define this to be:
\e
\begin{split}
\ti q_j=\begin{cases} \sum_{i=1}^k\ti x_j^{-1}x_i\frac{\pd\ti x_j}{\pd x_i}\,q_i+\sum_{i=k+1}^m\ti x_j^{-1}\frac{\pd\ti x_j}{\pd x_i}\,q_i, & j\le l, \\[4pt]
\sum_{i=1}^kx_i\frac{\pd\ti x_j}{\pd x_i}\,q_i+\sum_{i=k+1}^m\frac{\pd\ti x_j}{\pd x_i}\,q_i, & j>l. \end{cases}
\end{split}
\label{ac4eq8}
\e
This comes from $\frac{\pd}{\pd x_i}=\sum_{j=1}^m\frac{\pd\ti x_j}{\pd x_i}\cdot\frac{\pd}{\pd\ti x_j}$, applied to the bases $x_1\frac{\pd}{\pd x_1},\ab\ldots,\ab x_k\frac{\pd}{\pd x_k},\ab\frac{\pd}{\pd x_{k+1}},\ldots,\frac{\pd}{\pd x_m}$ for ${}^bTU$ and 
$\ti x_1\frac{\pd}{\pd\ti x_1},\ldots,\ti x_l\frac{\pd}{\pd\ti x_l},\frac{\pd}{\pd\ti x_{l+1}},\ldots,\frac{\pd}{\pd\ti x_m}$ for ${}^bT\ti U$.

The important thing about \eq{ac4eq8} is that the functions appearing in it are all a-smooth, as we will show. Let $j=1,\ldots,l$. Then $\ti x_j:\dot U\ra\lb 0,\iy)$ is a-smooth and interior, so by Definition \ref{ac3def3}(b)(i), locally in $\dot U$ we may write
\e
\ti x_j(x_1,\ldots,x_m)=F_j(x_1,\ldots,x_m)\cdot x_1^{a_{1,j}}\cdots x_k^{a_{k,j}},
\label{ac4eq9}
\e
where $a_{i,j}\in[0,\iy)$ and $F_j:\dot U\ra(0,\iy)$ is a-smooth. From \eq{ac4eq9} we see that
\e
\begin{aligned}
\ts\ti x_j^{-1}x_i\frac{\pd\ti x_j}{\pd x_i}&=\ts \bigl[F_j^{-1}\bigr]\cdot \bigl[x_i\frac{\pd F_j}{\pd x_i}\bigr]+a_{i,j}, 
& i&=1,\ldots,k,\\
\ts\ti x_j^{-1}\frac{\pd\ti x_j}{\pd x_i}&=\ts \bigl[F_j^{-1}\bigr]\cdot \bigl[\frac{\pd F_j}{\pd x_i}\bigr], & i&=k+1,\ldots,m.
\end{aligned}
\label{ac4eq10}
\e

Each term $[\cdots]$ in \eq{ac4eq10} is an a-smooth map to $\R$, by Definition \ref{ac3def2}. Thus the top line of \eq{ac4eq8} is an a-smooth function of $(x_1,\ldots,x_m,q_1,\ldots,q_m)$. The bottom line is clearly a-smooth. Hence the change of coordinates $(x_1,\ab\ldots,\ab x_m,\ab q_1,\ab\ldots,\ab q_m)\ab\rightsquigarrow(\ti x_1,\ldots,\ti x_m,\ti q_1,\ab\ldots,\ab\ti q_m)$ is a-smooth, with a-smooth inverse.

These coordinate changes between charts on ${}^bTX$ compose correctly under composition of coordinate changes on $X$. Thus they define a manifold with a-corners structure on ${}^bTX$, with a-smooth projection $\pi:{}^bTX\ra X$ acting in coordinates by $(x_1,\ldots,x_m, q_1,\ldots,q_m)\mapsto(x_1,\ldots,x_m)$, making ${}^bTX$ into a rank $m$ vector bundle over $X$. Elements of $\Ga^\iy({}^bTX)$ are called {\it b-vector fields}. 

There is a natural bilinear map ${}^b[\,,\,]:\Ga^\iy({}^bTX)\t \Ga^\iy({}^bTX)\ra \Ga^\iy({}^bTX)$ called the {\it b-Lie bracket}. If $u,v\in \Ga^\iy({}^bTX)$ and $w={}^b[u,v]$, and $(x_1,\ldots,x_m)\in\R^{k,m}$ are local coordinates on $\phi(U)\subseteq X$ as above, and on $\phi(U)$ we have
\begin{align*}
u\vert_{\phi(U)}&=\ts u_1x_1\frac{\pd}{\pd x_1}+\cdots+u_kx_k\frac{\pd}{\pd x_k}+u_{k+1}\frac{\pd}{\pd x_{k+1}}+\cdots+u_m\frac{\pd}{\pd x_m},\\
v\vert_{\phi(U)}&=\ts v_1x_1\frac{\pd}{\pd x_1}+\cdots+v_kx_k\frac{\pd}{\pd x_k}+v_{k+1}\frac{\pd}{\pd x_{k+1}}+\cdots+v_m\frac{\pd}{\pd x_m},\\
w\vert_{\phi(U)}&=\ts w_1x_1\frac{\pd}{\pd x_1}+\cdots+w_kx_k\frac{\pd}{\pd x_k}+w_{k+1}\frac{\pd}{\pd x_{k+1}}+\cdots+w_m\frac{\pd}{\pd x_m},
\end{align*}
for $u_i,v_i,w_i:\phi(U)\ra\R$ a-smooth, then $w_1,\ldots,w_m$ are given by
\begin{equation*}
w_i=\ts\sum_{j=1}^k\bigl(u_jx_j\frac{\pd v_i}{\pd x_j}-v_jx_j\frac{\pd u_i}{\pd x_j}\bigr)+
\sum_{j=k+1}^m\bigl(u_j\frac{\pd v_i}{\pd x_j}-v_j\frac{\pd u_i}{\pd x_j}\bigr).\end{equation*}
The b-Lie bracket has the usual properties of Lie brackets on manifolds, in particular ${}^b[u,v]={}^b[-v,u]$ and ${}^b[u,{}^b[v,w]]+{}^b[v,{}^b[w,u]]+{}^b[w,{}^b[u,v]]=0$ for all $u,v,w\in\Ga^\iy({}^bTX)$. If $X\in\Man\subset\Manac$ then ${}^bTX=TX$ and~${}^b[\,,\,]=[\,,\,]$.

The {\it b-cotangent bundle\/} ${}^bT^*X$ is the dual vector bundle of ${}^bTX$. If $(U,\phi)$ is a chart on $X$, with $U\subseteq\R^{k,m}$ open, and $(x_1,\ldots,x_m)$ are the coordinates on $U$, then ${}^bT^*X$ has basis of sections $x_1^{-1}\d x_1,\ldots,x_k^{-1}\d x_k,\d x_{k+1},\ldots,\d x_m$ over~$\phi(U)$.

\label{ac4def4}
\end{dfn}

B-tangent bundles are functorial under interior maps:

\begin{dfn} Suppose $f:X\ra Y$ is an interior map of manifolds with a-corners. Then there is a natural interior map ${}^bTf:{}^bTX\ra{}^bTY$. To define ${}^bTf$, let $(U,\phi)$ and $(V,\psi)$ be a-charts on $X,Y$ with $U\subseteq\R^{k,m}$, $V\subseteq\R^{l,n}$, with coordinates $(x_1,\ldots,x_m)\in U$ and $(y_1,\ldots,y_n)\in V$, and let $({}^bTU,{}^bT\phi)$, $({}^bTV,{}^bT\psi)$ be the corresponding a-charts on ${}^bTX,{}^bTY$, with coordinates $(x_1,\ab\ldots,\ab x_m,\ab q_1,\ab\ldots,\ab q_m)\in {}^bTU$ and $(y_1,\ldots,y_n,r_1,\ldots,r_n)\in {}^bTV$. Then \eq{ac3eq6} defines $\psi^{-1}\ci f\ci\phi:\dot U\ra V$, where $\dot U=(f\ci\phi)^{-1}(\psi(V))\subseteq U$ is open.

Write $\psi^{-1}\ci f\ci\phi=(f_1,\ldots,f_n)$, for $f_j=f_j(x_1,\ldots,x_m)$. Then the corresponding ${}^bT\psi^{-1}\ci {}^bTf\ci {}^bT\phi$ maps
\e
\begin{split}
&{}^bT\psi^{-1}\!\ci\! {}^bTf\!\ci\! {}^bT\phi:(x_1,\ldots,x_m,q_1,\ldots,q_m)\!\longmapsto\! (y_1,\ldots,y_n,r_1,\ldots,r_n),\!\!\!\!\!\!\!\!{}
\\
&\text{where}\quad y_j=f_j(x_1,\ldots,x_m),\quad j=1\ldots,n,\\
&\text{and}\quad r_j=
\begin{cases} \sum_{i=1}^kf_j^{-1}x_i\frac{\pd f_j}{\pd x_i}\,q_i+\sum_{i=k+1}^mf_j^{-1}\frac{\pd f_j}{\pd x_i}\,q_i, & j\le l, \\[4pt]
\sum_{i=1}^kx_i\frac{\pd f_j}{\pd x_i}\,q_i+\sum_{i=k+1}^m\frac{\pd f_j}{\pd x_i}\,q_i, & j>l. \end{cases}
\end{split}
\label{ac4eq11}
\e

As for \eq{ac4eq8}--\eq{ac4eq10}, we claim that the functions $f_j^{-1}x_i\frac{\pd f_j}{\pd x_i}$ for $i\le k$, $j\le l$ and $f_j^{-1}\frac{\pd f_j}{\pd x_i}$ for $i>k$, $j\le l$ occurring in \eq{ac4eq11} extend uniquely to a-smooth functions of $(x_1,\ldots,x_m)$ in $\dot U$ where $f_j=0$. To see this, note that as $f$ is interior, by Definition \ref{ac3def3}(b)(i) for $j=1,\ldots,l$, locally in $\dot U$ we may write
\e
f_j(x_1,\ldots,x_m)=F_j(x_1,\ldots,x_m)\cdot x_1^{a_{1,j}}\cdots x_k^{a_{k,j}},
\label{ac4eq12}
\e
where $a_{i,j}\in[0,\iy)$ and $F_j:\dot U\ra(0,\iy)$ is a-smooth. From \eq{ac4eq12} we see that
\e
\begin{aligned}
\ts f_j^{-1}x_i\frac{\pd f_j}{\pd x_i}&=\ts \bigl[F_j^{-1}\bigr]\cdot \bigl[x_i\frac{\pd F_j}{\pd x_i}\bigr]+a_{i,j}, 
& i&=1,\ldots,k,\\
\ts f_j^{-1}\frac{\pd f_j}{\pd x_i}&=\ts \bigl[F_j^{-1}\bigr]\cdot \bigl[\frac{\pd F_j}{\pd x_i}\bigr], & i&=k+1,\ldots,m.
\end{aligned}
\label{ac4eq13}
\e

Each term $[\cdots]$ in \eq{ac4eq13} is an a-smooth map to $\R$, by Definition \ref{ac3def2}. Thus the third line of \eq{ac4eq11} is an a-smooth function of $(x_1,\ldots,x_m,q_1,\ldots,q_m)$. The fourth line is clearly a-smooth. Hence \eq{ac4eq11} defines an interior map
\begin{equation*}
{}^bT\psi^{-1}\ci {}^bTf\ci {}^bT\phi:{}^bT\dot U\longra{}^bTV.
\end{equation*}
These maps are compatible with changes of a-charts $(U,\phi)$ on $X$, and $(V,\psi)$ on $Y$. Thus they define a global interior map ${}^bTf:{}^bTX\ra{}^bTY$, which fits into a commutative diagram:
\begin{equation*}
\xymatrix@C=80pt@R=15pt{ *+[r]{{}^bTX} \ar[d]^\pi \ar[r]_{{}^bTf} &
*+[l]{{}^bTY} \ar[d]_\pi \\ *+[r]{X} \ar[r]^f & *+[l]{Y.} }
\end{equation*}
We can also regard ${}^bTf$ as a vector bundle morphism ${}^b\d f:{}^bTX\ra f^*({}^bTY)$ on $X$, which has dual morphism ${}^b\d f:f^*({}^bT^*Y)\ra {}^bT^*X$.

If $g:Y\ra Z$ is another interior map of manifolds with a-corners, it is easy to check from \eq{ac4eq11} that ${}^bT(g\ci f)={}^bTg\ci{}^bTf:{}^bTX\ra {}^bTZ$. Thus the assignment $X\mapsto {}^bTX$, $f\mapsto {}^bTf$ is a functor, the {\it b-tangent functor\/} ${}^bT:\Manacin\ra\Manacin$. It extends to ${}^bT:\cManacin\ra\cManacin$ in the obvious way.

Note that if $f:X\ra Y$ is a morphism in $\Manac$ then $C(f):C(X)\ra C(Y)$ is interior in $\cManac$, so ${}^bTC(f):{}^bTC(X)\ra{}^bTC(Y)$ is well defined, and we can use this as a substitute for ${}^bTf:{}^bTX\ra{}^bTY$ when $f$ is not interior.
\label{ac4def5}
\end{dfn}

\begin{rem}{\bf(a)} We cannot define ordinary tangent bundles $TX$ for manifolds with a-corners. This is because if $(U,\phi),(\ti U,\ti\phi)$ are a-charts on $X$ with $U\subseteq\R^{k,m}$ open, $\ti U\subseteq\R^{l,m}$ open, and we tried to define a-charts $(TU,T\phi)$ with coordinates $(x_1,\ab\ldots,\ab x_m,\ab q_1,\ab\ldots,\ab q_m)$ on $TX$ and $(T\ti U,T\ti\phi)$ on $TX$ with coordinates $(\ti x_1,\ldots,\ti x_m,\ti q_1,\ab\ldots,\ab\ti q_m)$, then in the analogue of \eq{ac4eq8}, the change of coordinates $(x_1,\ab\ldots,\ab x_m,\ab q_1,\ab\ldots,\ab q_m)\ab\rightsquigarrow(\ti x_1,\ldots,\ti x_m,\ti q_1,\ab\ldots,\ab\ti q_m)$ should be
\begin{equation*}
\ti q_j=\ts\sum_{i=1}^m\frac{\pd\ti x_j}{\pd x_i}\,q_i. 
\end{equation*}
But when $i=1,\ldots,k$, the function $\frac{\pd\ti x_j}{\pd x_i}$ need not be a-smooth at $x_i=0$, so this change of coordinates is not a-smooth, and does not define an a-atlas on $TX$.
\smallskip

\noindent{\bf(b)} The reason that we only define ${}^bTf:{}^bTX\ra{}^bTY$ for interior $f$ is that if $f$ were not interior, we could have $f_j(x_1,\ldots,x_m)=0$ instead of \eq{ac4eq12}, and then the terms $f_j^{-1}x_i\frac{\pd f_j}{\pd x_i}$, $f_j^{-1}\frac{\pd f_j}{\pd x_i}$ in \eq{ac4eq11} do not make sense.
\smallskip

\noindent{\bf(c)} The definitions of r-smooth and a-smooth functions in \S\ref{ac21} were designed so that b-(co)tangent bundles, and b-derivatives below, would be well-behaved.

\label{ac4rem2}
\end{rem}

Here is the analogue of Definition \ref{ac2def8}:

\begin{dfn} Let $f:X\ra Y$ be an interior map of manifolds with a-corners. We call $f$ a {\it b-submersion\/} if ${}^b\d f:{}^bTX\ra f^*({}^bTY)$ is a surjective morphism of vector bundles on $X$. We call $f$ a {\it b-fibration\/} if $f$ is b-normal, in the sense of Definitions \ref{ac3def3}(d) and \ref{ac3def4}, and a b-submersion.
\label{ac4def6}
\end{dfn}

B-fibrations $f:X\ra Y$ are important as they are a good notion of `family of manifolds with a-corners $X_y$ over a base $Y$'. We will discuss them further in \S\ref{ac54} and \S\ref{ac62}. We can give examples of b-fibrations as in Example~\ref{ac2ex3}.

Here is the analogue of Definition~\ref{ac2def9}:

\begin{dfn} Let $X$ be a manifold with a-corners. As in \eq{ac2eq11}, there is an exact sequence 
\e
\xymatrix@C=22pt{ 0 \ar[r] & {}^bN_{\pd X} \ar[rr]^(0.4){{}^bi_T} &&
i_X^*({}^bTX) \ar[rr]^{{}^b\pi_T} && {}^bT(\pd X) \ar[r] & 0 }
\label{ac4eq14}
\e
of vector bundles on $\pd X$, where ${}^bN_{\pd X}$ is a line bundle on $\pd X$, defined as a vector subbundle of $i_X^*({}^bTX)$, which we call the {\it b-normal bundle\/} of $\pd X$ in $X$, and ${}^bi_T:{}^bN_{\pd X}\hookra i_X^*({}^bTX)$ is the inclusion.

We define ${}^bN_{\pd X}$ and ${}^b\pi_T$ in \eq{ac4eq14} as follows: if $(U,\phi)$ is an a-chart on $X$, with $U\subseteq\R^{k,m}$ open for $k\ge 1$, and $(x_1,\ldots,x_m)$ are the coordinates on $U$, and $(x,\be)\in\pd X$ with $x=\phi(\ti x_1,\ldots,\ti x_m)$ with $\ti x_1=0$, and $\be$ is the local boundary component $x_1=0$ of $X$ at $x$, then ${}^bN_{\pd X}\vert_{(x,\be)}=\langle x_1\frac{\pd}{\pd x_1}\rangle_\R$, and ${}^b\pi_T$ maps 
\begin{align*}
{}^b\pi_T&:\ts\bigl((x,\be),\sum_{i=1}^mc_i\cdot x_i\frac{\pd}{\pd x_i}\bigr)\longmapsto\bigl((x,\be),\sum_{i=2}^mc_i\cdot x_i\frac{\pd}{\pd x_i}\bigr),
\end{align*}
for $c_i\in\R$, using $(x_2,\ldots,x_m)\in\R^{k-1,m-1}$ as the local coordinates on $\pd X$ near~$x$.

The b-normal bundle ${}^bN_{\pd X}$ has a natural orientation. We say that a section $\la\in\Ga^\iy({}^bN_{\pd X})$ is {\it positive}, written $\la>0$, if in local coordinates $(x_2,\ldots,x_m)$ on $\pd X$ as above we have $\la=c(x_2,\ldots,x_m)\cdot x_1\frac{\pd}{\pd x_1}$ for a local smooth function $c:\pd X\ra\R$ with $c>0$. Similarly we define when $\la$ is {\it negative}, {\it nonnegative}, and {\it nonpositive}, written $\la<0$, $\la\ge 0$, $\la\le 0$ respectively, in the obvious way.
\label{ac4def7}
\end{dfn}

\begin{rem}{\bf(a)} For manifolds with corners $X$ in Definition \ref{ac2def9}, we defined ${}^bN_{\pd X}$ to be the trivial line bundle $\pd X\t\R\ra\pd X$. But for manifolds with a-corners, ${}^bN_{\pd X}$ is {\it not canonically trivial}, although it is non-canonically isomorphic to the trivial line bundle. To see this, note that if $X,(U,\phi)$ are as in Definition \ref{ac4def7}, and $\al>0$, then we can define another a-chart $(\hat U,\hat\phi)$ on $X$ with
\begin{align*}
\hat U&=\bigl\{(x_1^{1/\al},x_2,\ldots,x_m):(x_1,\ldots,x_m)\in U\bigr\}\subseteq\R^{k,m},\\
\hat\phi(\hat x_1,\ldots,\hat x_m)&=\phi(\hat x_1^\al,\hat x_2,\ldots,\hat x_m).
\end{align*}
Then the a-charts $(U,\phi),(\hat U,\hat\phi)$ give the bases $x_1\frac{\pd}{\pd x_1}$ and $\hat x_1\frac{\pd}{\pd\hat x_1}$ over $\phi(U)=\hat\phi(\hat U)$, where $\hat x_1\frac{\pd}{\pd\hat x_1}=\al\cdot x_1\frac{\pd}{\pd x_1}$ as $x_1=\hat x_1^\al$. Hence if $\al\ne 1$ then $(U,\phi),(\hat U,\hat\phi)$ give different trivializations of ${}^bN_{\pd X}$, and there is no canonical trivialization.
\smallskip

\noindent{\bf(b)} The material of this section extends to $\Mancac$ in~\S\ref{ac35}. For a manifold $X$ with corners and a-corners, we can define two natural tangent bundles: the {\it b-tangent bundle\/} ${}^bTX$, and the {\it mixed tangent bundle\/} ${}^mTX$. Here if $X\in\Manc\subset\Mancac$ then ${}^bTX$ is as in \S\ref{ac23} and ${}^mTX=TX$ from \S\ref{ac23}, and if $X\in\Manac\subset\Mancac$ then ${}^bTX,{}^mTX$ are both ${}^bTX$ above. If $(U,\phi)$ is a chart on $X$, with $U\subseteq\R^{k,m}_l$ open, and $(x_1,\ldots,x_m)$ are the coordinates on $U$, then over $\phi(U)$, ${}^bTX$ has basis $x_1\frac{\pd}{\pd x_1},\ldots,x_{k+l}\frac{\pd}{\pd x_{k+l}},\frac{\pd}{\pd x_{k+l+1}},\ldots,\frac{\pd}{\pd x_m}$, and ${}^mTX$ has basis~$x_1\frac{\pd}{\pd x_1},\ldots,x_k\frac{\pd}{\pd x_k},\frac{\pd}{\pd x_{k+1}},\ldots,\frac{\pd}{\pd x_m}$. Here ${}^bTX$ is functorial over interior maps in $\Mancac$, and ${}^mTX$ is functorial over smooth maps which are interior in the a-corner variables. We generally prefer to use ${}^mTX,{}^mT^*X$.
\label{ac4rem3}
\end{rem}

\subsection{The b-derivative, and b-de Rham cohomology}
\label{ac43}

Next we generalize exterior forms, the de Rham differential, and de Rham cohomology to manifolds with a-corners.

\begin{dfn} Let $X$ be a manifold with a-corners, of dimension $m$. We will define the {\it b-de Rham differential\/} ${}^b\d:C^\iy(X)\ra\Ga^\iy({}^bT^*X)$, an $\R$-linear map. Let $(U,\phi)$ be an a-chart on $X$, with $U\subseteq\R^{k,m}$ open, and $(x_1,\ldots,x_m)$ be the coordinates on $U$. Then as for ${}^bTX$ in Definition \ref{ac4def4}, we have an a-chart $({}^bT^*U,{}^bT^*\phi)$ on ${}^bT^*X$, where ${}^bT^*U=U\t\R^m\subseteq\R^{k,2m}$, such that $(x_1,\ldots,x_m,q_1,\ab\ldots,\ab q_m)$ in ${}^bTU$ represents the b-covector
\e
\ts q_1(x_1^{-1}\d x_1)+\cdots+q_k(x_k^{-1}\d x_k)+q_{k+1}\d x_{k+1}+\cdots+q_m\d x_m
\label{ac4eq15}
\e
in ${}^bT^*_{\phi(x_1,\ldots,x_m)}X$. Here $x_1^{-1}\d x_1,\ldots,x_k^{-1}\d x_k,\d x_{k+1},\ldots,\d x_m$ are just notation for the basis of sections of ${}^bT^*X$ over $\phi(U)$ dual to the basis of sections $x_1\frac{\pd}{\pd x_1},\ab\ldots,\ab x_k\frac{\pd}{\pd x_k},\ab\frac{\pd}{\pd x_{k+1}},\ldots,\frac{\pd}{\pd x_m}$ for ${}^bTX$ over $\phi(U)$ (which are also just notation). In particular, $x_i^{-1}\d x_i$ are well defined for $i=1,\ldots,k$ even when~$x_i=0$.

Then ${}^b\d:C^\iy(X)\ra\Ga^\iy({}^bT^*X)$ for $c\in C^\iy(X)$ acts in local coordinates by
\begin{align*}
{}^b\d:c\longmapsto\ts \bigl(x_1\frac{\pd c}{\pd x_1}\bigr)\cdot \bigl(x_1^{-1}\d x_1\bigr)+\cdots+\bigl(x_k\frac{\pd c}{\pd x_k}\bigr)\cdot \bigl(x_k^{-1}\d x_k\bigr)&\\
\ts+\bigl(\frac{\pd c}{\pd x_{k+1}}\bigr)\cdot\bigl(\d x_{k+1}\bigr)+\cdots+\bigl(\frac{\pd c}{\pd x_m}\bigr)\cdot\bigl(\d x_m\bigr)&.
\end{align*}
That is, in local coordinates $(x_1,\ldots,x_m)$ on $X$ and $(x_1,\ldots,x_m,q_1,\ldots,q_m)$ representing the b-covector \eq{ac4eq15} in ${}^bT^*X$ as above, ${}^b\d c:X\ra {}^bT^*X$ acts by
\e
{}^b\d c:(x_1,\ldots,x_m)\longmapsto \ts\bigl(x_1,\ldots,x_m,x_1\frac{\pd c}{\pd x_1},\ldots,x_k\frac{\pd c}{\pd x_k},\frac{\pd c}{\pd x_{k+1}},\ldots,\frac{\pd c}{\pd x_m}\bigr).
\label{ac4eq16}
\e
These ${}^b\d c$ in \eq{ac4eq16} are a-smooth by \S\ref{ac31}, and are compatible with changes of a-charts $(U,\phi)$ on $X$ and $({}^bT^*U,{}^bT^*\phi)$ on ${}^bT^*X$. Thus they define a global a-smooth map ${}^b\d c:X\ra {}^bT^*X$, with $\pi\ci{}^b\d c=\id_X:X\ra X$, so ${}^b\d c\in\Ga^\iy({}^bT^*X)$.

The interior $X^\ci$ is an ordinary manifold, with ${}^bT^*X\vert_{X^\ci}=T^*X^\ci$, and ${}^b\d$ restricts to the usual de Rham differential $\d:C^\iy(X^\ci)\ra C^\iy(T^*X^\ci)$ on $X^\ci$.

From \eq{ac4eq16} we see that ${}^b\d$ satisfies
\begin{equation*}
{}^b\d(c_1c_2)=c_1\cdot{}^b\d c_2+c_2\cdot{}^b\d c_1
\end{equation*}
for all $c_1,c_2\in C^\iy(X)$, that is, ${}^b\d$ is a derivation on the $\R$-algebra~$C^\iy(X)$.

As for exterior derivatives on manifolds, we can generalize ${}^b\d$ to unique $\R$-linear maps ${}^b\d:C^\iy(\La^l({}^bT^*X))\ra\Ga^\iy(\La^{l+1}({}^bT^*X))$ for $l\ge 0$ called {\it b-exterior derivatives}, given by \eq{ac4eq16} when $l=0$, such that ${}^b\d\ci{}^b\d=0$ and ${}^b\d(\al\w\be)=({}^b\d\al)\w\be+(-1)^k\al\w({}^b\d\be)$ for all $\al\in C^\iy(\La^k({}^bT^*X))$ and~$\be\in C^\iy(\La^l({}^bT^*X))$.

As for ordinary de Rham cohomology, as in Bredon \cite[\S V]{Bred1}, define the {\it b-de Rham cohomology groups\/} of $X$ for $l=0,1,\ldots$ to be
\begin{equation*}
{}^bH^l_{\rm dR}(X;\R)=\frac{\Ker\bigl({}^b\d:C^\iy(\La^l({}^bT^*X))\ra\Ga^\iy(\La^{l+1}({}^bT^*X))\bigr)}{\Im\bigl({}^b\d:C^\iy(\La^{l-1}({}^bT^*X))\ra\Ga^\iy(\La^l({}^bT^*X))\bigr)}\,.
\end{equation*}
If $X\in\Man\subset\Manac$ then ${}^bH^l_{\rm dR}(X;\R)=H^l_{\rm dR}(X;\R)$ is ordinary de Rham cohomology, and so isomorphic to the (topological) cohomology of $X$ over~$\R$.
\label{ac4def8}
\end{dfn}

\begin{ex} Let $X=\lb 0,1\rb$. There is an isomorphism $C^\iy(X)\cong\Ga^\iy({}^bT^*X)$ identifying $c'\in C^\iy(X)$ with $c'\cdot x^{-1}(1-x)^{-1}\d x$, with $x$ the coordinate on $\lb 0,1\rb$. Using this we see that
\e
\begin{split}
{}^bH^0_{\rm dR}(X;\R)&=\ts\bigl\{c\in C^\iy(X):x(1-x)\frac{\pd c}{\pd x}=0\bigr\}=\langle 1\rangle_\R\cong\R,\\
{}^bH^1_{\rm dR}(X;\R)&\cong \frac{C^\iy(X)}{\bigl\{x(1-x)\frac{\pd c}{\pd x}:c\in C^\iy(X)\bigr\}}\,.
\end{split}
\label{ac4eq17}
\e

Let $c'\in C^\iy(X)$. If $c'=x(1-x)\frac{\pd c}{\pd x}$ for some $c\in C^\iy(X)$ then $c'(0)=c'(1)=0$. Conversely, suppose $c'(0)=c'(1)=0$, and define $c:\lb 0,1\rb\ra\R$ by
\e
c(x)=\int_{1/2}^x t^{-1}(1-t)^{-1}c'(t)\,\d t.
\label{ac4eq18}
\e
Since $c'$ is a-smooth with $c'(0)=c'(1)=0$, from \S\ref{ac31} there exist $C,\al>0$ and $\ep\in[0,1)$ such that $\md{c'(t)}\le Ct^\al$ for $t\in[0,\ep)$ and $\md{c'(t)}\le C(1-t)^\al$ for $t\in(1-\ep,1]$. It follows that $c$ in \eq{ac4eq18} is well defined, including at $x=0,1$, with $c(x)-c(0)=O(x^\al)$ as $x\ra 0$ and $c(x)-c(1)=O((1-x)^\al)$ as $x\ra 1$. Also $x(1-x)\frac{\pd c}{\pd x}=c'$, which is a-smooth, so $c\in C^\iy(X)$. Therefore \eq{ac4eq17} gives
\begin{equation*}
{}^bH^0_{\rm dR}(\lb 0,1\rb;\R)\cong\R\qquad\text{and}\qquad {}^bH^1_{\rm dR}(\lb 0,1\rb;\R)\cong\R^2,
\end{equation*}
where the second equation identifies $\bigl[c'\cdot x^{-1}(1-x)^{-1}\d x\bigr]$ with~$(c'(0),c'(1))\in\R^2$.

\label{ac4ex1}
\end{ex}

\begin{rem}{\bf(a)} In \S\ref{ac31}--\S\ref{ac32} we chose to define manifolds with a-corners using a-smooth functions, including the $O(x^\al)$ decay condition \eq{ac3eq2} near boundary hypersurfaces $x=0$, rather than the weaker r-smooth functions.

If we had used r-smooth functions, then $c:\lb 0,1\rb\ra\R$ in \eq{ac4eq18} need not be defined at $x=0,1$, and in fact ${}^bH^1_{\rm dR}(\lb 0,1\rb;\R)$ would be infinite-dimensional. A-smoothness is essential for b-de Rham cohomology to be well behaved.
\smallskip

\noindent{\bf(b)} Consider b-de Rham cohomology from the point of view of sheaf theory, as in Bredon \cite{Bred2}. Let $X$ be a manifold with a-corners, and write ${}^b\Om^k$ for the sheaf of a-smooth sections of $\La^k({}^bT^*X)$ for $k=0,1,\ldots,$ and ${}^b\d:{}^b\Om^k\ra{}^b\Om^{k+1}$ for the sheaf morphism induced by b-exterior derivatives. 

Then ${}^b\Om^k$ is a soft sheaf of $\R$-vector spaces on $X$, since partitions of unity exist in $C^\iy(X)$, and ${}^b\Om^\bu=({}^b\Om^*,{}^b\d)$ is a complex of soft sheaves on $X$, regarded as an object of the derived category $D\mathop{\rm Sh}_\R(X)$ of sheaves of $\R$-vector spaces on $X$. The hypercohomology ${\mathbb H}^l({}^b\Om^\bu)$ is isomorphic to ${}^bH^l_{\rm dR}(X;\R)$, as ${}^b\Om^k$ is soft.

In the analogue for ordinary manifolds $X$ we would have $\Om^\bu\simeq\R_X$ in $D\mathop{\rm Sh}_\R(X)$, where $\R_X$ is the constant sheaf on $X$ with fibre $\R$, so that
\begin{equation*}
H^l_{\rm dR}(X;\R)\cong {\mathbb H}^l(\Om^\bu)\cong {\mathbb H}^l(\R_X)\cong H^l(X;\R),
\end{equation*}
that is, de Rham cohomology is isomorphic to ordinary cohomology.

In the b-de Rham case we will {\it not\/} have ${}^b\Om^\bu\simeq\R_X$. However, we should be able to describe the cohomology sheaves $h^k({}^b\Om^\bu)$ explicitly. The author expects that $h^0({}^b\Om^\bu)\cong\R_X$, and $h^k({}^b\Om^\bu)$ for $k\ge 0$ is supported on $\ov{S^k(X)}\subset X$, and is the pushforward $\Pi_*(\cS_{C_k(X)})$ along $\Pi:C_k(X)\ra X$ of a locally constant sheaf $\cS_{C_k(X)}$ on $C_k(X)$ with fibre $\R$ (that is, a flat line bundle on $C_k(X)$). When $k=1$, this is the dual line bundle ${}^bN_{\pd X}^*$ in \S\ref{ac42}, with flat connection from \S\ref{ac44}. There should be a spectral sequence $\bigl(H^j(C_k(X),\cS_{C_k(X)})\bigr)_{j,k}\Longra {}^bH^{j+k}_{\rm dR}(X;\R)$ computing b-de Rham cohomology.

\smallskip

\noindent{\bf(c)} Part {\bf(b)} would imply that ${}^bH^*(X;\R)$ is finite-dimensional if $X$ is compact, or if $X$ is noncompact but `topologically finite' in a suitable sense (e.g.\ $X$ should not have infinitely many noncompact ends).
\smallskip

\noindent{\bf(d)} In the case when $X$ is a manifold with a-boundary, the prediction in {\bf(b)} yields a long exact sequence
\begin{equation*}
\xymatrix@C=12pt{ \cdots \ar[r] & H^l(X;\R) \ar[r] & {}^bH^l_{\rm dR}(X;\R) \ar[r] & H^{l-1}(\pd X,{}^bN_{\pd X}^*) \ar[r] & H^{l+1}(X;\R) \ar[r] & \cdots, }
\end{equation*}
for $H^*(\pd X,{}^bN_{\pd X}^*)$ the cohomology of $\pd X$ twisted by the flat line bundle~${}^bN_{\pd X}^*$.

\smallskip

\noindent{\bf(e)} Let $X=(\lb 0,\iy)\t\R)/\Z$ depending on $\al>0$ be as in Example \ref{ac3ex3}. Then
\begin{equation*}
H^l(X;\R)=\begin{cases} \R, & l=0,1, \\ 0, & \text{otherwise,} \end{cases} \quad
H^l(\pd X,{}^bN_{\pd X}^*)=\begin{cases} \R, & \text{$l=0,1$ and $\al=1$,} \\ 0, & \text{otherwise.} \end{cases}
\end{equation*}
Thus {\bf(d)} predicts that if $\al=1$ then ${}^bH^l_{\rm dR}(X;\R)$ is $\R$ if $l=0,2$, and $\R^2$ if $l=1$, and 0 otherwise; but if $\al\ne 1$ then ${}^bH^l_{\rm dR}(X;\R)$ is $\R$ if $l=0,1$, and 0 otherwise. If this prediction is confirmed, it shows that ${}^bH^*_{\rm dR}(X;\R)$ {\it depends on the a-smooth structure of\/} $X$, and not just on the (stratified) topological space.

\label{ac4rem4}
\end{rem}

\subsection{B-connections on vector bundles}
\label{ac44}

Let $X$ be a manifold, and $E\ra X$ a vector bundle. One way to define a {\it connection\/} $\nabla$ on $E$ is an $\R$-linear map $\nabla:\Ga^\iy(E)\ra\Ga^\iy(E\ot T^*X)$ satisfying
\begin{equation*}
\nabla(c\cdot e)=c\cdot \nabla e+e\ot \d c\qquad\text{for all $c\in C^\iy(X)$ and $e\in\Ga^\iy(E)$.}
\end{equation*}
Here is the analogue for manifolds with a-corners.

\begin{dfn} Let $X$ be a manifold with a-corners, and $E\ra X$ a vector bundle. A {\it b-connection\/} on $E$ is an $\R$-linear map ${}^b\nabla:\Ga^\iy(E)\ra\Ga^\iy(E\ot {}^bT^*X)$ satisfying, for ${}^b\d c$ as in \S\ref{ac43},
\begin{equation*}
{}^b\nabla(c\cdot e)=c\cdot {}^b\nabla e+e\ot {}^b\d c\qquad\text{for all $c\in C^\iy(X)$ and $e\in\Ga^\iy(E)$.}
\end{equation*}

If $v\in \Ga^\iy({}^bTX)$ is a b-vector field, and $e\in\Ga^\iy(E)$, we write ${}^b\nabla_ve\in\Ga^\iy(E)$ for $v\cdot{}^b\nabla e$, where $v\ot{}^b\nabla e$ is a section of $({}^bTX)\ot(E\ot {}^bT^*X)$, and `$\cdot$' applies the dual pairing between ${}^bTX,{}^bT^*X$ to map~${}^bTX\ot E\ot {}^bT^*X\ra E$.

As for connections in ordinary differential geometry, a b-connection ${}^b\nabla$ has a natural {\it b-curvature\/} $R_{{}^b\nabla}\in \Ga^\iy\bigl(\End(E)\ot \La^2({}^bT^*X)\bigr)$, characterized by the property that if $v,w\in\Ga^\iy({}^bTX)$ and $e\in\Ga^\iy(E)$ then
\begin{equation*}
{}^b\nabla_v({}^b\nabla_we)-{}^b\nabla_w({}^b\nabla_ve)-{}^b\nabla_{{}^b[v,w]}e=R_{{}^b\nabla}\cdot (e\ot(v\w w)),
\end{equation*}
where ${}^b[v,w]$ is the b-Lie bracket from \S\ref{ac42}. We call ${}^b\nabla$ {\it flat\/} if $R_{{}^b\nabla}=0$.

On the interior $X^\ci$, which is an ordinary manifold, ${}^b\nabla$ is an ordinary connection on $E\vert_{X^\ci}$, and $R_{{}^b\nabla}$ is the usual curvature.
\label{ac4def9}
\end{dfn}

Many of the usual properties of connections extend to b-connections. For example, b-connections ${}^b\nabla^E,{}^b\nabla^F$ on vector bundles $E,F\ra X$ induce natural b-connections on all the associated vector bundles $E^*,E\op F,E\ot F,\La^kE,\ldots.$

\begin{ex}{\bf(a)} Let $\ul{\R}^k$ be the trivial vector bundle $X\t\R^k\ra X$ over a manifold with a-corners $X$. Then $\ul{\R}^k$ has a natural, flat b-connection ${}^b\nabla_0$ given by ${}^b\nabla_0(e_1,\ldots,e_k)=({}^b\d e_1,\ldots,{}^b\d e_k)$, for $e_1,\ldots,e_k\in C^\iy(X)$.

If ${}^b\nabla$ is a b-connection on $E\ra X$, we call $(E,{}^b\nabla)$ {\it locally trivial\/} if we can cover $X$ by open subsets on which $(E,{}^b\nabla)\cong(\ul{\R}^k,{}^b\nabla_0)$, for $k=\rank E$.
\smallskip

\noindent{\bf(b)} Let $X=\lb 0,\iy)$, with coordinate $x$, and $E$ be the trivial line bundle $X\t\R\ra X$. Define a b-connection ${}^b\nabla$ on $E$ by ${}^b\nabla e={}^b\d e+e\ot(x^{-1}\d x)$ for $e\in\Ga^\iy(E)=C^\iy(X)$. Then ${}^b\nabla$ is trivially flat, as $\La^2({}^bT^*X)=0$. However, near $0\in X$ there exists no nonzero section $e\in \Ga^\iy(E)$ with ${}^b\nabla e=0$, as solving the o.d.e.\ locally gives $e(x)=c x^{-1}$ for $x>0$ and $c\in\R$, which does not extend to $x=0$. So $E$ is not locally trivial near $0\in X$.
\smallskip

In conventional differential geometry, by the Frobenius Theorem a vector bundle with connection $(E,\nabla)$ is flat if and only if it is locally trivial. But {\bf(a)},{\bf(b)} show that locally trivial b-connections are flat, but {\it flat b-connections need not be locally trivial}. That is, {\it the Frobenius Theorem fails for b-connections}. 

\label{ac4ex2}
\end{ex}

\begin{dfn} Let $X$ be a manifold with a-corners. Then Definition \ref{ac4def7} defined a line bundle ${}^bN_{\pd X}$ on $\pd X$. We will define a natural flat b-connection ${}^b\nabla_{\pd X}$ on ${}^bN_{\pd X}$, which is locally trivial. Suppose $(U,\phi)$ is an a-chart on $X$, with $U\subseteq\R^{k,m}$ open for $k\ge 1$, and $(x_1,\ldots,x_m)$ are the coordinates on $U$. Then we have an a-chart $(U',\phi')$ on $\pd X$ with
\begin{align*}
U'&=\bigl\{(x_2,\ldots,x_m)\in\R^{k-1,m-1}:(0,x_2,\ldots,x_m)\in U\subseteq\R^m_k\bigr\},\\
\phi'(x_2,\ldots,x_m)&=\bigl(\phi(0,x_2,\ldots,x_m),\phi_*(\{x_1=0\})\bigr),
\end{align*}
and over $\phi'(U')$ we have ${}^bN_{\pd X}=\langle x_1\frac{\pd}{\pd x_1}\rangle_\R$. We require that ${}^b\nabla_{\pd X}(x_1\frac{\pd}{\pd x_1})=0$ over $\phi'(U')$, for all such $(U,\phi)$.

To see that this is well defined, we have to show that if $(\ti U,\ti\phi)$ is an a-chart on $X$, with $\ti U\subseteq\R^{l,m}$ open for $l\ge 1$, and $(\ti x_1,\ldots,\ti x_m)$ are the coordinates on $\ti U$, and $(\ti U',\ti\phi')$ is the corresponding a-chart on $\pd X$, then
${}^b\nabla_{\pd X}(x_1\frac{\pd}{\pd x_1})=0$ and ${}^b\nabla_{\pd X}(\ti x_1\frac{\pd}{\pd\ti x_1})=0$ induce the same connection on ${}^bN_{\pd X}$ over $\phi'(U')\cap\ti\phi'(\ti U')$. Over $\phi(U)\cap\ti\phi(\ti U)$ we have a change of coordinates $(x_1,\ldots,x_m)\rightsquigarrow(\ti x_1,\ldots,\ti x_m)$. Using Definition \ref{ac3def3}(b)(i), locally near $\phi'(U')\cap\ti\phi'(\ti U')$ we have 
\begin{equation*}
\ti x_1=F_1(x_1,\ldots,x_m)\cdot x_1^{a_{1,1}},
\end{equation*}
for $a_{1,1}\in(0,\iy)$ and $F_1(x_1,\ldots,x_m):U\ra(0,\iy)$ a-smooth. Here Definition \ref{ac3def3}(b)(i) also includes factors $x_2^{a_{2,1}}\cdots x_k^{a_{k,1}}$, but we must have $a_{2,1}=\cdots=a_{k,1}=0$ in this case.

It then follows as in Remark \ref{ac4rem3}(a) that $\ti x_1\frac{\pd}{\pd\ti x_1}=a_{1,1}^{-1}\cdot x_1\frac{\pd}{\pd x_1}$ locally near $\phi'(U')\cap\ti\phi'(\ti U')$ in $\pd X$. Thus ${}^b\nabla_{\pd X}(x_1\frac{\pd}{\pd x_1})=0$ and ${}^b\nabla_{\pd X}(\ti x_1\frac{\pd}{\pd\ti x_1})=0$ induce the same connection on ${}^bN_{\pd X}$, and ${}^b\nabla_{\pd X}$ is well defined. It is clearly locally trivial, as it is trivial over each $\phi'(U')$, and so is flat.

We call the manifold with a-corners {\it untwisted\/} if $({}^bN_{\pd X},{}^b\nabla_{\pd X})$ is globally isomorphic to the trivial flat line bundle $(\ul{\R},{}^b\nabla_0)$ on $\pd X$, and {\it twisted\/} otherwise.

\label{ac4def10}
\end{dfn}

\begin{ex}{\bf(a)} Let $W$ be a manifold with corners, and $\ti W=F_\Manc^\Manac(W)$, for $F_\Manc^\Manac$ as in Definition \ref{ac3def5}. Now the normal line bundle ${}^bN_{\pd W}$ in Definition \ref{ac2def9} is canonically trivial. But $F_\Manc^\Manac$ takes ${}^bN_{\pd W}$ to ${}^bN_{\pd\ti W}$, so ${}^bN_{\pd\ti W}$ is canonically trivial, and this trivialization is compatible with the flat connection on ${}^bN_{\pd\ti W}$ from Definition \ref{ac4def10}. Hence the manifold with a-corners $\ti W$ is untwisted.

One can in fact show that a manifold with a-corners $W'$ is untwisted if and only if $W'\cong F_\Manc^\Manac(W)$ for some manifold with corners~$W$.
\smallskip

\noindent{\bf(b)} Let $X=(\lb 0,\iy)\t\R)/\Z$ depending on $\al>0$ be as in Example \ref{ac3ex3}, where $(x,y)$ are the coordinates on $\lb 0,\iy)\t\R$ and $\Z$ acts by $n:(x,y)\mapsto (x^{\al^n},y+n)$. Then $\pd X\cong\R/\Z$ is a circle, where $y$ is the coordinate on $\R$ and $\Z$ acts by $n:y\mapsto y+n$. We may write ${}^bN_{\pd X}\ra\pd X$ as $(\R\t\R)/\Z\ra\R/\Z$, where $(y,z)$ are the coordinates on $\R^2$, and represent $z\cdot x\frac{\pd}{\pd x}$ over $(0,y)$ in $\lb 0,\iy)\t\R$. 

Since $x\mapsto x^{\al^n}$ maps $x\frac{\pd}{\pd x}\mapsto \al^{-n}\cdot x\frac{\pd}{\pd x}$, we see that $\Z$ acts on $\R^2$ as $n:(y,z)\mapsto (y+n,\al^n\cdot z)$. The connection ${}^b\nabla_{\pd X}$ on ${}^bN_{\pd X}$ has $z=\,\,$constant as its local constant sections. Therefore the holonomy of ${}^b\nabla_{\pd X}$ around $\R/\Z\cong\cS^1$ is multiplication by $\al$. So $X$ is untwisted if $\al=1$, and twisted otherwise.
\smallskip

Combining {\bf(a)},{\bf(b)} we see that if $X$ is as in Example \ref{ac3ex3} and $\al\ne 1$ then $X\not\cong F_\Manc^\Manac(W)$ for any $W\in\Manc$. Therefore $F_\Manc^\Manac:\Manc\ra\Manac$ {\it is not essentially surjective}. This also shows that {\it manifolds with a-corners up to a-diffeomorphism can depend nontrivially on continuous parameters}.
\label{ac4ex3}
\end{ex}

Next we define {\it weights\/} $\la$ on a manifold with a-corners $X$, and a line bundle $L_\la$ with flat b-connection ${}^b\nabla_\la$ for each weight $\la$.

\begin{dfn} Let $X$ be a manifold with a-corners. A {\it weight\/} on $X$ is a constant section $\la$ of $({}^bN_{\pd X},{}^b\nabla_{\pd X})$ in Definition \ref{ac4def10}, that is, $\la\in\Ga^\iy({}^bN_{\pd X})$ with ${}^b\nabla_{\pd X}\la=0$ in $\Ga^\iy({}^bN_{\pd X}\ot {}^bT^*(\pd X))$. Write $W(X)$ for the set of weights $\la$ for $X$. Then $W(X)$ is a real vector space.

Write $\pd X=\coprod_{i\in I}\pd_iX$, where $\pd_iX$ for $i\in I$ are the connected components of $\pd X$. If $X$ is untwisted then $({}^bN_{\pd X},{}^b\nabla_{\pd X})$ is noncanonically isomorphic to the trivial flat line bundle, so $W(X)$ is noncanonically isomorphic to the set of locally constant functions $\pd X\ra\R$, which is $\R^I$. 

If $X$ is twisted then $({}^bN_{\pd X},{}^b\nabla_{\pd X})$ is nontrivial on $\pd_iX$ for at least one $i\in I$, and then any weight $\la$ is zero on $\pd_iX$. In general $W(X)\cong\R^J$, where $J\subseteq I$ is the set of $i\in I$ such that $({}^bN_{\pd X},{}^b\nabla_{\pd X})$ is trivial on $\pd_iX$.

As in Definition \ref{ac4def7}, the orientation on ${}^bN_{\pd X}$ gives a notion of when a weight $\la$ has $\la>0$, $\la<0$, $\la\ge 0$, or $\la\le 0$. More generally, if $\la,\mu\in W(X)$ we write $\la>\mu$ if $\la-\mu>0$, and so on.

For each $\la\in W(X)$ we will define a real line bundle $L_\la$ on $X$ with a flat b-connection ${}^b\nabla_\la$, where we have a canonical identification
\e
(L_\la,{}^b\nabla_\la)\vert_{X^\ci}\cong (\ul{\R},{}^b\nabla_0)\vert_{X^\ci}
\label{ac4eq19}
\e
with the trivial flat line bundle $\ul{\R}$ on $X^\ci$. This $(L_\la,{}^b\nabla_\la)$ is not locally trivial on $X$ if $\la\ne 0$.

Fix a weight $\la\in W(X)$. Let $(U,\phi)$ be an a-chart on $X$, with $U\subseteq\R^{k,m}$ open, and $(x_1,\ldots,x_m)$ are the coordinates on $U$. Then for each $i=1,\ldots,k$ we have an a-chart $(U_i',\phi_i')$ on $\pd X$, where 
\begin{gather*}
U'_i=\bigl\{(y_1,\ldots,y_{m-1})\in\R^{k-1,m-1}:(y_1,\ldots,y_{i-1},0,y_i,\ldots,y_{m-1})\in U\subseteq\R^m_k\bigr\},\\
\phi_i'(y_1,\ldots,y_{m-1})=\bigl(\phi(y_1,\ldots,y_{i-1},0,y_i,\ldots,y_{m-1}),\phi_*(\{x_i=0\})\bigr).
\end{gather*}
Suppose for simplicity that $U_i'$ is connected for $i=1,\ldots,k$. Then there exist unique $\la_1,\ldots,\la_k\in\R$ such that $\la=\la_i\cdot x_i\frac{\pd}{\pd x_i}$ over $\phi_i'(U_i')\subseteq\pd X$ for $i=1,\ldots,k$, where we take $\la_i=0$ if $U_i'=\es$. Over $\phi(U)\subseteq X$ we define $L_\la$ to have basis the nonvanishing section $e_{U,\phi}$, with b-derivative
\e
{}^b\nabla_\la e_{U,\phi}=
\ts\sum_{i=1}^k \la_i\,e_{U,\phi}\ot(x_i^{-1}\d x_i).
\label{ac4eq20}
\e

Over $\phi(U^\ci)\subseteq X^\ci$ the isomorphism $(L_\la,{}^b\nabla_\la)\vert_{X^\ci}\cong (\ul{\R},{}^b\nabla_0)\vert_{X^\ci}$ identifies 
\e
e_{U,\phi}\vert_{\phi(U^\ci)}\cong x_1^{\la_1}\cdots x_k^{\la_k}\in C^\iy(\phi(U^\ci))=\Ga^\iy(\ul{\R}\vert_{\phi(U^\ci)}).
\label{ac4eq21}
\e
Note that \eq{ac4eq20}--\eq{ac4eq21} imply that $x_1^{-\la_1}\cdots x_k^{-\la_k}e_{U,\phi}$ in $\Ga^\iy(L_\la\vert_{\phi(U^\ci)})$ is identified with $1\in\Ga^\iy(\ul{\R}\vert_{\phi(U^\ci)})$, where ${}^b\nabla_\la(x_1^{-\la_1}\cdots x_k^{-\la_k}e_{U,\phi})=0$, so this isomorphism $L_\la\vert_{\phi(U^\ci)}\cong\ul{\R}\vert_{\phi(U^\ci)}$ does identify ${}^b\nabla_\la$ with ${}^b\nabla_0$.

Now suppose $(\ti U,\ti\phi)$ is another a-chart on $X$, with $\ti U\subseteq\R^{l,m}$ open, and define $(\ti x_1,\ldots,\ti x_m)$ and $(\ti U_i',\ti\phi_i'),\ti\la_i$ for $i=1,\ldots,l$ as above. Then we have a coordinate change $(x_1,\ldots,x_m)\rightsquigarrow(\ti x_1,\ldots,\ti x_m)$ from $(U,\phi)$ to $(\ti U,\ti\phi)$. By Definition \ref{ac3def3}(b)(i) for $j=1,\ldots,l$ locally over $\phi(U)\cap\ti\phi(\ti U)$ we may write
\e
\ti x_j(x_1,\ldots,x_m)=F_j(x_1,\ldots,x_m)\cdot x_1^{a_{1,j}}\cdots x_k^{a_{k,j}}
\label{ac4eq22}
\e
where $a_{i,j}\in[0,\iy)$ and $F_j:U\ra(0,\iy)$ is a-smooth. 

Over $\phi(U)\cap\ti\phi(\ti U)$, we define the sections $e_{U,\phi},e_{\ti U,\ti\phi}$ of $L_\la$ to be related by
\e
e_{\ti U,\ti\phi}=\ti x_1^{\ti\la_1}\cdots\ti x_l^{\ti\la_l}\cdot x_1^{-\la_1}\cdots x_k^{-\la_k}\cdot e_{U,\phi}.
\label{ac4eq23}
\e
This is compatible with the isomorphisms $L_\la\vert_{\phi(U^\ci)}\cong\ul{\R}\vert_{\phi(U^\ci)}$ and $L_\la\vert_{\ti\phi(\ti U^\ci)}\cong\ul{\R}\vert_{\ti\phi(\ti U^\ci)}$ above. Now \eq{ac4eq23} may not be defined at points in $\phi(U)\cap\ti\phi(\ti U)$ where some $x_i=0$ or $\ti x_j=0$. However, by \eq{ac4eq22} we may rewrite \eq{ac4eq23} as
\e
e_{\ti U,\ti\phi}=F_1^{\ti\la_1}(x_1,\ldots,x_m)\cdots F_l^{\ti\la_l}(x_1,\ldots,x_m)\cdot \ts\prod_{i=1}^k x_i^{\sum_{j=1}^la_{i,j}\ti\la_j-\la_i}\cdot e_{U,\phi}.
\label{ac4eq24}
\e
If $x_i=0$ at some point $z$ in $\phi(U)\cap\ti\phi(\ti U)$, then by considering the expressions for $\la$ in the coordinates $(x_1,\ldots,x_m)$ and $(\ti x_1,\ldots,\ti x_m)$ at $(z,\{x_i=0\})\in\pd X$, we see that $\la_i=\sum_{j=1}^la_{i,j}\ti\la_j$, so that the power of $x_i$ in \eq{ac4eq24} is zero. Hence the coefficient of $e_{U,\phi}$ in \eq{ac4eq24} is well defined and nonzero everywhere in~$\phi(U)\cap\ti\phi(\ti U)$.

It is easy to check that on triple overlaps $\phi(U)\cap\ti\phi(\ti U)\cap\hat\phi(\hat U)$ between a-charts $(U,\phi),(\ti U,\ti\phi),(\hat U,\hat\phi)$, the relations between sections $e_{U,\phi},e_{\smash{\ti U,\ti\phi}},e_{\smash{\hat U,\hat\phi}}$ are consistent. Thus we have defined a line bundle $L_\la$ and flat connection ${}^b\nabla_\la$ on $X$ for each weight~$\la\in W(X)$.
\label{ac4def11}
\end{dfn}

\begin{rem} For $X,\la,L_\la$ as above, smooth sections $s$ of $L_\la$ may be regarded as smooth functions $c$ on $X^\ci$ with prescribed growth rates near $\pd X$. In a chart $(U,\phi)$ as above, equation \eq{ac4eq21} implies that over $\phi(U)$ a section $s\in\Ga^\iy(L_\la)$ corresponds to a function $c$ on $X^\ci$ with growth of order~$O(x_1^{\la_1}\cdots x_k^{\la_k})$.

The line bundles $L_\la$ are intended to be used in problems involving prescribed growth rates at infinity. Rather than considering sections $e$ of a vector bundle $E$ over $X^\ci$ with growth $e=O(x_1^{\la_1}\cdots x_k^{\la_k})$ at infinity, one should consider sections $\bar e$ of the vector bundle $E\ot L_\la$ over $X$. It is often useful to take tensor products $E\ot L_\la$ of vector bundles $E\ra X$ with line bundles~$L_\la$.

If $X$ is twisted then there are fewer weights $\la$ on $X$, possibly only $\la=0$. So it may be helpful to suppose $X$ is untwisted for problems of this kind.
\label{ac4rem5}
\end{rem}

We can pull back weights $\la$ and line bundles $L_\la$ under interior maps:

\begin{dfn} Let $f:X\ra Y$ be an interior map of manifolds with a-corners, and $\la\in W(Y)$, with line bundle $L_\la\ra Y$ and canonical identification \eq{ac4eq19} from Definition \ref{ac4def11}. We will show that there is a unique weight $f^*(\la)$ in $W(X)$, and a unique isomorphism
\e
f^*(L_\la)\cong L_{f^*(\la)}
\label{ac4eq25}
\e
of line bundles on $X$, such that on restriction to $X^\ci$, the composition
\e
\ul{\R}\cong f\vert_{X^\ci}^*(\ul{\R})\cong f\vert_{X^\ci}^*(L_\la\vert_{Y^\ci})\cong
f^*(L_\la)\vert_{X^\ci}\cong L_{f^*(\la)}\vert_{X^\ci}\cong \ul{\R}
\label{ac4eq26}
\e
is the identity $\ul{\R}\cong\ul{\R}$. Here the second isomorphism of \eq{ac4eq26} comes from \eq{ac4eq19} for $\la$, the third holds as $f(Y^\ci)\subseteq X^\ci$ since $f$ is interior, the fourth is from \eq{ac4eq25}, and the fifth from \eq{ac4eq19} for~$f^*(\la)$.

It is easy to see from Definition \ref{ac4def11} that if $\mu,\mu'\in W(X)$ with an isomorphism $L_\mu\cong L_{\mu'}$ compatible with the canonical isomorphisms $L_\mu\vert_{X^\ci}\cong\ul{\R}\cong L_{\mu'}\vert_{X^\ci}$ from \eq{ac4eq19} then $\mu=\mu'$. Hence $f^*(\la)$ is unique if it exists. Also \eq{ac4eq26} determines the isomorphism \eq{ac4eq25} on $X^\ci$, and hence on $X=\ov{X^\ci}$ by continuity. Thus the isomorphism \eq{ac4eq25} is unique if it exists. Therefore it is enough to construct $\la$ and the isomorphism \eq{ac4eq25} locally on $\pd X$ and~$X$. 

Let $x\in S^k(X)\subseteq X$ with $f(x)=y\in S^l(Y)\subseteq Y$, and choose local coordinates $(x_1,\ldots,x_m)\in\R^{k,m}$ on open $x\in U\subseteq X$ near $x$ with $x=(0,\ldots,0)$ and $(y_1,\ldots,y_n)\in\R^{l,n}$ on open $y\in f(U)\subseteq V\subseteq Y$ with $y=(0,\ldots,0)$. Then by Definition \ref{ac3def3} we may write $f=(f_1,\ldots,f_n)$ with $f_j=f_j(x_1,\ldots,x_m)$ in coordinates near $x$, where for $j=1,\ldots,l$ we have
\begin{equation*}
f_j(x_1,\ldots,x_m)=F_j(x_1,\ldots,x_m)\cdot x_1^{a_{1,j}}\cdots x_k^{a_{k,j}}
\end{equation*}
near $x$, for $F_j:U\ra(0,\iy)$ a-smooth and $a_{i,j}\in[0,\iy)$. 

Write $\la\vert_{y_j=0}=\la_j\cdot y_j\frac{\pd}{\pd y_j}$ on the boundary face $y_j=0$ near $y$ for $j=1,\ldots,l$ and $\la_j\in\R$. Then $f^*(\la)$ is given by $f^*(\la)\vert_{x_i=0}=\bigl(\sum_{j=1}^la_{i,j}\la_j\bigr)\cdot x_i\frac{\pd}{\pd x_i}$ on the boundary face $x_i=0$ near $x$ for $i=1,\ldots,k$. One can check using Definition \ref{ac4def11} that with this local definition of $f^*(\la)$ there is a unique isomorphism \eq{ac4eq25} near $x$ satisfying the desired conditions. Thus $f^*(\la)$ and \eq{ac4eq25} are well defined. Pullbacks $f^*(\la)$ and the isomorphisms \eq{ac4eq25} are contravariantly functorial.
\label{ac4def12}
\end{dfn}

\begin{rem} To extend the material of this section to manifolds with corners and a-corners $X$ in \S\ref{ac35}, with $\pd X=\pd^{\rm c}X\amalg\pd^{\rm ac}X$ as in Remark \ref{ac4rem1}, we define $({}^bN_{\pd^{\rm ac}X},{}^b\nabla_{\pd^{\rm ac}X})$ over the a-boundary $\pd^{\rm ac}X$ only, and let $W(X)$ be the vector space of constant sections of ${}^bN_{\pd^{\rm ac}X}\ra\pd^{\rm ac}X$, and flat line bundles $(L_\la,{}^b\nabla_\la)$ for $\la\in W(X)$, where $(L_\la,{}^b\nabla_\la)$ is trivial near the ordinary boundary~$\pd^{\rm c}X$.
\label{ac4rem6}
\end{rem}

\section{Analysis on manifolds with a-corners}
\label{ac5}

We now discuss aspects of the analysis of partial differential operators on a manifold with a-corners $X$, regarded as acting on suitable Sobolev spaces of functions or sections of vector bundles on $X$. We focus in particular on elliptic operators and Fredholm properties, although the author expects that our theory may also have interesting applications to other classes of p.d.e.s, for example parabolic equations such as Ricci Flow or Mean Curvature Flow.

Most of this section is not very original: it is a translation into our language of ideas by other authors, especially those of Richard Melrose \cite{Melr1,Melr2,Melr3} and his school, but also those of Lockhart and McOwen \cite{Lock1,Lock2,LoMc} and others. This translation introduces some new changes in mathematical content (for example, a-smooth functions versus polyhomogeneous conormal functions, as we discuss in \S\ref{ac55}), and also a new point of view on well studied problems.

We hope to persuade the reader here and in \S\ref{ac6} that our `a-corners' language provides a useful and economical new way of thinking about, and writing about, many important topics in geometry and analysis.

In \S\ref{ac55} we will explain the relation of our theory to Melrose's theory of analysis on manifolds with corners~\cite{Melr1,Melr2,Melr3,Melr4}.

\subsection{Riemannian b-metrics}
\label{ac51}

This section is based on ideas of Melrose \cite{Melr1,Melr2,Melr3}, discussed in~\S\ref{ac55}.

\begin{dfn} Let $X$ be a manifold with a-corners, usually supposed compact. A {\it b-metric\/} $g$ on $X$ is $g\in \Ga^\iy(S^2({}^bT^*X))$ such that $g\vert_x$ is a positive definite quadratic form on ${}^bT_xX$ for all $x\in X$.

More generally, let $\la\in W(X)$ be a weight on $X$, and $L_\la$ the corresponding flat line bundle, as in \S\ref{ac44}. A {\it weighted b-metric\/ $g$ on\/ $X$ with weight\/} $\la$ is $g\in \Ga^\iy(S^2({}^bT^*X)\ot L_{-2\la})$ such that $g\vert_x$ is a positive definite quadratic form on ${}^bT_xX\ot L_\la\vert_x$ for all $x\in X$. An ordinary b-metric is weighted with weight~0.

\label{ac5def1}
\end{dfn}

Let $X$ be a compact manifold with a-corners with $\pd X\ne\es$, and $g$ a b-metric or weighted b-metric on $X$. The interior $X^\ci\subset X$ is an ordinary manifold, which is noncompact as $\pd X\ne\es$. Since ${}^bTX\vert_{X^\ci}=TX^\ci$ and $L_\la\vert_{X^\ci}=\ul{\R}$, we see that $g^\ci=g\vert_{X^\ci}$ is a Riemannian metric on $X^\ci$, as usual in differential geometry.

Now the fact that $g^\ci$ on $X^\ci$ extends to a (weighted) b-metric $g$ on $X$ implies that $g^\ci$ satisfies asymptotic conditions on the noncompact ends of $X^\ci$, and this is in fact a useful way of specifying asymptotic conditions on $(X^\ci,g^\ci)$ in many interesting problems. 

Here is the analogue of the fundamental theorem of Riemannian geometry. It may be proved by considering local coordinates $(x_1,\ldots,x_m)\in U\subseteq\R^{k,m}=\lb 0,\iy)^k\t\R^{m-k}$ on $X$, so that $(x_1,\ldots,x_m)\in U^\ci\subseteq(0,\iy)^k\t\R^{m-k}$ are local coordinates on $X^\ci$, considering the usual expression for the Levi-Civita connection $\nabla$ of $g\vert_{X^\ci}$ on $U^\ci$ in these coordinates, rewriting this using $x_i\frac{\pd}{\pd x_i},x_i^{-1}\d x_i$ rather than $\frac{\pd}{\pd x_i},\d x_i$ for $i=1,\ldots,k$, and then noting that this rewritten expression extends uniquely to an a-smooth b-connection on~$U\supseteq U^\ci$.

\begin{prop} Let\/ $g$ be a weighted b-metric on a manifold with a-corners $X,$ with weight\/ $\la$. Then there is a unique b-connection ${}^b\nabla$ on the vector bundle ${}^bTX\ot L_\la$ over $X$ called the \begin{bfseries}Levi-Civita b-connection\end{bfseries}, such that ${}^b\nabla\vert_{X^\ci}$ is identified with the usual Levi-Civita connection $\nabla$ of\/ $g\vert_{X^\ci}$ under the natural isomorphism $({}^bTX\ot L_\la)\vert_{X^\ci}\cong TX^\ci$. Also ${}^b\nabla$ is torsion-free with ${}^b\nabla g=0$.
\label{ac5prop1}
\end{prop}

Here is a simple example:

\begin{ex} Let $X$ be a compact manifold with a-boundary of dimension $m>0$, and suppose $X$ is untwisted, and $\pd X$ is nonempty and connected. Then there exist a compact subset $K\subset X^\ci$ and an a-diffeomorphism
\e
X\sm K\cong \lb 0,\ep)\t\pd X
\label{ac5eq1}
\e
for $\ep>0$, where $\pd X$ is a compact $(m-1)$-manifold. Write $t:X\sm K\ra \lb 0,\ep)$ for projection to $\lb 0,\ep)$ in \eq{ac5eq1}, so that $t$ is a-smooth. If $(x_1,\ldots,x_{m-1})\in\R^{m-1}$ are local coordinates on $\pd X$, then $(t,x_1,\ldots,x_{m-1})\in\lb 0,\iy)\t\R^{m-1}=\R^{1,m}$ are local coordinates on $X\sm K\subset X$ under~\eq{ac5eq1}.

The flat line bundle ${}^bN_{\pd X}$ of Definition \ref{ac4def10} is trivial (as $X$ is untwisted), with constant basis section $t\frac{\pd}{\pd t}$. Thus, as $\pd X$ is connected, weights $\la\in W(X)$ in Definition \ref{ac4def11} are $\la=\la'\cdot t\frac{\pd}{\pd t}$ for $\la'\in\R$, so that $W(X)\cong\R$. Note that this trivialization $W(X)\cong\R$ depends on the choice of a-diffeomorphism \eq{ac5eq1}. If we replace $t:X\sm K\ra \lb 0,\iy)$ by $\ti t=t^\al$ for $\al>0$ then $\ti t\frac{\pd}{\pd\ti t}=\al^{-1}\cdot t\frac{\pd}{\pd t}$, so $\la=\la'\cdot t\frac{\pd}{\pd t}=\ti \la'\cdot\ti t\frac{\pd}{\pd\ti t}$ for~$\ti\la'=\al\la'$.

We can identify a-smooth sections $e$ of the line bundle $L_\la$ of Definition \ref{ac4def11} with smooth functions $f:X^\ci\ra\R$ such that the function $t^{-\la'}\cdot f\vert_{X^\ci\sm K}:X^\ci\sm K\ra\R$ extends to an a-smooth function~$X\sm K\ra\R$.

Now suppose that $g$ is a weighted b-metric on $X$ with weight $\la$. Then under the identification \eq{ac5eq1} we may write
\e
g\vert_{X\sm K}\cong t^{-2\la'}\bigl(h_t+2h'_t t^{-1}\d t+h''_t t^{-2}(\d t)^2\bigr),
\label{ac5eq2}
\e
where for $t\in\lb 0,\ep)$ we have $h_t\in \Ga^\iy(S^2T^*(\pd X))$, $h'_t\in \Ga^\iy(T^*(\pd X))$ and $h_t''\in C^\iy(\pd X)$, which are a-smooth in $t\in\lb 0,\ep)$, and such that $h_t+2h'_t t^{-1}\d t+h''_t t^{-2}(\d t)^2$ is positive definite on $T(\pd X)\op\langle t\frac{\pd}{\pd t}\rangle_\R$, so that in particular $h_t$ is a Riemannian metric on $\pd X$ for $t\in\lb 0,\ep)$, and $h''_t>0$.

By definition of a-smooth functions in \S\ref{ac31}, and using the compactness of $\pd X$, we see that there exists small $\al>0$ such that
\e
h_t=h_0+O(t^\al),\quad h'_t=h'_0+O(t^\al),\quad h''_t=h''_0+O(t^\al)\quad\text{as $t\ra 0$.}
\label{ac5eq3}
\e
Combining \eq{ac5eq2}--\eq{ac5eq3} gives
\e
g\vert_{X\sm K}\cong t^{-2\la'}\bigl(h_0+2h'_0 t^{-1}\d t+h''_0 t^{-2}(\d t)^2+O(t^\al)\bigr)\quad\text{as $t\ra 0$.}
\label{ac5eq4}
\e

Let us now assume that $h_0'=0$ and $h_0''=1$. Divide into three cases (a) $\la'=0$; (b) $\la'>0$; and (c) $\la'<0$, and change variables for $t$ as follows:
\begin{itemize}
\setlength{\itemsep}{0pt}
\setlength{\parsep}{0pt}
\item[(a)] When $\la'=0$, so $g$ is unweighted, change variables from $t$ in $(0,\ep)$ to $r=-\log t$ in $(R,\iy)$, where $R=-\log\ep$. Then \eq{ac5eq1} induces a diffeomorphism
\e
X^\ci\sm K\cong (R,\iy)\t\pd X,
\label{ac5eq5}
\e
and under \eq{ac5eq5}, equation \eq{ac5eq4} becomes
\e
g^\ci\vert_{X^\ci\sm K}\cong h_0+(\d r)^2+O(e^{-\al r})\quad\text{as $r\ra \iy$.}
\label{ac5eq6}
\e
That is, the noncompact Riemannian manifold $(X^\ci,g^\ci)$ is {\it Asymptotically Cylindrical\/} ({\it ACy\/}), with exponential decay as $r\ra\iy$ to the Riemannian cylinder $\R\t\pd X$ with metric $h_0+(\d r)^2$.
\item[(b)] When $\la'>0$, change variables from $t$ in $(0,\ep)$ to $r=(\la')^{-1}t^{-\la'}$ in $(R,\iy)$, where $R=(\la')^{-1}\ep^{-\la'}$. Then \eq{ac5eq1} induces a diffeomorphism
\e
X^\ci\sm K\cong (R,\iy)\t\pd X,
\label{ac5eq7}
\e
and under \eq{ac5eq7} we have
\e
g^\ci\vert_{X^\ci\sm K}\cong r^2[(\la')^2h_0]+(\d r)^2+O(r^{-\al/\la'})\quad\text{as $r\ra \iy$.}
\label{ac5eq8}
\e
That is, $(X^\ci,g^\ci)$ is {\it Asymptotically Conical\/} ({\it ACo\/}), with polynomial decay to the Riemannian cone $(0,\iy)\t\pd X$ with metric $r^2[(\la')^2h_0]+(\d r)^2$ as $r\ra\iy$, where the error $O(r^{-\al/\la'})$ in \eq{ac5eq8} is measured using $g^\ci$.
\item[(c)] When $\la'<0$, change variables from $t$ in $(0,\ep)$ to $r=-(\la')^{-1}t^{-\la'}$ in $(0,R)$, where $R=-(\la')^{-1}\ep^{-\la'}$. Then \eq{ac5eq1} induces a diffeomorphism
\e
X^\ci\sm K\cong (0,R)\t\pd X,
\label{ac5eq9}
\e
and under \eq{ac5eq9} we have
\begin{equation*}
g^\ci\vert_{X^\ci\sm K}\cong r^2[(\la')^2h_0]+(\d r)^2+O(r^{-\al/\la'})\quad\text{as $r\ra 0$.}
\end{equation*}
That is, $(X^\ci,g^\ci)$ has a {\it conical singularity}, with polynomial decay as $r\ra 0$ to the Riemannian cone $(0,\iy)\t\pd X$ with metric $r^2[(\la')^2h_0]+(\d r)^2$.

It makes geometric sense to compactify $(X^\ci,g^\ci)$ by adding a single point at $r=0$, the vertex of the cone.
\end{itemize}
Thus, weighted b-metrics include Asymptotically Cylindrical metrics, Asymptotically Conical metrics, and metrics with conical singularities.
\label{ac5ex1}
\end{ex}

\subsection{Weighted Sobolev spaces on manifolds with a-corners}
\label{ac52}

Let $X$ be a compact manifold with a-corners with $\pd X\ne\es$, so that the interior $X^\ci$ is a noncompact manifold, and $E\ra X$ be a vector bundle. We now define Banach spaces $L^p(X),L^p_k(X),L^p(X)_\la,L^p_k(X)_\la$ of functions on $X^\ci$, and Banach spaces $L^p(E),L^p_k(E),L^p(E)_\la,L^p_k(E)_\la$ of sections of $E\vert_{X^\ci}$. Here although the geometry happens on the noncompact manifold $X^\ci$, the compactification $X$ is used to control the asymptotic behaviour at infinity in~$X^\ci$.

As we will see in \S\ref{ac53}--\S\ref{ac54}, the definition is designed to ensure that linear elliptic operators on $X^\ci$ with suitable extension properties to $X$ induce Fredholm maps between these Banach spaces. There is a very similar story for weighted H\"older spaces $C^{k,\al}(E),C^{k,\al}(E)_\la$ on $X^\ci$, but for brevity we will not explain it.

In many important special cases, our weighted Sobolev spaces are equivalent to those considered by other authors, as we explain in Remark \ref{ac5rem2}(c), and the Fredholm results in \S\ref{ac53}--\S\ref{ac54} will be deduced from these authors' work.

\begin{dfn} Let $X$ be a compact manifold with a-corners, and $E\ra X$ a vector bundle. For each real number $p\in(1,\iy)$ and integer $k\ge 0$ we will define the {\it Lebesgue space\/} $L^p(E)$ and {\it Sobolev space\/} $L^p_k(E)$, which are Banach spaces of sections of the vector bundle $E\vert_{X^\ci}$ over the ordinary manifold $X^\ci$. Also, if $\la\in W(X)$ is a weight, we will define the {\it weighted Lebesgue space\/} $L^p(E)_\la$ and the {\it weighted Sobolev space\/} $L^p_k(E)_\la$ by
\e
L^p(E)_\la=L^p(E\ot L_\la)\quad\text{and}\quad L^p_k(E)_\la=L^p_k(E\ot L_\la),
\label{ac5eq10}
\e
where $\la\in W(X)$ and $L_\la$ are as in \S\ref{ac44}.

To define $L^p(E),L^p_k(E)$ we need to make some additional choices. Let $g$ be a b-metric on $X$ ({\it not\/} a weighted b-metric), as in \S\ref{ac51}. Let $h^E$ be an a-smooth metric on the fibres of $E$, that is, $h^E\in\Ga^\iy(S^2E^*)$ with $h^E\vert_x$ a positive definite quadratic form on $E_x$ for all $x\in X$, and let ${}^b\nabla^E$ be a b-connection on~$E$.

Combining ${}^b\nabla^E$ with the Levi-Civita b-connection ${}^b\nabla$ on ${}^bTX$ from \S\ref{ac51} induces b-connections ${}^b\nabla^E$ on $E\ot\bigot^j{}^bT^*X$ for $j=1,2,\ldots.$ Thus, for any $e\in\Ga^\iy(E)$ we may form the $k^{\rm th}$ derivative $({}^b\nabla^E)^ke\in\Ga^\iy(E\ot\bigot^k{}^bT^*X)$ for all $k=0,1,2,\ldots.$ Also, combining $h^E$ with the inverse of $g$ induces metrics on the fibres of $E\ot\bigot^k{}^bT^*X$, so we can take norms~$\bmd{({}^b\nabla^E)^ke}\in C^\iy(X)$.

For $p>1$, define the {\it Lebesgue space\/} $L^p(E)$ to be the vector space of locally integrable sections $e$ of $E\vert_{X^\ci}$ over $X^\ci$ for which the norm
\begin{equation*}
\nm{e}_{L^p}=\left(\int_{X^\ci}\md{e}^p\d V_{g^\ci}\right)^{1/p}
\end{equation*}
is finite. Here $\d V_{g^\ci}$ is the volume form of the metric $g^\ci=g\vert_{X^\ci}$. For $p>1$ and $k\ge 0$ define the {\it Sobolev space\/} $L^p_k(E)$ to be the vector space of locally integrable 
sections $e$ of $E\vert_{X^\ci}$ such that $e$ is $k$ times weakly differentiable and the norm
\begin{equation*}
\nm{e}_{L^p_k}=\biggl(\sum_{j=0}^k\int_{X^\ci}\bmd{({}^b\nabla^E)^je}^p\d V_{g^\ci}\biggr)^{1/p}
\end{equation*}
is finite. Then $L^p(E),L^p_k(E)$ are Banach spaces w.r.t.\ the norms $\nm{\,.\,}_{L^p},\nm{\,.\,}_{L^p_k}$, and $L^2(E),L^2_k(E)$ are Hilbert spaces. For each $\la\in W(X)$, define the {\it weighted Lebesgue space\/} $L^p(E)_\la$ and the {\it weighted Sobolev space\/} $L^p_k(E)_\la$ by~\eq{ac5eq10}.

When $E$ is the trivial line bundle $\ul{\R}=\R\t X\ra X$ we will write $L^p(X),\ab L^p_k(X),\ab L^p(X)_\la,L^p_k(X)_\la$ in place of $L^p(E),L^p_k(E),L^p(E)_\la,L^p_k(E)_\la$, and regard $L^p(X),\ldots,L^p_k(X)_\la$ as Sobolev spaces of functions on~$X$.
\label{ac5def2}
\end{dfn}

The next lemma is easy to prove. It is important because it shows that the spaces $L^p(E),\ab L^p_k(E),\ab L^p(E)_\la,L^p_k(E)_\la$ are intrinsic to~$X,E$.

\begin{lem} In Definition\/ {\rm\ref{ac5def2},} noting that $X$ is compact, if we consider $L^p(E),\ab L^p_k(E),\ab L^p(E)_\la,L^p_k(E)_\la$ just as topological vector spaces of sections of\/ $E\vert_{X^\ci},$ then they depend only on $X,E,p,k,\la,$ and are independent of the auxiliary choices of\/ $g,h^E,\nabla^E$. The norms $\nm{\,.\,}_{L^p},\nm{\,.\,}_{L^p_k}$ do depend on $g,h^E,\nabla^E,$ but only up to bounded factors.
\label{ac5lem1}
\end{lem}

\begin{rem}{\bf(a)} The definitions above also give good Banach spaces $L^1(E),\ab\ldots,L^1_k(E)_\la$ when $p=1$, and Lemma \ref{ac5lem1} holds. But the Fredholm results in \S\ref{ac53} fail when $p=1$ even for ordinary manifolds, so we exclude this case.
\smallskip

\noindent{\bf(b)} If $X$ is not compact above then Lemma \ref{ac5lem1} fails, and $L^p(E),\ab L^p_k(E),\ab L^p(E)_\la,\ab L^p_k(E)_\la$ depend on $g,h^E,\nabla^E$, which makes them less interesting.

\smallskip

\noindent{\bf(c)} It is important in Definition \ref{ac5def2} that we define $L^p(E),\ldots,L^p_k(E)_\la$ using a b-metric $g$, {\it not\/} a weighted b-metric. If allow $g$ to be a weighted b-metric, the the Fredholm properties of elliptic operators in \S\ref{ac53} may no longer hold. 

For example, take $X$ to be an $n$-ball $B^n$, with $X^\ci\cong\R^n$, and give $\R^n$ the Euclidean metric $g_0$. Lockhart \cite{Lock1} explains that the Laplacian $\De$ on $\R^n$ gives a non-Fredholm map $\De:L^p_2(\R^n)\ra L^p(\R^n)$ between Sobolev spaces defined using $g_0$, and says that this is because $L^p_2(\R^n)$ is the `wrong' space, and one should instead consider certain `weighted Sobolev spaces' on~$\R^n$.

In our language, this is because $(\R^n,g_0)$ is Asymptotically Conical, as in Example \ref{ac5ex1}(b), and so comes from a weighted b-metric on $X$, and $L^p_k(\R^n)$ is defined using a weighted b-metric. But Lockhart's weighted Sobolev spaces on $\R^n$ are equivalent to those in Definition \ref{ac5def2} using an unweighted b-metric on~$X$.

\smallskip

\noindent{\bf(d)} Let $X$ be a compact, untwisted manifold with a-boundary, with $\pd X\ne\es$. Then our $L^p_k(X)_\la$ are equivalent to the weighted Sobolev spaces on $X^\ci$ defined by Lockhart and McOwen \cite{Lock1,Lock2,LoMc}. When $p=2$, they are equivalent to the weighted Sobolev spaces considered by Melrose and Mendoza \cite[\S 2]{MeMe}, \cite[p.~10]{Melr3}.

Now let $X$ be a compact, untwisted manifold with a-corners. When $p=2$, our weighted Sobolev spaces $L^2_k(X)_\la$ are equivalent to those discussed by Loya \cite{Loya1,Loya2}, and are also implicit in the work of Melrose, Nistor and Piazza~\cite{MeNi,MePi}.

\label{ac5rem2}
\end{rem}

\subsection{Elliptic operators and Fredholm properties}
\label{ac53}

We discuss partial differential operators on manifolds with a-corners. These are related to (but less general than) the `b-pseudodifferential operators' on a compact manifold with corners studied by Melrose, Nistor and Piazza~\cite{MeNi,MePi}.

\begin{dfn} Let $X$ be a compact manifold with a-corners, and $E,F\ra X$ be vector bundles. A linear map $P:\Ga^\iy(E)\ra\Ga^\iy(F)$ is a {\it partial differential operator\/} ({\it p.d.o.}) {\it of order\/} $l\ge 0$ if $P(e)\vert_x\in F_x$ depends linearly only on $e$ and its first $l$ b-derivatives at $x$, for each $x\in X$ and all~$e\in\Ga^\iy(E)$.

We may write this more explicitly as follows. Choose b-connections ${}^b\nabla^E$ on $E$ and ${}^b\nabla$ on ${}^bTX$. Then ${}^b\nabla^E,{}^b\nabla$ induce b-connections ${}^b\nabla^E$ on $E\ot\bigot^j{}^bT^*X$ for $j\ge 1$ so given $e\in\Ga^\iy(E)$ we can form $({}^b\nabla^E)^je\in\Ga^\iy(E\ot\bigot^j{}^bT^*X)$ for $j\ge 0.$ Then $P:\Ga^\iy(E)\ra\Ga^\iy(F)$ is a p.d.o.\ of order $l$ if it is of the form
\e
P(e)=\ts\sum_{j=0}^l a_j\cdot ({}^b\nabla^E)^je,
\label{ac5eq11}
\e
where $a_j\in \Ga^\iy\bigl(E^*\ot\bigot^j{}^bTX\ot F\bigr)$ for $j=0,\ldots,l$. The notion of p.d.o.\ of order $l$, and the leading term $a_l$, are both independent of the choices of ${}^b\nabla^E,{}^b\nabla$, although $a_0,a_1,\ldots,a_{l-1}$ do depend on~${}^b\nabla^E,{}^b\nabla$.

For any $p>1$ and $k\ge 0$, $P$ induces a morphism of Banach spaces
\e
P^p_k:L^p_{k+l}(E)\longra L^p_k(F),
\label{ac5eq12}
\e
given by \eq{ac5eq11} restricted to $X^\ci$, with $e\in L^p_{k+l}(E)$.

Let $\la$ be a weight on $X$, so that \S\ref{ac44} defines a line bundle $L_\la$ on $X$ with flat connection ${}^b\nabla_\la$. For $j=0,1,\ldots,$ write ${}^b\nabla^E$ for the b-connection on $E\ot L_\la\ot\bigot^j{}^bT^*X$ given by combining the b-connections ${}^b\nabla^E$ on $E\ot\bigot^j{}^bT^*X$ and ${}^b\nabla_\la$ on $L_\la$. Now define a p.d.o.\ $P_\la:\Ga^\iy(E\ot L_\la)\ra\Ga^\iy(F\ot L_\la)$ of order $l$ by
\begin{equation*}
P_\la(e)=\ts\sum_{j=0}^l a_j\cdot ({}^b\nabla_\la^E)^je,
\end{equation*}
where we now regard $a_j$ as lying in $\Ga^\iy\bigl((E\ot L_\la)^*\ot\bigot^j{}^bTX\ot (F\ot L_\la)\bigr)$, since $L_\la^*\ot L_\la$ is trivial. Then $P_\la$ depends only on $X,E,F,P,\la$, and not on the choices of ${}^b\nabla^E,{}^b\nabla$. This shows that any p.d.o.\ $P:\Ga^\iy(E)\ra\Ga^\iy(F)$ can be `twisted by $L_\la$' to give $P_\la:\Ga^\iy(E\ot L_\la)\ra\Ga^\iy(F\ot L_\la)$, for all~$\la\in W(X)$.

Applying \eq{ac5eq12} to $P_\la$ rather than $P$ gives a morphism of Banach spaces
\e
P^p_{k,\la}:L^p_{k+l}(E)_\la\longra L^p_k(F)_\la,
\label{ac5eq13}
\e
for all $p>1$, $k\ge 0$ and $\la\in W(X)$.

We may restrict a p.d.o.\ $P$ on $X$ to the boundary $\pd X$ and $k$-corners $C_k(X)$, as follows. For $P\vert_{\pd X}$, from \S\ref{ac41}--\S\ref{ac42} we have an a-smooth map $i_X:\pd X\ra X$ and a projection ${}^b\pi_T:i_X^*({}^bTX)\ra{}^bT(\pd X)$ as in \eq{ac4eq14}. Thus we may form
\begin{equation*}
\ts\bigot^j{}^b\pi_T[i_X^*(a_j)]\in \Ga^\iy\bigl(E^*\vert_{\pd X}\ot\bigot^j{}^bT(\pd X)\ot F\vert_{\pd X}\bigr).
\end{equation*}
Writing $E\vert_{\pd X},F\vert_{\pd X}$ for $i_X^*(E),i_X^*(F)$, define $P\vert_{\pd X}:\Ga^\iy(E\vert_{\pd X})\ra \Ga^\iy(F\vert_{\pd X})$~by
\e
P\vert_{\pd X}(e')=\ts\sum_{j=0}^l \ts\bigot^j{}^b\pi_T[i_X^*(a_j)]\cdot ({}^b\nabla_\pd^E)^je',
\label{ac5eq14}
\e
where ${}^b\nabla_\pd^E$ are the b-connections on $E\vert_{\pd X}\ot\bigot^j{}^bT^*(\pd X)\ra\pd X$ induced by ${}^b\nabla^E$ on $E\ra X$ and ${}^b\nabla$ on ${}^bTX\ra X$. Then $P\vert_{\pd X}$ is a p.d.o.\ of order $l$ on~$X$.

Similarly we may define $P\vert_{C_k(X)}:\Ga^\iy(E\vert_{C_k(X)})\ra \Ga^\iy(F\vert_{C_k(X)})$ and $P\vert_{\pd^kX}:\Ga^\iy(E\vert_{\pd^kX})\ra \Ga^\iy(F\vert_{\pd^kX})$ for all~$k=0,1,\ldots,\dim X$.
\label{ac5def3}
\end{dfn}

We can also consider {\it nonlinear\/} partial differential operators, in which $P(e)\vert_x$ is a nonlinear a-smooth function of $x$ and $({}^b\nabla^E)^je\vert_x$ for $j=0,\ldots,l$. These are very important, but for brevity we will not discuss them here.

\begin{rem} The canonical restrictions $P\vert_{\pd X}$ above work because our partial differential operators $P$ are defined using b-(co)tangent bundles ${}^bTX,{}^bT^*X$ rather than (co)tangent bundles $TX,T^*X$. To see this, compare equations \eq{ac2eq10}, \eq{ac2eq11} and \eq{ac4eq14}: the natural morphisms ${}^b\pi_T:i_X^*({}^bTX)\ra{}^bT(\pd X)$ and $\d i_X:T(\pd X)\ra i_X^*(TX)$ go in opposite directions, and $\d i_X$ goes the wrong way to be able to define $P\vert_{\pd X}$ in~\eq{ac5eq14}.
\label{ac5rem3}
\end{rem}

Next we discuss linear elliptic operators on manifolds with a-corners.

\begin{dfn} Let $X$ be a manifold with a-corners, and $P:\Ga^\iy(E)\ra\Ga^\iy(F)$ a partial differential operator of order $l\ge 0$. For all $x\in X$ and $\xi\in {}^bT_x^*X$, define a linear map $\si_P(x,\xi):E_x\ra F_x$ by $\si_P(x,\xi):e_x\mapsto a_l\vert_x\cdot (e_x\ot \bigot^l\xi)$. We call $\si_P$ the {\it principal symbol\/} of $P$. We call $P$ {\it elliptic\/} if $\si_P(x,\xi):E_x\ra F_x$ is an isomorphism for all $x\in X$ and~$0\ne\xi\in {}^bT_x^*X$.

As in Definition \ref{ac5def3} we may form the restriction $P\vert_{\pd X}$. The principal symbols $\si_P,\si_{P\vert_{\pd X}}$ of $P$ and $P\vert_{\pd X}$ are related by
\begin{equation*}
\si_{P\vert_{\pd X}}(x',\xi')=\si_P\bigl(i_X(x'),{}^b\pi_T^*\vert_{x'}(\xi')\bigr),
\end{equation*}
where $x'\in\pd X$ and $\xi'\in {}^bT_{x'}^*(\pd X)$, and ${}^b\pi_T^*\vert_{x'}:{}^bT_{x'}^*(\pd X)\ra {}^bT_{i_X(x')}^*X$ is dual to ${}^b\pi_T$ in \eq{ac4eq14} at $x'$. Since ${}^b\pi_T^*\vert_{x'}$ is injective, we see that $P$ elliptic implies that $P\vert_{\pd X}$ is elliptic, and similarly $P\vert_{C_k(X)},P\vert_{\pd^kX}$ are elliptic for~$k\ge 0$.
\label{ac5def4}
\end{dfn}

Let $P:\Ga^\iy(E)\ra\Ga^\iy(F)$ be an elliptic p.d.o.\ of order $l$ on a compact manifold with a-corners $X$. We are interested in the question of when the Banach space morphisms \eq{ac5eq12}--\eq{ac5eq13} are Fredholm. When $\pd X=\es$, so that $X$ is a compact manifold, it is well known that \eq{ac5eq12} is Fredholm for any $p>1$ and $k\ge 0$, with index given by the Atiyah--Singer Index Theorem. Much work has been done on similar problems by many authors, including Lockhart and McOwen \cite{Lock1,Lock2,LoMc}, Loya \cite{Loya1,Loya2}, Melrose et al.\ \cite{Melr3,MeMe,MeNi,MePi}, and Piazza \cite{Piaz}.

The next theorem summarizes results of Lockhart--McOwen \cite{LoMc}, Melrose--Mendoza \cite{MeMe,Melr3} and Piazza \cite{Piaz} for compact manifolds with boundary, translated into our language.

\begin{thm} Let\/ $X$ be a compact, untwisted manifold with a-boundary. Write $\pd_1X,\ldots,\pd_NX$ for the connected components of\/ $\pd X$. By choosing a trivialization \eq{ac5eq1} near $\pd X$ we may identify $W(X)\cong\R^N,$ where $\la\in W(X)$ is identified with $(\la_1,\ldots,\la_N)\in\R^N$ with\/ $\la\cong \la_jt\frac{\pd}{\pd t}$ on $\pd_jX$ for\/~$j=1,\ldots,N$.

Let\/ $E,F\ra X$ be vector bundles with\/ $\rank E=\rank F,$ and\/ $P:\Ga^\iy(E)\ra\Ga^\iy(F)$ be an elliptic partial differential operator of order\/ $l>0$. Then:
\begin{itemize}
\setlength{\itemsep}{0pt}
\setlength{\parsep}{0pt}
\item[{\bf(a)}] There exist discrete subsets\/ $\cD_j\subset\R$ for\/ $j=1,\ldots,N$ such that 
\begin{equation*}
(P_\la\vert_{\pd_j X})^p_k:L^p_{k+l}(E\ot L_\la\vert_{\pd_j X})\longra L^p_k(F\ot L_\la\vert_{\pd_j X})
\end{equation*}
is an isomorphism if and only if\/ $\la_j\notin\cD_j,$ for all\/ $p>1,$ $k\ge 0$ and\/ $\la\in W(X),$ where $\la$ corresponds to $(\la_1,\ldots,\la_N)\in \R^N$ as above.
\item[{\bf(b)}] For all\/ $p>1,$ $k\ge 0$ and\/ $\la\in W(X),$ the operator $P^p_{k,\la}$ in \eq{ac5eq13} is Fredholm if and only if\/ $\la_j\notin\cD_j$ for $j=1,\ldots,N$. When it is Fredholm, the kernel\/ $\Ker P^p_{k,\la}$ and index\/ $\ind P^p_{k,\la}$ depend only on\/ $X,E,F,P$ and the connected components of\/ $\R\sm\cD_j$ containing\/ $\la_j$ for\/ $j=1,\ldots,N$. 

If\/ $\la,\la'\in W(X)$ with $\la\ge\la'$ in the sense of Definition\/ {\rm\ref{ac4def11},} then
\e
\Ker P^p_{k,\la}\subseteq \Ker P^p_{k,\la'}.
\label{ac5eq15}
\e
\item[{\bf(c)}] When $P^p_{k,\la}$ is Fredholm, there is an `Index Theorem' that writes
\e
\ind P^p_{k,\la}=\mathop{\text{\rm t-ind}}(P)+\ts\sum_{j=1}^N\mathop{\text{\rm b-ind}}(P\vert_{\pd_jX},\la_j),
\label{ac5eq16}
\e
where $\mathop{\text{\rm t-ind}}(P)$ is a `topological index' similar to that in the Atiyah--Singer Index Theorem, and\/ $\mathop{\text{\rm b-ind}}(P\vert_{\pd_jX},\la_j)$ is a `boundary index', roughly a count of eigenvalues in the interval\/ $[0,\la_j]$ of an operator on\/~$\pd_jX$.
\item[{\bf(d)}] By choosing a b-metric $g$ on $X$ and metrics $h^E,h^F$ on the fibres of\/ $E,F,$ we may define the \begin{bfseries}formal adjoint\end{bfseries}\/ $P^*:\Ga^\iy(F)\ra\Ga^\iy(E)$. It is an elliptic p.d.o.\ of order $l,$ which is characterized by the property that if\/ $\la\in W(X),$ $e\in L^2_l(E)_\la$ and\/ $f\in L^2_l(F)_{-\la}$ then
\begin{equation*}
\int_{X^\ci} h^F\bigl(P^2_{0,\la}(e),f\bigr)\,\d V_{g^\ci}=\int_{X^\ci} h^E\bigl(e,(P^*)^2_{0,-\la}(f)\bigr)\,\d V_{g^\ci},
\end{equation*}
where the integrals are well-defined.

Then for all\/ $p>1,k\ge 0$ and\/ $\la\in W(X),$ the operator $P^p_{k,\la}$ in \eq{ac5eq13} is Fredholm if and only if\/ $(P^*)^p_{k,-\la}$ is Fredholm, and in this case
\begin{gather*}
\Coker P^p_{k,\la}\cong \Ker(P^*)^p_{k,-\la},\qquad
\Coker (P^*)^p_{k,-\la}\cong \Ker P^p_{k,\la},\\
\text{and\/}\quad \ind P^p_{k,\la}= -\ind (P^*)^p_{k,-\la}=\dim\Ker P^p_{k,\la}-\dim \Ker(P^*)^p_{k,-\la}.
\end{gather*}
\end{itemize}
\label{ac5thm1}
\end{thm}

\begin{proof} Equation \eq{ac5eq15} is obvious as $L^p_{k+l}(E)_\la\subseteq L^p_{k+l}(E)_{\la'}$ if $\la\ge\la'$. Parts (a)--(c) but without computing $\mathop{\text{\rm t-ind}}(P)$ are proved in Lockhart and McOwen \cite[Th.~8.1]{LoMc}. Part (d) also follows from results in \cite[\S 3 \& \S 7]{LoMc}. When $P$ is self-adjoint they also prove \eq{ac5eq17} with $\mathop{\text{\rm t-ind}}(P)=0,$ \cite[Th.~7.4]{LoMc}. 

Melrose and Mendoza \cite[Th.~6.17]{MeMe} also independently prove (a)--(c) when $p=2$, but without computing $\mathop{\text{\rm t-ind}}(P)$. When $P$ is a twisted Dirac operator and $p=2$, parts (a)--(c) are proved at length by Melrose \cite{Melr3}, see in particular \cite[Th.s 5.40, 6.5, 9.1]{Melr3} and the formula \cite[(9.5)]{Melr3} for $\mathop{\text{\rm t-ind}}(P)$. Assuming \cite{MeMe,Melr3}, part (c) was proved for general elliptic $P$ by Piazza \cite{Piaz} when $p=2$, which by (b) also implies (c) for general~$p>1$.
\end{proof}

Theorem \ref{ac5thm1} was generalized to manifolds with corners by Loya and Melrose \cite{Loya1,Loya2,LoMe}, using material on b-pseudodifferential operators on manifolds with corners due to Melrose, Nistor and Piazza \cite{MeNi,MePi}. Parts (a),(b) of the next result are taken from \cite[Th.~B.2]{LoMe}, \cite[Th.~1.2]{Loya2}, translated into our language.

\begin{thm} Let\/ $X$ be a compact, untwisted manifold with a-corners. Write $\pd_1X,\ldots,\pd_NX$ for the connected components of\/ $\pd X$. Suppose that\/ $i_X\vert_{\pd_jX}:\pd_jX\ra X$ is injective for each\/ $j=1,\ldots,N$. By choosing local coordinates $t_j\in\lb 0,\iy)$ on $X$ near $i_X(\pd_jX)$ normal to $\pd_jX$ we identify $W(X)\cong\R^N,$ where $\la\in W(X)$ is identified with $(\la_1,\ldots,\la_N)\in\R^N$ with\/ $\la\cong \la_jt_j\frac{\pd}{\pd t_j}$ on $\pd_jX$ for\/ $j=1,\ldots,N$. Write\/ $e_j\in W(X)$ for the vector corresponding to $(\la_1,\ldots,\la_N)$ in $\R^N$ with\/ $\la_i=\de_{ij},$ for all\/~$i,j=1,\ldots,N$.

Let\/ $E,F\ra X$ be vector bundles with\/ $\rank E=\rank F,$ and\/ $P:\Ga^\iy(E)\ra\Ga^\iy(F)$ be an elliptic partial differential operator of order\/ $l>0$. Then:
\begin{itemize}
\setlength{\itemsep}{0pt}
\setlength{\parsep}{0pt}
\item[{\bf(a)}] For all\/ $k\ge 0$ and\/ $\la\in W(X),$ the operator on $X$
\e
P^2_{k,\la}:L^2_{k+l}(E)_\la\longra L^2_k(F)_\la
\label{ac5eq17}
\e
is Fredholm if and only if the operator on $\pd_jX$
\e
(P_{\la+\mu e_j}\vert_{\pd_j X})^2_k:L^2_{k+l}(E\ot L_{\la+\mu e_j}\vert_{\pd_j X})\longra L^2_k(F\ot L_{\la+\mu e_j}\vert_{\pd_j X})
\label{ac5eq18}
\e
is invertible for all\/ $i=1,\ldots,N$ and all\/ $\mu\in\R$. Whether \eq{ac5eq17}--\eq{ac5eq18} are Fredholm, or invertible, is independent of\/~$k$.
\item[{\bf(b)}] When $P^2_{k,\la}$ in \eq{ac5eq17} is Fredholm, there is an `Index Theorem' for\/ $\ind P^2_{k,\la}$ of the form \eq{ac5eq16}.
\end{itemize}
\label{ac5thm2}
\end{thm}

Loya \cite[Th.~1.1]{Loya2} also explains how the Fredholmness of \eq{ac5eq17} for {\it generic\/} $\la$ in $W(X)$ depends on invertibility of operators similar to \eq{ac5eq18} over~$C_2(X)$. 

Note that if \eq{ac5eq18} is invertible, then it is Fredholm, so we can apply Theorem \ref{ac5thm2} to $P_{\la+\mu e_j}\vert_{\pd_j X}$. Therefore questions about elliptic operators on manifolds with a-corners have an inductive flavour: it helps to understand the operators $P_\la\vert_{\pd^kX}$, or more-or-less equivalently $P_\la\vert_{C_k(X)}$, by reverse induction on the codimension~$k=\dim X,\dim X-1,\ldots,0$.

\subsection{B-fibrations and families of elliptic operators}
\label{ac54}

Definition \ref{ac4def6} defined {\it b-fibrations\/} $f:X\ra Y$ of manifolds with a-corners, following Melrose \cite{Melr1,Melr2,Melr4} for manifolds with ordinary corners. As we discuss in \S\ref{ac55}, b-fibrations play an important r\^ole in Melrose's b-calculus, in that the `Pushforward Theorem' Theorem \ref{ac5thm3}(b) works for b-fibrations. 

\begin{dfn} Let $X,Y$ be manifolds with a-corners with $\dim X\ge\dim Y$, and $f:X\ra Y$ be a proper b-fibration. Then ${}^b\d f:{}^bTX\ra f^*({}^bTY)$ is surjective by definition. Define the {\it relative tangent bundle\/ ${}^bT(X/Y)$ of\/} $f:X\ra Y$ to be $\Ker({}^b\d f)$. It is a vector subbundle of ${}^bTX$, of rank $\dim X-\dim Y$. 

If $y\in Y^\ci$ then $X_y=f^{-1}(y)$ is an embedded submanifold of $X$ of dimension $\dim X-\dim Y$, which is compact as $f$ is proper, and ${}^bT(X_y)\cong {}^bT(X/Y)\vert_{X_y}$. We loosely think of $f:X\ra Y$ as a {\it family of manifolds with a-corners $X_y$ over the base\/} $Y$, and of ${}^bT(X/Y)$ as the {\it family of b-tangent bundles\/}~${}^bT(X_y)$. 

If $y\in Y\sm Y^\ci$ then things are more complicated: the fibre $X_y=f^{-1}(y)$ may not be a submanifold, but only a union of submanifolds intersecting along boundary faces, and the topology of $X_y$ can change discontinuously as $y$ moves between different strata $S^k(Y)$ of $Y$. However, this is a feature, not a bug: as we explain in \S\ref{ac6}, this kind of topology change in b-fibrations is a very useful way of studying problems involving `gluing', `neck-stretching' and `bubbling'.

Now let $E,F\ra X$ be vector bundles and $P:\Ga^\iy(E)\ra\Ga^\iy(F)$ be a p.d.o.\ on $X$ of order $l$, as in Definition \ref{ac5def3}. Then after choosing b-connections ${}^b\nabla^E,{}^b\nabla$ on $E,{}^bTX$ we may write $P$ in the form \eq{ac5eq11} with $a_j\in \Ga^\iy\bigl(E^*\ot\bigot^j{}^bTX\ot F\bigr)$ for $j=0,\ldots,l$. Let ${}^b\nabla$ be chosen to preserve the subbundle ${}^bT(X/Y)\subseteq{}^bTX$. Then we call $P$ a {\it family of partial differential operators on\/ $X/Y$ of order\/} $l$ if 
\e
\ts a_j\in \Ga^\iy\bigl(E^*\ot\bigot^j{}^bT(X/Y)\ot F\bigr)\subseteq\Ga^\iy\bigl(E^*\ot\bigot^j{}^bTX\ot F\bigr),\;\> 0\le j\le l.
\label{ac5eq19}
\e

Generalizing Definition \ref{ac5def4}, for all $x\in X$ and $\xi\in {}^bT_x^*(X/Y)$, define a linear map $\si_{P/Y}(x,\xi):E_x\ra F_x$ by $\si_{P/Y}(x,\xi):e_x\mapsto a_l\vert_x\cdot (e_x\ot \bigot^l\xi)$, noting that $a_l\vert_x\in E^*_x\ot\bigot^l{}^bT_x(X/Y)\ot F_x$ by \eq{ac5eq19}. We call $\si_{P/Y}$ the {\it principal family symbol\/} of $P$. We call $P$ a {\it family of elliptic operators on\/} $X/Y$ if $\si_{P/Y}(x,\xi):E_x\ra F_x$ is an isomorphism for all $x\in X$ and~$0\ne\xi\in {}^bT_x^*(X/Y)$.

If $y\in Y^\ci$ then $P$ has a natural restriction $P\vert_{X_y}:\Ga^\iy(E\vert_{X_y})\ra\Ga^\iy(F\vert_{X_y})$, which is a p.d.o.\ of order $l$, and is elliptic if $P$ is an elliptic family. It is given by
\begin{equation*}
P\vert_{X_y}(e)=\ts\sum_{j=0}^l a_j\vert_{X_y}\cdot ({}^b\nabla_{X_y}^E)^je,
\end{equation*}
where we consider $a_j\vert_{X_y}$ as an element of $\Ga^\iy\bigl(E^*\vert_{X_y}\ot\bigot^j{}^bT(X_y)\ot F\vert_{X_y}\bigr)$ using ${}^bT(X_y)\cong {}^bT(X/Y)\vert_{X_y}$, and ${}^b\nabla_{X_y}^E$ are b-connections on $E\vert_{X_y}\ot\bigot^j{}^bT^*(X_y)$ constructed from the restrictions of ${}^b\nabla^E,{}^b\nabla$ to~$E\vert_{X_y},{}^bT(X/Y)\vert_{X_y}$.
\label{ac5def5}
\end{dfn}

Although we will not state any results to justify this, based on Melrose's theory, and on experience with `gluing' and `bubbling' problems in geometric analysis, we claim that families of elliptic operators over a proper b-fibration $f:X\ra Y$ in $\Manac$ are an interesting thing to study. 

If $P:\Ga^\iy(E)\ra\Ga^\iy(F)$ is a family of elliptic operators of order $l$ over a proper b-fibration $f:X\ra Y$, and $\la\in W(X)$, then for fixed $p>1$, $k\ge 0$ and all $y\in Y^\ci$ we have Banach space morphisms
\begin{equation*}
(P_\la\vert_{X_y})^p_k:L^p_{k+l}(E\ot L_\la\vert_{X_y})\longra L^p_k(F\ot L_\la\vert_{X_y}),
\end{equation*}
which as in \S\ref{ac53} are often Fredholm for suitable $\la$. Melrose's theory, and examples from the literature, suggest that we can understand these operators well as a family, including their asymptotic behaviour as $y$ approaches~$Y\sm Y^\ci$.

As we explain in \S\ref{ac6}, there are many important geometric problems involving `gluing', `bubbling' or `neck-stretching' which can be written in terms of families of {\it nonlinear\/} elliptic equations over a proper b-fibration $f:X\ra Y$ in~$\Manac$.

\subsection{\texorpdfstring{Relation to Melrose's b-calculus}{Relation to Melrose\textquoteright s b-calculus}}
\label{ac55}

We now discuss the theory of analysis on manifolds with corners developed by Richard Melrose and his collaborators \cite{Loya1,Loya2,LoMe,Melr1,Melr2,Melr3,Melr4,MeMe,MeNi,MePi,Piaz}, known as the `b-calculus'. For a survey we recommend Grieser \cite{Grie} and Melrose \cite{Melr1}. For manifolds with boundary the theory is developed in detail in Melrose's book \cite{Melr3}, and extensions to the corner case are discussed by Melrose, Nistor and Piazza~\cite{MeNi,MePi}.

Actually, as in \cite{Melr1} Melrose's framework applies to a larger class of theories; we will be concerned only with `b-problems', that is, those in which the initial Lie algebra of vector fields $\cV$ on a manifold with corners $X$ in \cite{Melr1} is~$\cV=\Ga^\iy({}^bTX)$.

Let $X$ be a compact manifold with corners, and $E,F\ra X$ be vector bundles. Then Melrose's theory defines a class $\Psi_b(X;E,F)$ of {\it b-pseudodifferential operators\/} from $E$ to $F$, which include {\it elliptic b-pseudodifferential operators}. The classical b-differential operators in $\Psi_b(X;E,F)$ are (modulo technicalities on the definition of a-smoothness) more-or-less the same as the differential operators discussed in \S\ref{ac53} on the corresponding manifold with a-corners $\ti X=F_\Manc^\Manac(X)$, and the notions of ellipticity in Melrose's theory and \S\ref{ac53} also agree.

B-pseudodifferential operators $P\in\Psi_b(X;E,F)$ act by convolution with a Schwartz kernel $K_P\in C^{-\iy}(\pi_1^*(E^*)\ot\pi_2^*(F))$ over $X\t X$, a distribution on $X\t X$ which is smooth in $(X^\ci\t X^\ci)\sm\De_X$, and has prescribed asymptotic behaviour near the boundary $\pd(X\t X)$ and diagonal $\De_X\subset X\t X$. Often one works on a blow-up $\,\,\,\,\,\widetilde{\!\!\!\!\!X\t X\!\!\!\!\!}\,\,\,\,\,$ of $X\t X$, to describe this asymptotic behaviour.

Given an elliptic b-pseudodifferential operator $P\in\Psi_b(X;E,F)$, one of the goals of the theory is to construct an approximate inverse (parametrix) $Q\in\Psi_b(X;F,E)$ for $P$, and understand the asymptotic behaviour of $K_Q$ at the boundary faces of $\,\,\,\,\,\widetilde{\!\!\!\!\!X\t X\!\!\!\!\!}\,\,\,\,\,$. This helps to understand the behaviour of solutions $e$ of $P(e)=f$, since $e\approx Q(f)$ up to smooth, fast-decaying errors.

Here is an important concept in Melrose's theory, \cite[\S 2.1.2]{Grie}, \cite[\S 4]{Melr4}:

\begin{dfn} Let $X$ be a manifold with corners, and $g:X^\ci\ra\R$ be a smooth function. We call $g$ {\it polyhomogeneous conormal\/} if whenever $x\in S^k(X)$ for $k>0$ and $(x_1,\ldots,x_m)\in\R^m_k$ are local coordinates on open $x\in U\subseteq X$ with $x=(0,\ldots,0)$, then near $x$ we have an asymptotic expansion
\e
g(x_1,\ldots,x_m)\sim \sum_{a_1,\ldots,a_k=0}^\iy\sum_{b_1=0}^{n_1^{a_1}}\cdots \sum_{b_k=0}^{n_k^{a_k}}\begin{aligned}[t]&G_{a_1\cdots a_k}^{b_1\cdots b_k}(x_{k+1},\ldots,x_m)\cdot{} \\
& x_1^{\al_1^{a_1}}\log^{b_1}x_1\cdots x_k^{\al_k^{a_k}}\log^{b_k}x_k
\end{aligned}
\label{ac5eq20}
\e
in $U^\ci$ as $x_1,\ldots,x_k\ra 0$, with all derivatives. Here for each $i=1,\ldots,k$, $(\al_i^j)_{j=0}^\iy$ is a sequence in $\C$ ordered such that $\Re\al_i^0\le \Re\al_i^1\le \cdots,$ with $\Re\al_i^j\ra\iy$ as $j\ra\iy$. Also $b_i^j\in\N$ for $i=1,\ldots,k$ and $j=0,1,\ldots,$ and $G_{a_1\cdots a_k}^{b_1\cdots b_k}$ is a smooth $\C$-valued function of $x_{k+1},\ldots,x_m$. The definition of `$\sim$' in \eq{ac5eq20} is explained in Melrose \cite[\S 4]{Melr2}, \cite[\S 5.10]{Melr3}, \cite[\S 4.13]{Melr4} and Grieser~\cite[\S 2.1.2]{Grie}.

The `index set' of $g$ assigns to each boundary face $\pd_iX$ of $X$ locally defined by $x_i=0$ for $x_i\in[0,\iy)$, the set $S_i$ of pairs $(\al,b)\in\C\t\N$ for which nonzero terms including the factor $x_i^\al\log^bx_i$ appear in the expansion \eq{ac5eq20} for $g$.
\label{ac5def6}
\end{dfn}

We define polyhomogeneous conormal sections of vector bundles in the same way. They are used by requiring the Schwartz kernels $K_P$ of $P\in\Psi_b(X;E,F)$ to be polyhomogeneous conormal near the boundary faces of $\,\,\,\,\,\widetilde{\!\!\!\!\!X\t X\!\!\!\!\!}\,\,\,\,\,$, and imposing conditions on the index sets for $K_P$.

Here are two important results for polyhomogeneous conormal functions \cite[Th.s 3 \& 5]{Melr2}, \cite[Th.s~3.10 \& 3.12]{Grie}, due to Melrose. In (b), a `b-density' is a section of $\La^{\rm top}({}^bT^*X)$, basically a volume form.

\begin{thm}{\bf(a) (Pullback Theorem)} Let\/ $f:X\ra Y$ be an interior map of manifolds with corners, and\/ $g:Y^\ci\ra\R$ be polyhomogeneous conormal on $Y$. Then the pullback\/ $f^*(g):X^\ci\ra\R$ is a well defined polyhomogeneous conormal function on $X,$ and its index set can be computed from that of\/~$g$.
\smallskip

\noindent{\bf(b) (Pushforward Theorem)} Let\/ $f:X\ra Y$ be a proper, cooriented b-fibration of manifolds with corners, and\/ $\de$ be a polyhomogeneous conormal b-density on $X$. Suppose the index set of\/ $\de$ satisfies $\Re\al_i>0$ for all\/ $(\al_i,b_i)$ in $S_i$ whenever $\pd_iX$ is a boundary face of\/ $X$ with\/ $f((\pd_iX)^\ci)\subseteq Y^\ci$. Then the pushforward\/ $f_*(\de)$ is a well defined polyhomogeneous conormal b-density on $Y,$ and its index set can be computed from that of\/~$\de$.
\label{ac5thm3}
\end{thm}

The Pullback and Pushforward Theorems are used in Melrose's theory for understanding operations on Schwarz kernels $K_P$ of $P\in\Psi_b(X;E,F)$ -- for instance, the kernel $K_{P\ci Q}$ of a composition $P\ci Q$ is determined from $K_P,K_Q$ using pullbacks and pushforwards along projections~$X\t X\t X\ra X\t X$.

The fact that the Pushforward Theorem holds for b-fibrations is relevant to \S\ref{ac54}. In particular, given a family of elliptic operators $P:\Ga^\iy(E)\ra\Ga^\iy(F)$ over a proper b-fibration $f:X\ra Y$, under good conditions Melrose's theory will construct a family of parametrices $Q$ for $P$ given by a polyhomogeneous conormal kernel on $X$ with known index set, and this index set will determine asymptotic features of the equations $P\vert_{X_y}(e)=f$ for $y\in Y^\ci$ as~$y\ra\pd Y$.

Now we will relate Melrose's theory to our manifolds with a-corners. The next lemma relates polyhomogeneous conormal functions with a-smooth sections of $L_\la$. It can be proved by comparing Definitions \ref{ac3def2} and~\ref{ac5def6}. 

\begin{lem} Let\/ $X$ be a manifold with corners, and\/ $g:X^\ci\ra\R$ a polyhomogeneous conormal function. Write\/ $\ti X=F_\Manc^\Manac(X)$ for the corresponding (untwisted) manifold with a-corners, and let\/ $\la\in W(\ti X)$. Suppose that for each connected component\/ $\pd_iX$ of\/ $\pd X$ defined locally by $x_i=0$ for $x_i\in[0,\iy),$ so that\/ $\la=\la_i\cdot x_i\frac{\pd}{\pd x_i}$ on $\pd_i\ti X$ for some $\la_i\in\R,$ we have either $\la_i=\al$ and\/ $b=0$ or $\la_i<\Re\al$ for all\/ $(\al,b)$ in the corresponding indexing set\/ $S_i$ of\/ $g$. Then there is a unique a-smooth section $\ti g\in\Ga^\iy(L_\la)$ with\/ $\ti g\vert_{\ti X^\ci}=g,$ noting that\/~$\ti X^\ci=X^\ci$.

\label{ac5lem2}
\end{lem}

Note that any $\la\in W(\ti X)$ which is `sufficiently negative' satisfies the conditions in the lemma, so every polyhomogeneous conormal function on $X^\ci$ extends to an a-smooth section of $L_\la\ra\ti X$ for some $\la\in W(\ti X)$. The converse of Lemma \ref{ac5lem2} is false: an a-smooth section of $L_\la$ need not have an asymptotic expansion \eq{ac5eq20}, so a-smoothness is weaker than polyhomogeneous conormality.

The author proposes to modify Melrose's b-calculus as follows:
\begin{itemize}
\setlength{\itemsep}{0pt}
\setlength{\parsep}{0pt}
\item[(i)] Manifolds with corners should be replaced by manifolds with a-corners (often supposed untwisted) throughout.
\item[(ii)] Polyhomogeneous conormal functions, or polyhomogeneous conormal sections of a vector bundle $E\ra X$, should be replaced by a-smooth sections of $L_\la$ or $E\ot L_\la$ for $\la\in W(X)$.
\end{itemize}

One difference between polyhomogeneous conormal functions $g:X^\ci\ra\R$, as in Melrose's theory, and a-smooth sections $\ti g$ of $L_\la$, as in our proposal, is that knowing $\ti g\in\Ga^\iy(L_\la)$ gives you a good understanding of the leading term $G_{0\cdots 0}^{0\cdots 0}(x_{k+1},\ldots,x_m) x_1^{\la_1}\cdots x_k^{\la_k}$ in the asymptotic expansion \eq{ac5eq20} --- the function $G_{0\cdots 0}^{0\cdots 0}$ corresponds to $\ti g\vert_{S^k(X)}$ --- but tells you less about subsequent terms. 

Melrose \cite[\S 5.16]{Melr3} also discusses functions satisfying `conormal bounds' like \eq{ac3eq3}, but without an asymptotic expansion \eq{ac5eq20}, and these seem close to our a-smooth sections of $L_\la$. Here is the analogue of Theorem \ref{ac5thm3}. It can be proved using similar methods to Melrose~\cite{Melr2}.

\begin{thm}{\bf(a) (Pullback Theorem)} Let\/ $f:X\ra Y$ be an interior map of manifolds with a-corners, and\/ $\la\in W(Y),$ and\/ $g\in\Ga^\iy(L_\la)$. Then Definition\/ {\rm\ref{ac4def12}} gives a weight\/ $f^*(\la)\in W(X)$ with a canonical isomorphism $L_{f^*(\la)}\cong f^*(L_\la),$ and using this we have\/~$f^*(g)\in\Ga^\iy(L_{f^*(\la)})$.
\smallskip

\noindent{\bf(b) (Pushforward Theorem)} Let\/ $f:X\ra Y$ be a proper, cooriented b-fibration of manifolds with a-corners, and\/ $\la\in W(X)$ with\/ $\la\vert_{\pd_iX}>0$ whenever $\pd_iX$ is a boundary face of\/ $X$ with\/ $f((\pd_iX)^\ci)\subseteq Y^\ci$. Then there is a unique weight\/ $\mu\in W(Y)$ with the universal property that\/ $f^*(\mu)\le\la,$ and for any $\mu'\in W(Y)$ with\/ $f^*(\mu')\le\la$ we have $\mu'\le\mu$. For any\/ $\de$ in $\Ga^\iy(\La^{\rm top}({}^bT^*X)\ot L_\la)$ the pushforward\/ $f_*(\de)$ exists in\/~$\Ga^\iy(\La^{\rm top}({}^bT^*Y)\ot L_\mu)$.
\label{ac5thm4}
\end{thm}

The fact that the Pullback and Pushforward Theorems hold means that much of Melrose's b-calculus transfers to our framework with little extra work.

\section{How to use manifolds with a-corners}
\label{ac6}

We now briefly outline a number of topics of current research in geometry which can be rewritten using the language of manifolds with a-corners. Our choice is strongly biased towards the author's interests. Sometimes we believe the `a-corners' point of view is particularly helpful, or opens up new possibilities or new methods of proof, and we discuss these in more detail.

\subsection{Asymptotic conditions on manifolds and submanifolds}
\label{ac61}

First we consider noncompact Riemannian manifolds $(X,g)$ and closed submanifolds $Y\subset X$ satisfying asymptotic conditions at the noncompact ends of $X$. We also suppose $(X,g)$ and $Y$ have some interesting property, for example $g$ could be Ricci-flat or have special holonomy, and $Y$ might be minimal or calibrated. Often this means that $g$ and $Y$ satisfy nonlinear elliptic p.d.e.s.

The idea is to write $X=\ti X^\ci$ for a compact manifold with a-corners $\ti X$, and $g=\ti g\vert_X$ for a weighted b-metric $\ti g$ on $\ti X$ (or something similar), and $Y=\ti Y\cap X$ for $\ti Y\subseteq\ti X$ a compact submanifold with a-corners. Then the asymptotic conditions on $g$ and $Y$ are encoded in the differential geometry of $\ti X,\ti Y$, and the elliptic p.d.e.s satisfied by $g,Y$ extend to elliptic p.d.e.s on $\ti g,\ti Y$ as in~\S\ref{ac53}.

See \cite{Joyc2,Joyc9} for background on holonomy groups and calibrated geometry.

\subsubsection{Asymptotically Cylindrical manifolds and submanifolds}
\label{ac611}

As in Example \ref{ac5ex1}(a), a Riemannian $m$-manifold $(X,g)$ is {\it Asymptotically Cylindrical\/} ({\it ACy\/}) if it has one noncompact end asymptotic as $r\ra\iy$ to a Riemannian cylinder $(\R\t M,h_0+(\d r)^2)$, with $O(e^{-\al r})$ decay in all derivatives for small $\al>0$, where $(M,h_0)$ is a compact Riemannian $(m\!-\!1)$-manifold, the {\it link}. Given an ACy manifold $(X,g)$, we can define {\it Asymptotically Cylindrical submanifolds\/} $Y\subset X$, asymptotic as $r\ra\iy$ to $\R\t N$ for $N\subset M$ a closed submanifold.

{\it Calabi--Yau manifolds\/} are Riemannian $2m$-manifolds $(X,g)$ with holonomy $\SU(m)$. They have additional geometric structures, the complex structure $J\in\Ga^\iy(TX\ot T^*X)$, K\"ahler form $\om\in\Ga^\iy(\La^2T^*X)$, and holomorphic volume form $\Om\in\Ga^\iy(\La^mT^*X\ot_\R\C)$. ACy Calabi--Yau manifolds are studied in detail by Kovalev \cite{Kova}, Corti--Haskins--Nordstr\"om--Pacini \cite{CHNP1,CHNP2}, Haskins--Hein--Nordstr\"om \cite{HHN}, and Conlon--Mazzeo--Rochon \cite{CMR}. They are used to construct compact 7-manifolds with holonomy $G_2$ by Kovalev \cite{Kova}, as in~\S\ref{ac621}.

The asymptotic conditions on $g$ are exactly those in equation \eq{ac5eq6} of Example \ref{ac5ex1}(a), with all derivatives. Therefore $X$ extends to a compact manifold with a-boundary $\ti X$ with $\ti X^\ci=X$ and $\pd\ti X\cong M$, and $g=\ti g\vert_X$ for a b-metric $\ti g$ on $\ti X$. Also $J,\om,\Om$ extend to $\ti J\in\Ga^\iy({}^bT\ti X\ot{}^bT^*\ti X)$, $\ti\om\in\Ga^\iy(\La^2({}^bT^*\ti X))$ and~$\ti\Om\in\Ga^\iy(\La^m({}^bT^*\ti X)\ot_\R\C)$.

Kovalev and Nordstr\"om \cite{KoNo,Nord} study ACy 7-manifolds $(X,g)$ with holonomy $G_2$. Again, we can write $X=\ti X^\ci$ and $g=\ti g\vert_X$ for $\ti X$ a compact 7-manifold with a-boundary and $\ti g$ a b-metric on $\ti X$. ACy coassociative 4-folds $Y\subset X$ in an ACy 7-manifold $(X,g)$ with holonomy contained in $G_2$ are considered by Salur and the author \cite{JoSa}. (Here {\it coassociative submanifolds\/} are a kind of calibrated submanifold in 7-manifolds with holonomy in $G_2$, see \cite{HaLa,Joyc9}.) We have $Y=\ti Y\cap X$ for a compact submanifold with a-boundary~$\ti Y\subset\ti X$.

\subsubsection{Asymptotically Conical manifolds and submanifolds}
\label{ac612}

As in Example \ref{ac5ex1}(b), a Riemannian $m$-manifold $(X,g)$ is {\it Asymptotically Conical\/} ({\it ACo\/}) if it has one noncompact end asymptotic as $r\ra\iy$ to a Riemannian cone $((0,\iy)\t M,r^2h_0+(\d r)^2)$, with $O(r^{-\al-k})$ decay in $k^{\rm th}$ derivatives for small $\al>0$, where $(M,h_0)$ is a compact Riemannian $(m-1)$-manifold. Given an ACo manifold $(X,g)$, we can define {\it Asymptotically Conical submanifolds\/} $Y\subset X$, asymptotic as $r\ra\iy$ to $(0,\iy)\t N$ for $N\subset M$ a closed submanifold.

An ACo manifold $(X,g)$ is called {\it Asymptotically Euclidean\/} ({\it AE\/}) if $M=\cS^{m-1}$ with the round metric, so that $X$ is asymptotic at infinity to $\R^m$ with the Euclidean metric, and {\it Asymptotically Locally Euclidean\/} ({\it ALE\/}) if $M=\cS^{m-1}/\Ga$ for finite $\Ga\subset{\rm O}(m)$ acting freely on $\cS^{m-1}$, so that $X$ is asymptotic to~$\R^m/\Ga$.

ACo Calabi--Yau manifolds are studied by Kronheimer \cite{Kron1,Kron2} and the author \cite[\S 8]{Joyc2}, \cite{Joyc3} in the ALE case, and by Conlon and Hein \cite{CoHe1,CoHe2} in the general case. Examples of ACo Riemannian 7- and 8-manifolds $(X,g)$ with holonomy $G_2$ and $\Spin(7)$ are constructed by Bryant and Salamon \cite{BrSa}. In both cases, as in Example \ref{ac5ex1}(b) we can write $X=\ti X^\ci$ and $g=\ti g\vert_X$ for $\ti X$ a compact manifold with a-boundary and $\ti g$ a weighted b-metric on $\ti X$, with weight~$\la>0$.

{\it Special Lagrangian submanifolds\/} ({\it SL\/ $m$-folds\/}) in Calabi--Yau $m$-folds are a kind of calibrated submanifold, as in \cite{HaLa,Joyc9}. ACo SL $m$-folds in $\C^m$ are studied by Marshall \cite{Mars} and Pacini \cite{Paci1}, who studied their deformation theory, and examples can be found in Harvey and Lawson \cite[\S III.3]{HaLa} and the author \cite{Joyc5,Joyc6}. ACo coassociative 4-folds $Y$ in $\R^7$ are studied by Lotay \cite{Lota2}. In both cases we can write $\C^m,\R^7$ as the interior of compact disc $\ti X$ with a-boundary, and then $Y=\ti Y\cap\C^m$ or $Y=\ti Y\cap\R^7$ for a compact submanifold with a-boundary~$\ti Y\subset\ti X$.

\subsubsection{Manifolds and submanifolds with conical singularities}
\label{ac613}

A Riemannian manifold $(X,g)$ has a {\it conical singularity\/} if it has a noncompact end asymptotic as $r\ra 0$ to a Riemannian cone $((0,\iy)\t M,r^2h_0+(\d r)^2)$, for $(M,h_0)$ a compact Riemannian manifold. Often one compactifies the noncompact end by adding a single point $\{x_0\}$, giving a singular manifold $\hat X=X\amalg\{x_0\}$. An alternative way to compactify the noncompact end is as a manifold with a-boundary $\ti X$, with $\ti X^\ci=X$ and $\pd\ti X\cong M$. Then as in Example \ref{ac5ex1}(c) we have $g=\ti g\vert_X$ for $\ti g$ a weighted b-metric on $\ti X$, with weight~$\la<0$.

Calabi--Yau 3-folds with conical singularities are studied by Chan \cite{Chan1}. SL $m$-folds $Y\subset X$ with conical singularities are considered in Calabi--Yau $m$-folds by the author \cite{Joyc7,Joyc8}, in $\C^m$ by Pacini \cite{Paci2,Paci3}, and in Calabi--Yau $m$-folds with conical singularities by Chan \cite{Chan2}. Coassociative 4-folds with conical singularities in $G_2$-manifolds are studied by Lotay \cite{Lota1,Lota3,Lota4}. In each case we can extend $Y$ to a submanifold with a-boundary $\ti Y$ in $X$ or~$\ti X$.

\subsubsection{Quasi-ALE and Quasi-Asymptotically Conical manifolds}
\label{ac614}

A Riemannian manifold $(X,g)$ is {\it Quasi-ALE\/} $(X,g)$ if it has one noncompact end asymptotic at infinity (in a sense) to $\R^m/\Ga$, where $\Ga\subseteq\mathbin{\rm O}(m)$ is a finite subgroup, and the singularities of $\R^m/\Ga$ may extend to infinity. The asymptotic conditions on $(X,g)$ are complicated, and inductive on dimension: if $\R^m/\Ga\simeq\R^k\t(\R^{m-k}/\De)$ locally near infinity in $\R^m/\Ga$, then $(X,g)\simeq(\R^k,g_0)\t(Y,h)$ for $(Y,h)$ a Quasi-ALE manifold asymptotic to~$\R^{m-k}/\De$.

Quasi-ALE manifolds were introduced by myself \cite[\S 8]{Joyc2}, \cite{Joyc4} to describe the natural K\"ahler metrics on resolutions $X\ra\C^m/\Ga$. I proved a Calabi Conjecture for such metrics, and so constructed many Quasi-ALE Calabi--Yau manifolds. I also found examples of Quasi-ALE $G_2$-manifolds \cite[\S 11.2]{Joyc2}. Quasi-ALE manifolds have since been studied by Carron \cite{Carr1,Carr2}. Degeratu and Mazzeo \cite{DeMa,Mazz} generalize Quasi-ALE manifolds to {\it Quasi-Asymptotically Conical\/} ({\it Quasi-ACo\/}) manifolds, which are to Quasi-ALE manifolds as ACo manifolds are to ALE manifolds. They prove Fredholm results for Laplacian-type operators on weighted Sobolev and H\"older spaces on Quasi-ACo manifolds.

Given a Quasi-ACo manifold $(X,g)$, there is a way to find a (possibly nonunique) compact manifold with a-corners $\ti X$ with $X=\ti X^\ci$ and a vector bundle $E\ra\ti X$ with $E\vert_{X}=TX$, such that $g=\ti g\vert_X$ for some positive definite $\ti g\in\Ga^\iy(S^2E^*)$. In general we have $\pd^2\ti X\ne\es$. If we had $E={}^bT\ti X\ot L_\la$ for $\la\in W(\ti X)$ then $\ti g$ would be a weighted b-metric, as in \S\ref{ac51}. However, things may be more complicated than this, for example, we may have ${}^bT\ti X=F_1\op F_2$ and $E=F_1\ot L_{\la_1}\op F_2\ot L_{\la_2}$ for $\la_1\ne\la_2$ in $W(\ti X)$. The Quasi-ACo conditions on $(X,g)$ are then encoded in the differential geometry of~$\ti X,E$.

Using manifolds with a-corners may help in this situation. For example, the Fredholm results on weighted Sobolev and H\"older spaces in \cite[\S 8]{Joyc2}, \cite{Joyc4} and \cite{DeMa} may be consequences of general facts about analysis of elliptic operators on compact manifolds with a-corners, as in~\S\ref{ac53}. 

\subsection{\texorpdfstring{Family problems involving `stretching necks'}{Family problems involving \textquoteleft stretching necks\textquoteright}}
\label{ac62}

Next we discuss geometric problems involving a family $\{(X_y,e_y):y\in Y^\ci\}$ of solutions $e_y$ of nonlinear elliptic p.d.e.s on compact manifolds $X_y$, parametrized by the interior $Y^\ci$ of a manifold with (a-)corners $Y$, such that as $y$ approaches $\pd Y$, the solution $(X_y,e_y)$ undergoes an asymptotic decay, often described as `neck stretching', `neck pinching', `bubbling', or `gluing', and for $\bar y\in Y\sm Y^\ci$, we have a `singular solution' $(X_{\bar y},e_{\bar y})$, with different topology to $(X_y,e_y)$ for~$y\in Y^\ci$.

In the language of manifolds with a-corners, these problems should be written in terms of a family of solutions $e$ of a family of nonlinear elliptic p.d.e.s over a proper b-fibration $f:X\ra Y$ of manifolds with a-corners, as in~\S\ref{ac54}.

Problems of this kind occur in the literature in two main ways:
\begin{itemize}
\setlength{\itemsep}{0pt}
\setlength{\parsep}{0pt}
\item[(a)] If we wish to construct examples of nonsingular solutions of a nonlinear elliptic p.d.e., one method is to start with a singular solution (which may be simpler), and deform it to a nonsingular solution. This corresponds to starting with $(X_{\bar y},e_{\bar y})$ for $\bar y\in Y\sm Y^\ci$, and constructing $(X_y,e_y)$ for $y\in Y^\ci$ close to $\bar y$. It is usually enough to take $Y=\lb 0,\ep)$ for small~$\ep>0$. 

As explained in \S\ref{ac621}, this method is used to construct examples of compact 7-manifolds with holonomy $G_2$ in \cite{CHNP1,Kova}, and compact 8-manifolds with holonomy $\Spin(7)$ in \cite{Joyc1}, \cite[\S 15]{Joyc2}, and compact SL $m$-folds in~\cite{Joyc7,Joyc8}.
\item[(b)] Consider the moduli space $\cM$ of solutions of some nonlinear elliptic p.d.e.\ on a compact manifold. We would like to form and study a compactification $\oM$, where points of the boundary $\oM\sm\cM$ correspond to some kind of singular solutions. This occurs, for example, for `breaking' of Morse flow-lines as in \S\ref{ac622}, for `bubbling' of moduli spaces of instantons on 4-manifolds in Donaldson theory \cite{DoKr}, and for moduli spaces of stable $J$-holomorphic curves in symplectic geometry, as in~\S\ref{ac63}.
\end{itemize}

The `a-corners' point of view may help in such situations:
\begin{itemize}
\setlength{\itemsep}{0pt}
\setlength{\parsep}{0pt}
\item[(i)] There may be general results on families of elliptic operators over a b-fibration, which could be applied in many different problems. 
\item[(ii)] We argue that the `natural' smooth structure on moduli spaces $\oM$ comes from manifolds with a-corners $Y$ in the b-fibrations~$f:X\ra Y$.
\item[(iii)] In \S\ref{ac634} we advocate studying such moduli problems using Grothendieck's method of representable functors, and universal families. This depends crucially on good properties of the category~$\Manac$.
\end{itemize}

\subsubsection{Some examples of gluing constructions}
\label{ac621}

In the next example we describe a construction of compact 7-manifolds with holonomy $G_2$ proposed by Donaldson (see \cite[\S 11.9]{Joyc2}), carried out by Kovalev \cite{Kova}, and developed further by Kovalev and Nordstr\"om \cite{KoNo} and Corti, Haskins, Nordstr\"om, and Pacini \cite{CHNP1,CHNP2}. We explain how to write the construction in terms of b-fibrations of manifolds with a-corners.

\begin{ex} See the author \cite{Joyc2,Joyc9} for background on $G_2$. A $G_2$-{\it manifold\/} $(X,\vp,g)$ is a 7-manifold $X$ with a torsion-free $G_2$-structure $(\vp,g)$, where $g$ is a Riemannian metric on $X$, and $\vp$ is a 3-form on $X$ with $\d\vp=\d^*\vp=0$ satisfying a pointwise compatibility with $g$. These imply that $\Hol(g)\subseteq G_2$, and $(\vp,g)$ determine a principal $G_2$-subbundle $P\subset F$ of the frame bundle $F$ of~$X$.

Let $(\R\t N,\vp_{\rm cyl},g_{\rm cyl})$ be a cylindrical $G_2$-manifold with $g_{\rm cyl}=g_N+\d r^2$, where $(N,g_N)$ is a compact Riemannian 6-manifold (usually $K3\t T^2$), and $r$ the coordinate on $\R$. Suppose $(X_0^+,\vp_0^+,g_0^+)$ and $(X_0^-,\vp_0^-,g_0^-)$ are Asymptotically Cylindrical (ACy) $G_2$-manifolds, as in \S\ref{ac611}, with $(X_0^+,\vp_0^+,g_0^+)$ asymptotic to $(\R\t N,\vp_{\rm cyl},g_{\rm cyl})$ as $r\ra\iy$, and $(X_0^-,\vp_0^-,g_0^-)$ asymptotic to $(\R\t N,\vp_{\rm cyl},g_{\rm cyl})$ as $r\ra-\iy$. Then we can choose compact $K_0^\pm\subset X_0^\pm$ and diffeomorphisms
\e
\begin{aligned}
&X_0^+\sm K_0^+\cong (0,\iy)\t N,\quad X_0^-\sm K_0^-\cong (-\iy,0)\t N,\quad\text{such that}\\ 
&\nabla^k(\vp_0^+-\vp_{\rm cyl},g_0^+-g_{\rm cyl})=O(e^{-\al r})\;\>
\text{as $r\ra\iy$, $k\ge 0,$ and}\\
&\nabla^k(\vp_0^--\vp_{\rm cyl},g_0^--g_{\rm cyl})=O(e^{\al r})\;\>\text{as $r\ra-\iy$, $k\ge 0,$ for small $\al>0$.}
\end{aligned}
\label{ac6eq1}
\e

Let $\ep\in(0,1)$ be small, and $y\in(0,\ep)$. Make a compact 7-manifold $X_y$ by
\begin{equation*}
X_y=\bigl(X_0^+\sm (-\log y,\iy)\t N\bigr)\amalg_{\{-\log y\}\t N=\{\log y\}\t N}\bigl(X_0^-\sm (-\iy,\log y)\t N\bigr).
\end{equation*}
That is, we cut off the cylindrical ends of $X_0^+$ at $r=-\log y$ and of $X_0^-$ at $r=\log y$, and glue the boundary hypersurfaces $\{-\log y\}\t N$ and $\{\log y\}\t N$ together in the obvious way.

Next we glue the $G_2$-structures $(\vp_0^+,g_0^+),(\vp_0^-,g_0^-),(\vp_{\rm cyl},g_{\rm cyl})$ together using a partition of unity to get a $G_2$-structure with torsion $(\ti\vp_y,\ti g_y)$ on $X_y$. Using \eq{ac6eq1} and $r=\pm\log y$ we see that the torsion of $(\ti\vp_y,\ti g_y)$ satisfies
\begin{equation*}
\bmd{\ti\nabla\ti\vp_y}_{\ti g_y}=O(y^\al).
\end{equation*}
By solving a nonlinear elliptic p.d.e., Kovalev \cite{Kova} shows that for small enough $y$, we can find a torsion-free $G_2$-structure $(\vp_y,g_y)$ on $X_y$ close to $(\ti\vp_y,\ti g_y)$, with
\begin{equation*}
\nabla^k(\vp_y-\ti\vp_y,g_y-\ti g_y)=O(y^\al).
\end{equation*}

If $\pi_1(X_y)$ is finite then $\Hol(g_y)=G_2$, so that $(X_y,g_y)$ is a compact Riemannian 7-manifold with holonomy $G_2$. Kovalev \cite{Kova} and Corti, Haskins, Nordstr\"om, and Pacini \cite{CHNP1,CHNP2} use this to make many new examples of compact 7-manifolds with holonomy $G_2$, in which $X_0^\pm$ are of the form $W^\pm_0\t\cS^1$, for $W_0^\pm$ ACy Calabi--Yau 3-folds, as in \S\ref{ac611}. Kovalev and Nordstr\"om \cite{KoNo} also give examples in which $X_0^\pm$ have holonomy~$G_2$.

Here is how to interpret all this in the language of manifolds with a-corners. We should define an 8-manifold with a-corners $X$ and a proper b-fibration $f:X\ra Y=\lb 0,\ep)$, where $f^{-1}(y)=X_y$ for $y\in(0,\ep)$, and $f^{-1}(0)=X_0^+\amalg N\amalg X_0^-$, where $X_0^\pm$ are two a-boundary faces of $X$, meeting in the codimension 2 a-corner $N$ of $X$. Also $X$ near $N$ is a-diffeomorphic to $\lb 0,\iy)^2\t N$ near $\{(0,0)\}\t N$, with $f$ identified with $(s,t,n)\mapsto st$. We illustrate this in Figure~\ref{ac6fig1}.
\begin{figure}[htb]
\centerline{$\splinetolerance{.8pt}
\begin{xy}
0;<1.2mm,0mm>:
,(-10,15)*{\bu}
,(-10,-15)*{\bu}
,(-30,0)*{\bu}
,(20,15)*{\ci}
,(20,-15)*{\ci}
,(-10,-23)*{\bu}
,(20,-23)*{\ci}
,(-1,-19)*{f\downarrow}
,(8,-17)*{X}
,(8,-21)*{Y}
,(-12,-23)*{0}
,(22,-23)*{\ep}
,(-9.6,15);(19.5,15)**\crv{(5,15)}
,(-9.6,-15);(19.5,-15)**\crv{(5,-15)}
,(-9.6,-23);(19.5,-23)**\crv{(5,-23)}
,(-10,15);(-30,0)**\crv{}
,(-10,-15);(-30,0)**\crv{}
,(19.7,14.4);(19.7,-14.4)**\crv{~**\dir{--}(10,0)}
,(15,15);(15,-15)**\crv{~*{.} (3,0)}
,(10,15);(10,-15)**\crv{~*{.} (-6,0)}
,(5,15);(5,-15)**\crv{~*{.} (-15,0)}
,(0,15);(0,-15)**\crv{~*{.} (-23,0)}
,(-5,15);(-5,-15)**\crv{~*{.}(-36,1),(-50,0),(-36,-1)}
,(-29,10)*{\text{a-boundary $X_0^+$}}
,(-29,-10)*{\text{a-boundary $X_0^-$}}
,(-38,0)*{\text{a-corner $N$}}
,(-17.7,0)*{X_{y_1}}
,(-9.4,0)*{X_{y_2}}
,(-2.2,0)*{X_{y_3}}
,(4.9,0)*{X_{y_4}}
,(11.9,0)*{X_{y_5}}
,(17.2,0)*{X_\ep}
,(-5,-23)*{\bu}
,(0,-23)*{\bu}
,(5,-23)*{\bu}
,(10,-23)*{\bu}
,(15,-23)*{\bu}
,(-5,-25.5)*{y_1}
,(0,-25.5)*{y_2}
,(5,-25.5)*{y_3}
,(10,-25.5)*{y_4}
,(15,-25.5)*{y_5}
\end{xy}$}
\caption{B-fibration $f:X\ra\lb 0,\ep)$ in Kovalev's $G_2$-manifold construction}
\label{ac6fig1}
\end{figure}

There should exist a b-metric $g$ on $X$, in the sense of \S\ref{ac51}, roughly $g=g_y+y^{-2}\d y^2$, with $b\vert_{X_y}=g_y$ for $y\in(0,\ep)$, and $g\vert_{X_0^\pm}=g_0^\pm$, and $g\vert_N=g_N$. Using the notation of \S\ref{ac54}, there should exist a-smooth $\vp\in\Ga^\iy\bigl(\La^3({}^bT^*(X/Y))\bigr)$ with $\vp\vert_{X_y}=\vp_y$ for $y\in(0,\ep)$, and $\vp\vert_{X_0^\pm}=\vp_0^\pm$. The $(\vp_t,g_t)$ for $t\in(0,\ep)$, and $(\vp_0^\pm,g_0^\pm)$ for $t=0$, are an a-smooth family of solutions of a family of nonlinear elliptic p.d.e.s over the proper b-fibration $f:X\ra Y$, in the sense of~\S\ref{ac54}.

This is an example of the following general question. Let $f:X\ra Y$ be a proper b-fibration in $\Manac$, and suppose we are given a family of nonlinear elliptic p.d.e.s $P$ on sections of $E\ra X$ over $f$, a point $y_0\in Y$, and $e_{y_0}$ in $\Ga^\iy(E\vert_{X_{y_0}})$ with $P\vert_{X_{y_0}}(e_{y_0})=0$. We want to know whether there exist $e_y$ in $\Ga^\iy(E\vert_{X_y})$ close to $e_{y_0}$ with $P\vert_{X_y}(e_y)=0$ for all $y\in Y$ sufficiently close to~$y_0$.

If $y_0$ were a nonsingular point of the b-fibration, this would be a standard question in deformation theory, and the answer would be yes provided the cokernel of the (Fredholm) linearization of $P\vert_{X_{y_0}}$ at $e_{y_0}$ is zero. 

In our case, $y_0=0$ in $Y=\lb 0,\ep)$ is a singular point of the b-fibration. But the author expects that in the `a-corners' theory, we should often be able to treat {\it all\/} fibres of a b-fibration in a uniform way, not just the nonsingular fibres, and that the Kovalev gluing construction should follow from the general deformation theory of families of nonlinear elliptic p.d.e.s over proper b-fibrations.
\label{ac6ex1}
\end{ex}

Here are some more examples of gluing constructions which work in a similar way, with the difference that rather than gluing together two Asymptotically Cylindrical objects, as in \S\ref{ac611}, they glue an Asymptotically Conical object $X_0^+$ (shrunk by a conformal factor) into an object with a conical singularity $X_0^-$, as in \S\ref{ac612}--\S\ref{ac613}. They can also be interpreted using b-fibrations of manifolds with a-corners. The shrinking by a conformal factor is achieved by twisting by a line bundle $L_\la$ from \S\ref{ac44}, with $\la$ nonzero on the a-boundary component $X_0^+$.
\begin{itemize}
\setlength{\itemsep}{0pt}
\setlength{\parsep}{0pt}
\item The author \cite{Joyc1}, \cite[\S 15]{Joyc2} constructs compact 8-manifolds with holonomy $\Spin(7)$ by gluing ALE $\Spin(7)$-manifolds into a $\Spin(7)$-orbifold.
\item The author \cite{Joyc7,Joyc8} and Pacini \cite{Paci3} construct SL $m$-folds in Calabi--Yau $m$-folds by gluing ACo SL $m$-folds into SL $m$-folds with conical singularities.
\item Lotay \cite{Lota3,Lota4} glues coassociative 4-folds in $G_2$-manifolds in the same way.
\end{itemize}

\subsubsection{Morse homology and moduli spaces of Morse flow-lines}
\label{ac622}

We recall some background from Morse homology, as in Schwarz \cite{Schw} or Austin and Braam \cite{AuBr}. Let $W$ be a compact manifold, and $f:W\ra\R$ be smooth. We call $f$ {\it Morse\/} if the critical locus $\Crit(f)$ consists of isolated points $w\in W$ with Hessian $\Hess f=\nabla^2f\in S^2T^*_wW$ nondegenerate. We call $f$ {\it Morse--Bott\/} if $\Crit(f)$ is a disjoint union of submanifolds $S\subset W$ with $\Hess_wf$ nondegenerate on the normal space $N_wS\subseteq T_wW$ for each~$w\in S$.

If $f$ is Morse--Bott, write $C$ for the set of connected components of $\Crit(f)$, and define $\mu:C\ra\N$ by $\mu(S)=k$ if $\Hess f$ has $k$ negative eigenvalues on $N_wS$ for each $w\in S$. Fix a Riemannian metric $g$ on $W$, and let $\nabla f\in\Ga^\iy(TW)$ be the associated gradient vector field of $W$, so that $\nabla^af=g^{ab}\d_bf$ in index notation.

For $S,T\in C$, consider the moduli space $\cM(S,T)$ of $\sim$-equivalence classes $[\ga]$ of smooth flow-lines $\ga:\R\ra W$ of $-\nabla f$ (that is, $\frac{\d}{\d t}\ga(t)=-\nabla f\vert_{\ga(t)}$ for $t\in\R$) with $\lim_{t\ra-\iy}\ga(t)\in S$ and $\lim_{t\ra\iy}\ga(t)\in T$, with equivalence relation $\ga\sim\ga'$ if $\ga(t)=\ga'(t+c)$ for some $c\in\R$ and all $t$. Define {\it evaluation maps\/} $\ev_-:\cM(S,T)\ra S$, $\ev_+:\cM(S,T)\ra T$ by $\ev_\pm:[\ga]\mapsto\lim_{t\ra\pm\iy}\ga(t)$.

Next, enlarge $\cM(S,T)$ to a compactification $\oM(S,T)$ by including `broken flow-lines', which are sequences $([\ga_1],\ldots,[\ga_k])$ for $[\ga_i]\in\cM(U_i,U_{i+1})$ with $U_i\in C$, $U_0=S$, $U_k=T$ and $\ev_+([\ga_i])=\ev_-([\ga_{i+1}])\in U_i$ for~$1\le i<k$.

As discussed by Wehrheim \cite[\S 1]{Wehr}, there is a widely-believed folklore result, which we state ({\it without\/} claiming it is true) in a particularly strong form:
\smallskip

\noindent{\bf `Folklore Theorem'.} {\it Let\/ $W$ be a compact manifold, $f:W\ra\R$ a Morse--Bott function, and\/ $g$ a generic Riemannian metric on $W$. Then with the notation above, $\oM(S,T)$ has the canonical structure of a compact manifold with corners of dimension $\dim T+\mu(T)-\mu(S)-1$ for all\/ $S,T\in J,$ with interior $\oM(S,T)^\ci=\cM(S,T),$ where $\oM(S,T)=\es$ if\/ $\mu(S)\ge\dim T+\mu(T)$. The evaluation maps $\ev_-:\oM(S,T)\ra S,$ $\ev_+:\oM(S,T)\ra T$ are smooth. There is a canonical diffeomorphism, where the fibre products on the right hand side are transverse:
\begin{equation*}
\pd\oM(S,T)\cong \ts\coprod_{U\in J}\oM(S,U)\t_{\ev_+,U,\ev_-}\oM(U,T),
\end{equation*}
compatible with the maps from both sides to $S,T$ defined using $\ev_\mp$.}
\smallskip

One can find a claim of this kind for Morse functions in Austin and Braam \cite[p.~130]{AuBr}, for instance. However, references such as Austin and Braam \cite[Lem.~2.5]{AuBr} and Schwarz \cite[Th.~3, p.~69]{Schw} do not actually construct a canonical smooth structure with corners on $\oM(S,T)$; instead, they show that $\cM(S,T)$ is a manifold, and then indicate how to glue the strata of $\oM(S,T)$ together topologically, but without smooth structures. This is enough for the application of constructing and studying the Morse homology groups $H_*^{\rm Mo}(W;\Z)$ of~$W$.

Actually proving the `Folklore Theorem' seems to be surprising difficult. Wehrheim \cite{Wehr} does it when $f$ is Morse and $(f,g)$ have a special normal form (a `Euclidean Morse--Smale pair') near each critical point, and deduces the result with a {\it non-canonical\/} smooth structure on $\oM(S,T)$ when $f$ is Morse. Infinite-dimensional versions of the `Folklore Theorem' are important in Floer theories, Fukaya categories, Symplectic Field Theory, etc., as in~\S\ref{ac63}.

We claim that manifolds with a-corners provide a good framework for setting up and proving the `Folklore Theorem'. In particular, the author expects:
\begin{itemize}
\setlength{\itemsep}{0pt}
\setlength{\parsep}{0pt}
\item[(a)] Make $\lb -\iy,\iy\rb$ into a manifold with a-boundary as in Example \ref{ac3ex2}. Let $f:W\ra\R$ be Morse--Bott, $g$ be any Riemannian metric on $W$, and $\ga:\R\ra W$ be a Morse flow-line as above. Define $\bar\ga:\lb -\iy,\iy\rb\ra W$ by $\bar\ga\vert_\R=\ga$ and $\bar\ga(\pm\iy)=\lim_{t\ra\pm\iy}\ga(t)$. Then $\bar\ga$ is a-smooth, and it satisfies a nonlinear elliptic p.d.e.\ on $\lb -\iy,\iy\rb$ in the sense of~\S\ref{ac53}.

If $f$ is not Morse--Bott, then $\bar\ga$ can be r-smooth but not a-smooth.
\item[(b)] Fix $W,f,g$ as above. Define a {\it family of Morse flow-lines\/} to be a quintuple $(X,Y,\pi,v,\de)$, where $X,Y$ are manifolds with a-corners, $\pi:X\ra Y$ is a proper b-fibration with $X_y:=\pi^{-1}(y)\cong\lb -\iy,\iy\rb$ for $y\in Y^\ci$, and $v\in\Ga^\iy({}^bT(X/Y))$ in the notation of \S\ref{ac54} with $v\vert_x\ne 0$ for all $x\in X$, and $\de:X\ra W$ is a-smooth with ${}^bT\de(v)=\de^*(-\nabla f)$ in~$\Ga^\iy(\de^*(TW))$.

Then for each $y\in Y^\ci$ there should exist a-diffeomorphisms $X_y\cong\lb -\iy,\iy\rb$ identifying $v\vert_{X_y}\cong\frac{\d}{\d t}$, and these a-diffeomorphisms are unique up to $t\mapsto t+c$ in $\lb -\iy,\iy\rb$ for $t\in\R$, as in the equivalence relation $\sim$ in the definition of $\cM(S,T)$ above. Such a-diffeomorphisms identify $\de\vert_{X_y}:X_y\ra W$ with $\bar\ga:\lb -\iy,\iy\rb\ra W$ in (a) for some Morse flow-line~$\ga:\R\ra W$.

For $y\in Y\sm Y^\ci$, $\de\vert_{X_y}:X_y\ra W$ corresponds to a `broken flow-line' as above. Families of Morse flow-lines provide a good way to describe broken flow-lines, and how they occur as limits of Morse flow-lines. Figure \ref{ac6fig1} serves as an illustration of $\pi:X\ra Y$, with a broken flow-line over~$0\in Y$.
\item[(c)] In the situation of the `Folklore Theorem', with $g$ generic, the most natural smooth structure on $\oM(S,T)$ is that of a manifold with a-corners.

The author expects that for $S,T\in C$ there should exist a {\it universal family of Morse flow-lines\/} $(X,Y,\pi,v,\de)$ with a universal property for all families of Morse flow-lines $\bar\ga:\lb -\iy,\iy\rb\ra W$ with $\bar\ga(-\iy)\in S$ and $\bar\ga(\iy)\in T$. This family should be unique up to canonical a-diffeomorphisms of $Y$, and up to canonical a-diffeomorphisms of $X$ plus fibre-dependent $\R$-translations, and in the `Folklore Theorem' we should take~$\oM(S,T)=Y$.
\end{itemize}

\begin{rem} If we define $\oM(S,T)$ in the `Folklore Theorem' as a manifold with a-corners, as in (c), then we can apply $F_\Manacst^\Mancst$ in \S\ref{ac33} to make the moduli space $\oM(S,T)$ into a manifold with ordinary corners. Implicitly this involves a choice of `gluing profile', as in Wehrheim~\cite[Rem.~4.9]{Wehr}.

Now the functor $F_\Manacst^\Mancst$ is defined only for strongly a-smooth maps. But the b-fibrations $\pi:X\ra Y$ in (b),(c) above will generally {\it not\/} be strongly a-smooth. Near a codimension 2 corner of $X$ where flow-lines break, we expect $\pi$ to be locally modelled on $h:\lb 0,\iy)^2\ra\lb 0,\iy)$, $h(x_1,x_2)=x_1x_2$ near $(0,0)$ in $\lb 0,\iy)^2$, and this is a b-fibration, but not strongly a-smooth.

Therefore, although we can make $\oM(S,T)$ into a manifold with corners, {\it we cannot make the universal family $\pi:X\ra\oM(S,T)$ in\/ {\rm(c)} into a smooth map of manifolds with corners}. That is, the author expects that {\it a universal family of Morse flow-lines will generally not exist in\/} $\Manc$. If true, this would be an important reason for using $\Manac$ rather than $\Manc$ in problems of this kind.

\label{ac6rem1}
\end{rem}

\subsection{\texorpdfstring{Symplectic geometry and $J$-holomorphic curves}{Symplectic geometry and J-holomorphic curves}}
\label{ac63}

Let $(S,\om)$ be a symplectic manifold. An {\it almost complex structure\/} $J$ on $S$ is a vector bundle isomorphism $J:TS\ra TS$ with $J^2=-\id$. We call $J$ {\it compatible with\/} $\om$ if defining $g(v,w)=\om(v,Jw)$ for all $v,w\in\Ga^\iy(TS)$ gives a Riemannian metric $g$ on $S$. A $J$-{\it holomorphic curve\/} in $S$ is a Riemann surface $(\Si,j)$ with a smooth map $u:\Si\ra S$ with~$u^*(J)\ci Tu=Tu\ci j:T\Si\ra u^*(TS)$. 

Many important areas of symplectic geometry involve forming moduli spaces $\oM$ of $J$-holomorphic curves in a symplectic manifold $(S,\om)$, `counting' the moduli spaces to get a number, or a homology class (a `virtual class'), or a (co)chain in some (co)homology theory (a `virtual chain'), and using this to define some interesting symplectic invariant which is independent of the choice of almost complex structure $J$. Such areas include Gromov--Witten invariants \cite{FuOn,HWZ1}, Lagrangian Floer cohomology \cite{Fuka,FOOO1}, Fukaya categories \cite{Seid}, and Symplectic Field Theory \cite{EGH}. See McDuff and Salamon \cite{McSa} for an introduction.

For this `counting' of moduli spaces $\oM$ to work, $\oM$ must be compact, and to make $\oM$ compact we must include singular curves with nodes. This involves the kind of `bubbling' and `neck-stretching' issues discussed in~\S\ref{ac62}.

In the general case it is not possible to make the moduli spaces $\oM$ smooth manifolds, even when $J$ is generic. So we need to put a geometric structure on $\oM$, which is strong enough to define virtual classes/virtual chains. There are two main candidates for this structure, both still under development.

Firstly, Fukaya--Oh--Ohta--Ono \cite{FOOO1,FOOO2,FuOn} make moduli spaces $\oM$ into {\it Kuranishi spaces}. The author \cite{Joyc12,Joyc14} found a new definition of Kuranishi space with better categorical properties, and argued that Kuranishi spaces are really {\it derived orbifolds with corners}, in the sense of Derived Algebraic Geometry. Secondly, Hofer--Wysocki--Zehnder \cite{Hofe,HWZ1,HWZ2} give $\oM$ a {\it polyfold Fredholm structure}, where `polyfolds' are, very roughly, a complicated generalization of Banach orbifolds.

The author proposes that manifolds with a-corners should be used in the foundations of both these theories, and believes that this will help solve some current problems in the area, and also make new approaches possible. This was the author's principal motivation for inventing manifolds with a-corners. The rest of this section discusses some of the details of how this might be done.

\subsubsection{Riemann surfaces with a-boundary}
\label{ac631}

Here are two examples of how to use a-boundaries in $J$-holomorphic curves.

\begin{ex} Define $\Si=\lb 0,1\rb\t(\R/2\pi\Z)$, with coordinates $(x,y)$, a compact cylinder with a-boundary. Then $x(1-x)\frac{\pd}{\pd x},\frac{\pd}{\pd y}$ are a basis of sections of ${}^bT\Si$. Define a complex structure (or `b-complex structure') $j:{}^bT\Si\ra {}^bT\Si$ by
\e
\ts j:x(1-x)\frac{\pd}{\pd x}\longmapsto \frac{\pd}{\pd y}\quad\text{and}\quad j:\frac{\pd}{\pd y}\longmapsto -x(1-x)\frac{\pd}{\pd x}.
\label{ac6eq2}
\e
Then $(\Si,j)$ is a Riemann surface with a-boundary.

Now suppose $(S,\om)$ is a symplectic manifold, $J$ is an almost complex structure on $S$ compatible with $\om$, and $u:\Si\ra S$ is a $J$-holomorphic map. That is, $u:\Si\ra S$ is a-smooth with $u^*(J)\ci {}^bTu={}^bTu\ci j:{}^bT\Si\ra u^*({}^bTS)=u^*(TS)$. This is an elliptic equation. As in \S\ref{ac53}, elliptic equations on manifolds with a-boundary $\Si$ restrict to elliptic equations on $\pd\Si$. In this case, on $\pd\Si$ we have
\begin{equation*}
\ts{}^bTu\bigl(\frac{\pd}{\pd y}\bigr)=J\bigl(x(1-x)\frac{\pd}{\pd x}\bigr)=0,
\end{equation*}
as $x(1-x)=0$ on $\pd\Si$, and this means that $u$ is constant on each boundary circle $x=0$ and $x=1$. That is, $u:\Si\ra S$ factors through the quotient $\Si/\!\sim$, where $\sim$ is the equivalence relation on $\Si$ given by $(0,y)\sim(0,y')$ and $(1,y)\sim(1,y')$ for all $y,y'\in \R/2\pi\Z$. But $\Si/\!\sim$ is a 2-sphere $\cS^2$, and there is a homeomorphism $(\Si/\!\sim)\cong\CP^1$ identifying $J$-holomorphic maps $u:\Si\ra S$ and $u':\CP^1\ra S$. We illustrate this in Figure~\ref{ac6fig2}.
\begin{figure}[htb]
\centerline{$\splinetolerance{.8pt}
\begin{xy}
0;<.7mm,0mm>:
,(-20,3);(20,3)**\crv{(-20,5.5)&(-15,8)&(-10,9.25)&(0,10.5)&(10,9.25)&(15,8)&(20,5.5)}
,(-20,3);(20,3)**\crv{(-20,.5)&(-15,-2)&(-10,-3.25)&(0,-4.5)&(10,-3.25)&(15,-2)&(20,.5)}
?(.2)="a"
?(.4)="b"
?(.6)="c"
?(.8)="d"
,(-20,3);(-20,-23)**\crv{(-20,-10)}
,(20,3);(20,-23)**\crv{(20,-10)}
,(-20,-23);(20,-23)**\crv{(-20,-25.5)&(-15,-28)&(-10,-29.25)&(0,-30.5)&(10,-29.25)&(15,-28)&(20,-25.5)}
?(.2)="aa"
?(.4)="bb"
?(.6)="cc"
?(.8)="dd"
,"a";"aa"**@{.}
,"b";"bb"**@{.}
,"c";"cc"**@{.}
,"d";"dd"**@{.}
,(50,-10);(90,-10)**\crv{(50,-2.5)&(55,5)&(60,8.75)&(70,12.5)&(80,8.75)&(85,5)&(90,-2.5)}
?(.2)="e"
?(.4)="f"
?(.6)="g"
?(.8)="h"
,(50,-10);(90,-10)**\crv{(50,-17.5)&(55,-25)&(60,-28.75)&(70,-32.5)&(80,-28.75)&(85,-25)&(90,-17.5)}
?(.05)="ee"
?(.35)="ff"
?(.65)="gg"
?(.95)="hh"
,(70,5)*{\bu}
,(70,5);"e"**\crv{~*{.}(67,6)}
,(70,5);"f"**\crv{~*{.}(67,9)}
,(70,5);"g"**\crv{~*{.}(73,9)}
,(70,5);"h"**\crv{~*{.}(73,6)}
,(70,5);"ee"**\crv{~*{.}(60,4)&(50,-10)}
,(70,5);"ff"**\crv{~*{.}(60,-9)}
,(70,5);"gg"**\crv{~*{.}(80,-9)}
,(70,5);"hh"**\crv{~*{.}(80,4)&(90,-10)}
,(36,-12)*{\longra}
,(36,5)*{\text{collapse}}
,(36,0)*{\text{a-boundary}}
,(36,-5)*{\text{circles}}
,(0,-17)*{\Si}
,(70,-17)*{\Si/\!\sim}
\end{xy}$}
\caption{Riemann surface $\Si$ with a-boundary, and $\Si/\!\sim\,\,\cong\CP^1$}
\label{ac6fig2}
\end{figure}

\label{ac6ex2}
\end{ex}

Next we describe how to model a family of nonsingular $J$-holomorphic curves converging to a singular $J$-holomorphic curve with an interior node, using the material on b-fibrations from \S\ref{ac54}.

\begin{ex} Define a 3-manifold with a-corners $X$ with a proper b-fibration $f:X\ra\lb 0,\ep)=Y$ such that $X_y:=f^{-1}(y)\cong\cS^2$ for $y\in(0,\ep)$, and $\pd X$ (which lies over $y=0$ in $Y$) is the disjoint union of closed 2-discs $D^2_+,D^2_-$, whose common a-boundary $\cS^1$ is a codimension 2 a-corner of $X$. We illustrate this in Figure \ref{ac6fig3}. We can think of $X$ as a family of $\cS^2$'s $X_y$ for $y$ in $Y=\lb 0,\ep)$, which at $y=0$ develop a `fold' along the equator~$\cS^1\subset\cS^2$.

\begin{figure}[htb]
\centerline{$\splinetolerance{.8pt}
\begin{xy}
0;<1.2mm,0mm>:
,(-15,15)*{\bu}
,(-15,-15)*{\bu}
,(-30,3)*{\bu}
,(-18,-3)*{\bu}
,(20,15)*{\ci}
,(20,-15)*{\ci}
,(-10,-23)*{\bu}
,(-15,-23)*{\bu}
,(20,-23)*{\ci}
,(-1,-19)*{f\downarrow}
,(8,-17)*{X}
,(8,-21)*{Y}
,(-17,-23)*{0}
,(22,-23)*{\ep}
,(-14.6,15);(19.5,15)**\crv{(0,15)}
,(-14.6,-15);(19.5,-15)**\crv{(0,-15)}
,(-14.6,-23);(19.5,-23)**\crv{(0,-23)}
,(-15,15);(-30,3)**\crv{(-17,15)&(-25,12)&(-30,5)}
,(-15,-15);(-30,3)**\crv{(-17,-15)&(-25,-10)&(-30,1)}
,(-15,15);(-18,-3)**\crv{(-13,15)&(-13,6)&(-16,-1)}
,(-15,-15);(-17.8,-3.2)**\crv{(-13,-15)&(-13,-8)&(-17,-4)}
,(-30,3);(-18,-3)**\crv{~**\dir{--}(-30,4)&(-24,4)&(-18,-2)}
,(-30,3);(-18,-3)**\crv{(-30,2)&(-24,-4)&(-18,-4)}
,(20,15);(20,-15)**\crv{~**\dir{--}(15,15)&(10,0)&(15,-15)}
,(20,15);(20,-15)**\crv{~**\dir{--}(25,15)&(30,0)&(25,-15)}
,(15,15);(15,-15)**\crv{~*{.} (28,0)}
,(10,15);(10,-15)**\crv{~*{.}(18,8)&(13,-3)&(16,-10)}
,(5,15);(5,-15)**\crv{~*{.}(13,8)&(5,-3)&(11,-10)}
,(0,15);(0,-15)**\crv{~*{.}(8,8)&(-2,-3)&(7,-10)}
,(-5,15);(-5,-15)**\crv{~*{.}(3,8)&(-13,-3)&(2,-10)}
,(-10,15);(-10,-15)**\crv{~*{.}(-3,8)&(-21,-4)&(-5,-10)}
,(-37,10)*{\text{a-boundary $D^2_+$}}
,(-34,-10)*{\text{a-boundary $D^2_-$}}
,(-38.3,3.3)*{\text{a-corner $\cS^1$}}
,(36,0)*{\text{fibres $X_y$}}
,(37,-3)*{\text{for $y>0$ are}}
,(37,-6)*{\text{2-spheres $\cS^2$}}
,(-5,-23)*{\bu}
,(0,-23)*{\bu}
,(5,-23)*{\bu}
,(10,-23)*{\bu}
,(15,-23)*{\bu}
,(-10,-23)*{\bu}
\end{xy}$}
\caption{B-fibration $f:X\ra Y$ for family of $\CP^1$'s developing interior node}
\label{ac6fig3}
\end{figure}
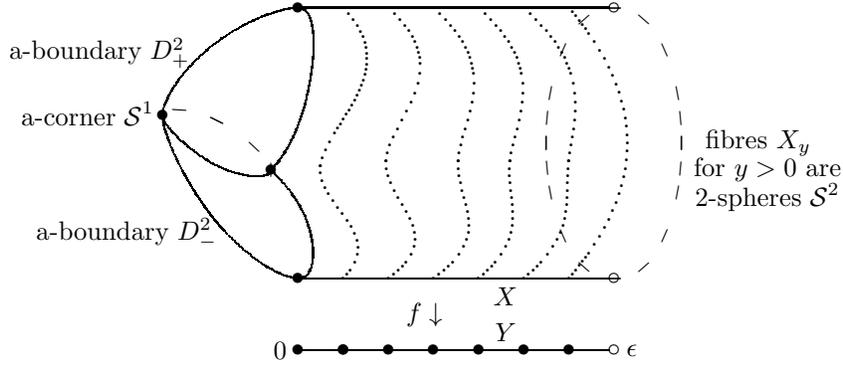

The relative tangent bundle ${}^bT(X/Y)$ of $f:X\ra Y$ from \S\ref{ac54} has rank 2, and we can choose an a-smooth morphism $j:{}^bT(X/Y)\ra {}^bT(X/Y)$ with $j^2=-\id$. Then $j_y=j\vert_{X_y}$ makes $(X_y,j_y)$ into a Riemann surface isomorphic to $\CP^1$ for all $y\in(0,\ep)$. Also $(X_0,j_0)$ is the union of $(D^2_+,j_{0,+})$ and $(D^2_-,j_{0,-})$ glued along $\cS^1$, where $(D^2_\pm,j_{0,\pm})$ are Riemann surfaces with a-boundary, which look near their ends like one end of the cylinder $(\Si,j)$ in Example~\ref{ac6ex2}.

Now let $(S,\om)$ be a symplectic manifold, $J$ an almost complex structure on $S$ compatible with $\om$, and $u:X\ra S$ be a family of $J$-holomorphic maps over $Y$. That is, $u:X\ra S$ is a-smooth with
\begin{equation*}
J\ci{}^bTu\vert_{{}^bT(X/Y)}={}^bTu\ci j:{}^bT(X/Y)\longra u^*({}^bTS)=u^*(TS).
\end{equation*}
Then $u\vert_{X_y}:X_y\ra S$ for $y\in(0,\ep)$ and $u\vert_{D^2_\pm}:D^2_\pm\ra S$ are $J$-holomorphic maps.

\begin{figure}[htb]
\centerline{$\splinetolerance{.8pt}
\begin{xy}
0;<1.2mm,0mm>:
,(-15,15)*{\bu}
,(-15,-15)*{\bu}
,(-25,0)*{\bu}
,(20,15)*{\ci}
,(20,-15)*{\ci}
,(-10,-23)*{\bu}
,(-15,-23)*{\bu}
,(20,-23)*{\ci}
,(-1,-19)*{\ti f\downarrow}
,(8,-17)*{\ti X}
,(8,-21)*{Y}
,(-17,-23)*{0}
,(22,-23)*{\ep}
,(-14.6,15);(19.5,15)**\crv{(0,15)}
,(-14.6,-15);(19.5,-15)**\crv{(0,-15)}
,(-14.6,-23);(19.5,-23)**\crv{(0,-23)}
,(-15,15);(-25,0)**\crv{(-17,15)&(-25,12)&(-29,1)}
,(-15,-15);(-25,0)**\crv{(-17,-15)&(-25,-10)&(-29,-1)}
,(-15,15);(-25,0)**\crv{(-13,15)&(-13,6)&(-20,2)}
,(-15,-15);(-25,0)**\crv{(-13,-15)&(-13,-8)&(-20,-2)}
,(20,15);(20,-15)**\crv{~**\dir{--}(15,15)&(10,0)&(15,-15)}
,(20,15);(20,-15)**\crv{~**\dir{--}(25,15)&(30,0)&(25,-15)}
,(15,15);(15,-15)**\crv{~*{.} (28,0)}
,(10,15);(10,-15)**\crv{~*{.}(16,8)&(13,0)&(16,-8)}
,(5,15);(5,-15)**\crv{~*{.}(11,8)&(3,0)&(11,-8)}
,(0,15);(0,-15)**\crv{~*{.}(7,8)&(-7,0)&(7,-8)}
,(-5,15);(-5,-15)**\crv{~*{.}(3,8)&(-18,0)&(2,-8)}
,(-10,15);(-10,-15)**\crv{~*{.}(-5,8)&(-28,0)&(-5,-8)}
,(-33,10)*{\text{$D^2_+/\!\sim\,\cong\cS^2$}}
,(-33,-10)*{\text{$D^2_-/\!\sim\,\cong\cS^2$}}
,(-37.3,.3)*{\text{$\cS^1/\!\sim$, a point $\ra$}}
,(36,0)*{\text{fibres $\ti X_y$}}
,(37,-3)*{\text{for $y>0$ are}}
,(37,-6)*{\text{2-spheres $\cS^2$}}
,(-5,-23)*{\bu}
,(0,-23)*{\bu}
,(5,-23)*{\bu}
,(10,-23)*{\bu}
,(15,-23)*{\bu}
,(-10,-23)*{\bu}
\end{xy}$}
\caption{B-fibration $f:X\ra Y$ with 2-corner $\cS^1\subset X$ collapsed to a point}
\label{ac6fig5}
\end{figure}
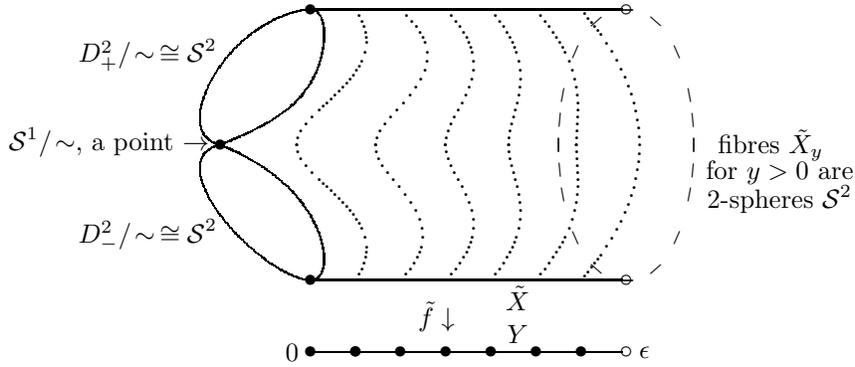

The argument of Example \ref{ac6ex2} applied to $u\vert_{D^2_\pm}:D^2_\pm\ra S$ shows that $u\vert_{\cS^1}$ is constant, say $u\vert_{\cS^1}=p\in S$. Thus $u$ factors through the quotient $\ti X=X/\sim$, where $\sim$ is the equivalence relation on $X$ which collapses $\cS^1$ to a point. Note that $\ti X$ is not a manifold with (a-)corners. Also $f:X\ra Y$ factors through $\ti f:\ti X\ra Y$, where $\ti f^{-1}(y)\cong\CP^1$ for $y\in(0,\ep)$, and $\ti f^{-1}(0)$ is two $\CP^1$'s joined at a node, a Deligne--Mumford prestable $\CP^1$. We illustrate this in Figure~\ref{ac6fig5}.

\label{ac6ex3}
\end{ex}

\begin{rem} Let $\ti X_0$ be a Deligne--Mumford prestable curve, in the sense of \cite{FuOn,HWZ1}, with a single interior node at $\ti x_0\in\ti X_0$. That is, $\ti X_0$ near $\ti x_0$ is locally modelled on $\bigl\{(s,t)\in\C^2:st=0\bigr\}$ near $(0,0)$. So $\ti X_0$ is not a manifold.

To model $\ti X_0$ using manifolds with a-corners, roughly speaking one should do a `real blow-up' at $\ti x_0$. This happens in two slightly different ways:
\begin{itemize}
\setlength{\itemsep}{0pt}
\setlength{\parsep}{0pt}
\item[(i)] If we are thinking of $\ti X_0$ just as a single curve, not part of a family, we define a proper morphism $\pi:\hat X_0\ra\ti X_0$ with $\hat X_0$ a manifold with a-boundary, where $\pi^{-1}(\ti x_0)=\cS^1_s\amalg\cS^1_t$ is two disjoint copies of $\cS^1$, the a-boundaries of the two local components of $\ti X_0$ meeting at $\ti x_0$, and $\pi$ is a diffeomorphism over $\ti X_0\sm\{\ti x_0\}$. In Example \ref{ac6ex3}, $\hat X_0$ corresponds to $D^2_+\amalg D^2_-$. Note that the two copies $\cS^1_s,\cS^1_t$ of $\cS^1$ in $\pi^{-1}(\ti x_0)$ are not identified.
\item[(ii)] If $\ti X_0$ is part of a family of curves including nonsingular curves, say $\ti f:\ti X\ra\lb 0,\ep)$ with $\ti X_y=\ti f^{-1}(y)$ locally modelled on $\bigl\{(s,t)\in\C^2:st=y\bigr\}$ near $\ti x_0$, then as in Example \ref{ac6ex3} we model the family as a proper b-fibration $f:X\ra\lb 0,\ep)$ in $\Manac$ with a real blow-up map $\pi:X\ra\ti X$ with $\pi^{-1}(x_0)\cong\cS^1$, and $\pi$ is a diffeomorphism over $\ti X\sm\{\ti x_0\}$. Note that $X_0=f^{-1}(0)$ is not a manifold with a-corners, though $\hat X_0=\pd X$ is. It is obtained by identifying the $\cS^1_s\cong\cS^1_t$ in $\hat X_0$ in (i), and in Example \ref{ac6ex3} corresponds to~$D^2_+\amalg_{\cS^1}D^2_-$.
\end{itemize}

In (i), the complex structure $j$ on $\hat X_0$ determines a nonvanishing vector field $v$ on $\pd\hat X_0=\cS^1_s\amalg\cS^1_t$, up to a positive constant, which we fix by requiring $v$ to have period $2\pi$ on each $\cS^1$, and regard as the vector field of a $\U(1)$-action on $\cS^1\amalg\cS^1$. The identification $\cS^1_s\cong\cS^1_t$ in (ii) is $\U(1)$-equivariant, but this does not determine it uniquely, the family of possible $\U(1)$-equivariant isomorphisms $\cS^1_s\cong\cS^1_t$ is a torsor for $\U(1)$. The actual identification $\cS^1_s\cong\cS^1_t$ used to build $X,X_0$ in (ii) depends not just on $\ti X_0$, but also on $\ti X_y$ as~$y\ra 0$.

This has important consequences for how moduli spaces of $J$-holomorphic curves should be described using a-corners. The pictures of \cite{FuOn,HWZ1} yield a moduli space $\ti{\cal M}$ with a real (virtual) codimension 2 subset $\ti{\cal M}_0\subseteq\ti{\cal M}$ corresponding to curves with nodes, where $\ti{\cal M}_0$ is not considered part of the boundary of~$\ti{\cal M}$. 

The a-corners picture naturally yields a moduli space $\cM$ with projection $\pi:\cM\ra\ti{\cal M}$, with $\pi^{-1}(p)\cong(\cS^1)^k$ when $p$ represents a curve with $k$ nodes. Points of $\cM$ correspond to a stable $J$-holomorphic curve $[X,j,u]$ in $\ti{\cal M}$ together with identifications $\smash{\cS^1_s\cong\cS^1_t}$ at each interior node of $X$. The boundary $\pd\cM$ has an extra  component $\pd_{\rm node}\cM$ lying over $\ti{\cal M}_0$, but there is a $\U(1)$-action on $\pd_{\rm node}\cM$ coming from changing the identification $\cS^1_s\cong\cS^1_t$, and using this $\U(1)$-action we can arrange that the virtual chain of $\pd_{\rm node}\cM$ is trivial, so this additional boundary will not cause problems in the virtual cycle theory.
\label{ac6rem2}
\end{rem}

\subsubsection{Including Lagrangian boundary conditions}
\label{ac632}

The theories of Lagrangian Floer cohomology \cite{Fuka,FOOO1} and Fukaya categories \cite{Seid} are built on studying moduli spaces of $J$-holomorphic curves $u:\Si\ra(S,\om)$ from a Riemann surface $\Si$ with corners whose boundary $\pd\Si$ is mapped by $u$ to one or more Lagrangians $L_1,\ldots,L_k$ in the symplectic manifold $S$. 

To study this in the `a-corners' picture, we must take $\Si$ to be a 2-manifold with both corners and a-corners as in \S\ref{ac35}, where the ordinary boundary $\pd^{\rm c}\Si$ is mapped to Lagrangians in $S$, and the a-boundary $\pd^{\rm ac}\Si$ does not satisfy explicit boundary conditions, but the $J$-holomorphic curve equation forces $u$ to be locally constant on $\pd^{\rm ac}\Si$. Here is an example.

\begin{ex} Define $\Si=\lb 0,1\rb\t[0,1]$, with coordinates $(x,y)$, a compact square with boundary and a-boundary. Then $x(1-x)\frac{\pd}{\pd x},\frac{\pd}{\pd y}$ are a basis of sections of ${}^mT\Si$ in Remark \ref{ac4rem3}(b). Define a complex structure $j:{}^mT\Si\ra {}^mT\Si$ by \eq{ac6eq2}. Then $(\Si,j)$ is a Riemann surface with boundary and a-boundary.

Now suppose $(S,\om)$ is a symplectic manifold, $L_0,L_1$ are Lagrangians in $S$, $J$ is an almost complex structure on $S$ compatible with $\om$, and $u:\Si\ra S$ is a $J$-holomorphic map with $u(x,0)\in L_0$ and $u(x,1)\in L_1$ for all~$x\in\lb 0,1\rb$.

Notice that on the boundaries $\pd^{\rm c}\Si$ at $y=0$ and $y=1$ we impose boundary conditions, that $\Si$ maps to $L_0$ and $L_1$. But on the a-boundaries $\pd^{\rm ac}\Si$ at $x=0$ and $x=1$ we impose no boundary conditions. This is necessary to get a well-behaved moduli problem, with Fredholm linearization.

As in Example \ref{ac6ex2}, the $J$-holomorphic curve equation forces ${}^mTu\bigl(\frac{\pd}{\pd y}\bigr)=0$ on the a-boundaries $x=0$ and $x=1$, so that they are mapped to points. As $u(0,0)\in L_0$ and $u(0,1)\in L_1$ we must have $u(0,y)=p\in L_0\cap L_1$ for all $y\in[0,1]$, and similarly $u(1,y)=q\in L_0\cap L_1$ for all $y\in[0,1]$.
Thus $u:\Si\ra S$ factors through the quotient $\Si/\!\sim$, where $\sim$ is the equivalence relation on $\Si$ given by $(0,y)\sim(0,y')$ and $(1,y)\sim(1,y')$ for all $y,y'\in[0,1]$. We illustrate this in Figure~\ref{ac6fig2}.

\begin{figure}[htb]
\centerline{$\splinetolerance{.8pt}
\begin{xy}
0;<1.2mm,0mm>:
,(-65,5)*{\bu}
,(-65,-5)*{\bu}
,(-35,5)*{\bu}
,(-35,-5)*{\bu}
,(-65,5);(-35,5)**\crv{(-50,5)}
?(.9)="aa"
?(.8)="bb"
?(.7)="cc"
?(.6)="dd"
?(.5)="ee"
?(.4)="ff"
?(.3)="gg"
?(.2)="hh"
?(.1)="ii"
,(-65,-5);(-35,-5)**\crv{(-50,-5)}
?(.9)="aaa"
?(.8)="bbb"
?(.7)="ccc"
?(.6)="ddd"
?(.5)="eee"
?(.4)="fff"
?(.3)="ggg"
?(.2)="hhh"
?(.1)="iii"
,"aa";"aaa"**@{.}
,"bb";"bbb"**@{.}
,"cc";"ccc"**@{.}
,"dd";"ddd"**@{.}
,"ee";"eee"**@{.}
,"ff";"fff"**@{.}
,"gg";"ggg"**@{.}
,"hh";"hhh"**@{.}
,"ii";"iii"**@{.}
,(-65,5);(-65,-5)**\crv{(-65,0)}
,(-35,5);(-35,-5)**\crv{(-35,0)}
,(-10,0);(10,0)**\crv{(0,10)}
?(.95)="a"
?(.85)="b"
?(.75)="c"
?(.65)="d"
?(.55)="e"
?(.45)="f"
?(.35)="g"
?(.25)="h"
?(.15)="i"
?(.05)="j"
?(.5)="y"
,(-10,0);(-15,-6)**\crv{(-15,-5)}
,(10,0);(15,-6)**\crv{(15,-5)}
,(-10,0);(10,0)**\crv{(0,-10)}
?(.95)="k"
?(.85)="l"
?(.75)="m"
?(.65)="n"
?(.55)="o"
?(.45)="p"
?(.35)="q"
?(.25)="r"
?(.15)="s"
?(.05)="t"
?(.5)="z"
,(-10,0);(-15,6)**\crv{(-15,5)}
,(10,0);(15,6)**\crv{(15,5)}
,"a";"k"**@{.}
,"b";"l"**@{.}
,"c";"m"**@{.}
,"d";"n"**@{.}
,"e";"o"**@{.}
,"f";"p"**@{.}
,"g";"q"**@{.}
,"h";"r"**@{.}
,"i";"s"**@{.}
,"j";"t"**@{.}
,"y"*{<}
,"z"*{>}
,(-10,0)*{\bu}
,(-10,-3)*{p}
,(10,0)*{\bu}
,(10,-3)*{q}
,(-50,0)*{\Si}
,(0,0)*{\Si/\!\sim}
,(-12,6)*{L_0}
,(-12,-6)*{L_1}
,(17,4)*{L_0}
,(17.5,-4)*{L_1}
,(-50,-7)*{\text{boundary in $L_0$}}
,(-50,7)*{\text{boundary in $L_1$}}
,(-73,0)*{\text{a-boundary}}
,(-24,-3.5)*{\longra}
,(-24,3.5)*{\text{collapse}}
,(-24,0)*{\text{a-boundary}}
\end{xy}$}
\caption{Holomorphic disc $\Si$ with boundary in $L_0\cup L_1$}
\label{ac6fig4}
\end{figure}
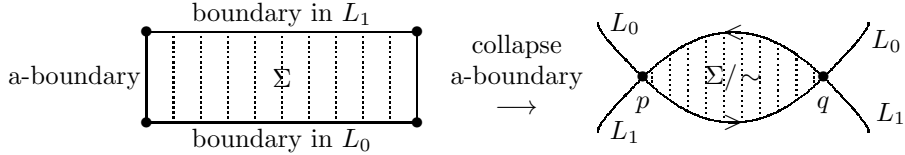

Now $J$-holomorphic curves of this type are used to define the Lagrangian Floer cohomology groups $HF^*(L_0,L_1)$ in \cite{Fuka,FOOO1,Seid}. We regard them as squares in $\Mancac$, as in the left hand picture of Figure \ref{ac6fig4}, rather than the conventional notion of 2-gons in $\Manc$, 
as in the right hand picture of Figure~\ref{ac6fig4}.

\label{ac6ex4}
\end{ex}

\begin{rem} For the deformation theory of curves in Example \ref{ac6ex4} to be well behaved, one usually requires $L_0,L_1$ to intersect transversely, or at least cleanly, at $p,q$. In our picture, this condition is {\it not local on\/} $\Si$, as the $L_0,L_1$ boundary conditions happen on disjoint parts of the curve. So we might guess that transverse intersection of $L_0,L_1$ is not necessary for the a-corners approach.

In fact it is needed. As in \S\ref{ac53}, elliptic operators $P^p_{k,\la}:L^p_{k+l}(E)_\la\ra L^p_k(F)_\la$ on compact manifolds with a-boundary $X$ are non-Fredholm at a subset of weights $\la$ determined by $P\vert_{\pd X}$. If $L_0,L_1$ do not intersect transversely at $p,q$, then the linearized $J$-holomorphic curve equation is not Fredholm when~$\la=0$.

\label{ac6rem3}
\end{rem}

We can also discuss families of $J$-holomorphic curves with Lagrangian boundary conditions using b-fibrations in $\Mancac$, as in Example \ref{ac6ex3} and Remark~\ref{ac6rem2}.

\subsubsection{Manifolds with a-corners, Kuranishi spaces, and polyfolds}
\label{ac633}

Setting up the foundations of $J$-holomorphic curve theory, moduli spaces, and virtual chains, in Symplectic Geometry in the general case is a mammoth task. Two rival groups who have been working on this for a long time are Fukaya, Oh, Ohta and Ono \cite{Fuka,FOOO1,FOOO2,FuOn}, who make moduli spaces into {\it Kuranishi spaces}, and Hofer, Wysocki and Zehnder \cite{Hofe,HWZ1,HWZ2} give moduli spaces a {\it polyfold Fredholm structure}. The two are related by Yang \cite{Yang}, who defines a `truncation functor' from polyfold Fredholm structures to Kuranishi spaces.

We want to point out that manifolds with a-corners can be incorporated very easily into both the Kuranishi space and the polyfold theories, and that there would be definite advantages to doing this. The definition of Kuranishi spaces in \cite[\S A.1]{FOOO2} involves {\it Kuranishi neighbourhoods\/} $(V,E,\Ga,\psi,s)$ in which $V$ is a manifold with corners, but we can use manifolds with a-corners instead.

The author \cite{Joyc12,Joyc14} produced a new definition of Kuranishi spaces, refining that of \cite[\S A.1]{FOOO1}, which form a well-behaved 2-category (Fukaya et al.\ \cite{Fuka,FOOO1,FOOO2,FuOn} do not define morphisms between Kuranishi spaces). Starting with any category of `manifolds' $\bf{\ti Man}$ satisfying certain conditions, the construction of \cite{Joyc12} produces an associated 2-category of `Kuranishi spaces' $\bf{\ti Kur}$. Our categories $\Manac$ and $\Mancac$ satisfy these conditions, so we immediately get 2-categories $\bf Kur^{ac}$ or $\bf Kur^{c,ac}$ of Kuranishi spaces with a-corners, or with corners and a-corners.

The definition of polyfolds involves `partial quadrants' \cite[Def.~1.11]{HWZ2}, of the form $[0,\iy)^k\t W$ for $W$ an sc-Banach space, which are an infinite-dimensional generalization of the local models $\R^m_k=[0,\iy)^k\t\R^{m-k}$ for manifolds with a-corners. By replacing these by $\lb 0,\iy)^k\t W$ and using a-smoothness, we can define `polyfolds with a-corners'. 

Here are some reasons for doing this. In both theories, when proving that moduli spaces including singular curves have a Kuranishi/polyfold Fredholm structure, the proofs work by first proving estimates of the type \eq{ac3eq3} (our definition of a-smooth function) on data such as the Kuranishi section $s:V\ra E$ in $(V,E,\Ga,\psi,s)$, in the natural coordinates. Then, using a `gluing profile', they change to different coordinates in which $s$ is actually smooth.

In our language, these proofs work by firstly constructing Kuranishi neighbourhoods $(V,E,\Ga,\psi,s)$ with $V$ a manifold with a-corners, $s:V\ra E$ a-smooth, etc., and secondly applying the functor $F_\Manacst^\Mancst:\Manacst\ra\Mancst$ from \S\ref{ac33} to get to manifolds with corners. For examples of the first step see 
\cite[\S A1.4, Lem.~A1.58]{FOOO1}, \cite[Th.~6.4]{FOOO2}, \cite[\S 4.4]{HWZ1}, and for the second step see \cite[\S A1.4, p.~777]{FOOO1}, \cite[\S 8]{FOOO2}, and \cite[\S 2.1 \& \S 2.6]{HWZ1}. If we used manifolds with a-corners from the outset, we would get a-smoothness in the natural coordinates, and the second step would be unnecessary.

As in \S\ref{ac33}, the functor $F_\Manacst^\Mancst$ works only for strongly a-smooth maps. Therefore we should expect that any operation in the theory of \cite{Fuka,FOOO1,FOOO2,FuOn,Hofe,HWZ1,HWZ2} that involves non-strongly a-smooth maps will be smooth in the a-corners picture, but will become non-smooth after applying gluing profiles.

For example, the Lagrangian Floer cohomology of \cite{FOOO1} involves moduli spaces $\oM_k$ of stable $J$-holomorphic discs $u:\Si\ra S$ with boundary in a Lagrangian $L\subset S$, with $k$ boundary marked points $z_1,\ldots,z_k\in\pd\Si$. We can consider the map $F_i:\oM_k\ra\oM_{k-1}$ which forgets the $i^{\rm th}$ marked point $z_i$ for $i=1,\ldots,k$. The local models for such $F_i$ on Kuranishi neighbourhoods can involve non-strongly a-smooth maps. So we expect that $F_i:\oM_k\ra\oM_{k-1}$ is an a-smooth 1-morphism of Kuranishi spaces with a-corners, in the sense of \cite{Joyc12}, but may not be smooth in the set-up of~\cite{Fuka,FOOO1,FOOO2,FuOn,Hofe,HWZ1,HWZ2}.

\subsubsection{Universal families and representable 2-functors}
\label{ac634}

The usual methods for constructing moduli spaces in Algebraic Geometry and Differential Geometry are very different. To construct moduli $\K$-schemes in Algebraic Geometry, one uses Grothendieck's method of representable functors. Writing $\Sch_\K$ for the category of schemes over a field $\K$, one defines a functor
\begin{equation*}
F:(\Sch_\K)^{\bf op}\longra\Sets,
\end{equation*}
where $F(S)$ is the set of isomorphism classes of families of objects in the moduli problem over a base $\K$-scheme $S$. For example, to study moduli of vector bundles $E\ra Y$ over a fixed smooth projective $\K$-scheme $Y$, we could take $F(S)$ to be the set of isomorphism classes of vector bundles~$E_S\ra Y\t S$.

Then we hope to prove that there exists a $\K$-scheme $\cM$ (necessarily unique up to isomorphism) with a natural isomorphism $F\cong \Hom(-,\cM)$. We call $\cM$ the {\it moduli $\K$-scheme}. There is a canonical family of objects $U\in F(\cM)$ over $\cM$ called the {\it universal family}, corresponding to~$\id_M\in\Hom(\cM,\cM)$.

We propose using this method in Differential Geometry. For problems in which the moduli space should be a smooth manifold with a-corners, as for the Morse flow-lines in \S\ref{ac622}, we can consider functors
\begin{equation*}
F:(\Manac)^{\bf op}\longra\Sets.
\end{equation*}
For problems where the moduli space should be a Kuranishi space with a-corners, we need a 2-categorical version, so we should consider weak 2-functors
\begin{equation*}
F:({\bf mKur^{ac}})^{\bf op}\longra{\bf Groupoids},
\end{equation*}
where $\bf mKur^{ac}$ is the 2-category of Kuranishi spaces with a-corners with trivial isotropy groups (i.e.\ the manifold version of Kuranishi spaces with a-corners).

For this method to work, it will be {\it essential\/} to use manifolds or Kuranishi spaces with a-corners, not with corners. This is because if we use ordinary corners then a universal family may not exist, as in Remark \ref{ac6rem1}, so the 2-functor $F$ will not be representable.

\medskip

\noindent{\small\sc The Mathematical Institute, Radcliffe Observatory Quarter, Woodstock Road, Oxford, OX2 6GG, U.K.}

\noindent{\small\sc E-mail: \tt joyce@maths.ox.ac.uk}


\begin{thebibliography}{99}
\addcontentsline{toc}{section}{References}


\bibitem{AuBr} D.M. Austin and P.M. Braam, {\it Morse--Bott theory and equivariant cohomology}, pages 123--183 in {\it The Floer memorial volume}, Progress in Math. 133, Birkh\"auser, Basel, 1995.

\bibitem{Bred1} G.E. Bredon, {\it Topology and Geometry}, Graduate
Texts in Math. 139, Springer, New York, 1993.

\bibitem{Bred2} G.E. Bredon, {\it Sheaf Theory}, second edition, Graduate Texts in Math. 170, Springer, New York, 1997.

\bibitem{BrSa} R.L. Bryant and S.M. Salamon, {\it On the construction of some complete metrics with exceptional holonomy}, Duke Math. J. 58 (1989), 829--850.

\bibitem{Carr1} G. Carron, {\it Cohomologie $L^2$ des vari\'et\'es QALE}, J. Reine Angew. Math. 655 (2011), 1--59. \href{http://arxiv.org/abs/math/0501290}{math.DG/0501290}. 

\bibitem{Carr2} G. Carron, {\it On the Quasi-Asymptotically Locally Euclidean geometry of Nakajima's metric}, J. Inst. Math. Jussieu 10 (2011), 119--147. \hfil\break \href{http://arxiv.org/abs/0811.3870}{arXiv:0811.3870}.

\bibitem{Cerf} J. Cerf, {\it Topologie de certains espaces de
plongements}, Bull. Soc. Math. France 89 (1961), 227--380.

\bibitem{Chan1} Y.-M. Chan, {\it Desingularizations of Calabi--Yau $3$-folds with a conical singularity}, Q. J. Math. 57 (2006), 151--181. \href{http://arxiv.org/abs/math/0410260}{math.DG/0410260}.

\bibitem{Chan2} Y.-M. Chan, {\it Simultaneous desingularizations of Calabi--Yau and special Lagrangian $3$-folds with conical singularities, I, II}, Ann. Global Anal. Geom. 35 (2009), 91--114 and 157--180. \href{http://arxiv.org/abs/math/0606399}{math.DG/0606399}.

\bibitem{CoHe1} R.J. Conlon and H.-J. Hein, {\it Asymptotically conical Calabi--Yau manifolds, I}, Duke Math. J. 162 (2013), 2855--2902. \href{http://arxiv.org/abs/1205.6347}{arXiv:1205.6347}.

\bibitem{CoHe2} R.J. Conlon and H.-J. Hein, {\it Asymptotically conical Calabi--Yau metrics on quasi-projective varieties}, Geom. Funct. Anal. 25 (2015), 517--552. \href{http://arxiv.org/abs/1301.5312}{arXiv:1301.5312}.

\bibitem{CMR} R.J. Conlon, R. Mazzeo and F. Rochon, {\it The moduli space of asymptotically cylindrical Calabi--Yau manifolds}, Comm. Math. Phys. 338 (2015), 953--1009. \href{http://arxiv.org/abs/1408.6562}{arXiv:1408.6562}.

\bibitem{CHNP1} A. Corti, M. Haskins, J. Nordstr\"om, and T. Pacini, {\it $G_2$-manifolds and associative submanifolds via semi-Fano $3$-folds}, Duke Math. J. 164 (2015), 1971--2092. \href{http://arxiv.org/abs/1207.4470}{arXiv:1207.4470}.

\bibitem{CHNP2} A. Corti, M. Haskins, J. Nordstr\"om, and T. Pacini, {\it Asymptotically cylindrical Calabi--Yau $3$-folds from weak Fano $3$-folds}, Geom. Topol. 17 (2013), 1955--2059. \href{http://arxiv.org/abs/1206.2277}{arXiv:1206.2277}.

\bibitem{DeMa} A. Degeratu and R. Mazzeo, {\it Fredholm theory for elliptic operators on quasi-asymptotically conical spaces}, \href{http://arxiv.org/abs/1406.3465}{arXiv:1406.3465}, 2014.

\bibitem{DoKr} S.K. Donaldson and P.B. Kronheimer, {\it The Geometry of Four-Manifolds}, OUP, Oxford, 1990.

\bibitem{EGH} Y. Eliashberg, A. Givental and H. Hofer, {\it
Introduction to Symplectic Field Theory}, Geom. Funct. Anal. 2000,
Special Volume, Part II, 560--673. \href{http://arxiv.org/abs/math/0010059}{math.SG/0010059}.

\bibitem{Fuka} K. Fukaya, {\it Floer homology of Lagrangian
submanifolds}, Sugaku Expositions 26 (2013), 99--127. \href{http://arxiv.org/abs/1106.4882}{arXiv:1106.4882}.

\bibitem{FOOO1} K. Fukaya, Y.-G. Oh, H. Ohta and K. Ono,
{\it Lagrangian intersection Floer theory --- anomaly and
obstruction}, Parts I \& II. AMS/IP Studies in Advanced Mathematics,
46.1 \& 46.2, A.M.S./International Press, 2009.

\bibitem{FOOO2} K. Fukaya, Y.-G. Oh, H. Ohta and K. Ono, {\it Exponential decay estimates and smoothness of the moduli space of pseudoholomorphic curves}, \href{http://arxiv.org/abs/1603.07026}{arXiv:1603.07026}, 2016.

\bibitem{FuOn} K. Fukaya and K. Ono, {\it Arnold Conjecture and
Gromov--Witten invariant}, Topology 38 (1999), 933--1048.

\bibitem{Grie} D. Grieser, {\it Basics of the b-calculus}, pages 30--84 in J.B. Gil, D. Grieser and M. Lesch, {\it Approaches to Singular Analysis}, Operator Theory 125, Birkh\"auser, Basel, 2001. \href{http://arxiv.org/abs/math/0010314}{math.AP/0010314}.

\bibitem{HaLa} R. Harvey and H.B. Lawson, {\it Calibrated geometries},
Acta Mathematica 148 (1982), 47--157.

\bibitem{HHN} M. Haskins, H.-J. Hein and J. Nordstr\"om, {\it Asymptotically cylindrical Calabi--Yau manifolds}, J. Diff. Geom. 101 (2015), 213--265. \hfil\break \href{http://arxiv.org/abs/1212.6929}{arXiv:1212.6929}.

\bibitem{HMM} A. Hassell, R. Mazzeo and R.B. Melrose, {\it Analytic surgery and the accumulation of eigenvalues}, Comm. Anal. Geom. 3 (1995), 115--222.

\bibitem{Hofe} H. Hofer, {\it Polyfolds and Fredholm Theory}, \href{http://arxiv.org/abs/1412.4255}{arXiv:1412.4255}, 2014.

\bibitem{HWZ1} H. Hofer, K. Wysocki and E. Zehnder, {\it
Applications of polyfold theory I: the polyfolds of Gromov--Witten
theory}, \href{http://arxiv.org/abs/1107.2097}{arXiv:1107.2097}, 2011.

\bibitem{HWZ2} H. Hofer, K. Wysocki and E. Zehnder, {\it Polyfold and Fredholm theory I: basic theory in M-polyfolds}, \href{http://arxiv.org/abs/1407.3185}{arXiv:1407.3185}, 2014.

\bibitem{Joyc1} D. Joyce, {\it A new construction of compact\/ $8$-manifolds with holonomy $\Spin(7)$}, J. Diff. Geom. 53 (1999), 89--130. \href{http://arxiv.org/abs/math/9910002}{math.DG/9910002}.

\bibitem{Joyc2} D. Joyce, {\it Compact Manifolds with Special Holonomy},
OUP, Oxford, 2000.

\bibitem{Joyc3} D. Joyce, {\it Asymptotically Locally Euclidean metrics with holonomy $SU(m)$}, Ann. Global
Anal. Geom. 19 (2001), 55--73. \href{http://arxiv.org/abs/math/9905041}{math.AG/9905041}.

\bibitem{Joyc4} D. Joyce, {\it Quasi-ALE metrics with holonomy $SU(m)$ and $Sp(m)$}, Ann. Global Anal. Geom. 19 (2001), 103--132. \href{http://arxiv.org/abs/math/9905043}{math.AG/9905043}.

\bibitem{Joyc5} D. Joyce, {\it Constructing special Lagrangian
$m$-folds in $\C^m$ by evolving quadrics}, Math. Ann. 320 (2001),
757--797. \href{http://arxiv.org/abs/math/0008155}{math.DG/0008155}.

\bibitem{Joyc6} D. Joyce, {\it Special Lagrangian $m$-folds in $\C^m$ 
with symmetries}, Duke Math. J. 115 (2002), 1--51. \href{http://arxiv.org/abs/math/0008021}{math.DG/0008021}.

\bibitem{Joyc7} D. Joyce, {\it Special Lagrangian submanifolds
with isolated conical singularities. I--IV}, Ann. Global
Anal. Geom. 25 (2004), 201--251; 25 (2004), 301--352; 26 (2004), 1--58; and 26 (2004), 117--174. \href{http://arxiv.org/abs/math/0211294}{math.DG/0211294}, \href{http://arxiv.org/abs/math/0211295}{math.DG/0211295}, \href{http://arxiv.org/abs/math/0302355}{math.DG/0302355} and \href{http://arxiv.org/abs/math/0302356}{math.DG/0302356}.

\bibitem{Joyc8} D. Joyce, {\it Special Lagrangian submanifolds with isolated conical singularities. V. Survey and applications}, J. Diff. Geom. 63 (2003), 279--347. \href{http://arxiv.org/abs/math/0303272}{math.DG/0303272}.

\bibitem{Joyc9} D. Joyce, {\it Riemannian holonomy groups and calibrated
geometry}, Oxford University Press, 2007.

\bibitem{Joyc10} D. Joyce, {\it On manifolds with corners}, pages
225--258 in S. Janeczko et al., editors, {\it
Advances in Geometric Analysis}, Advanced Lectures in Mathematics
21, International Press, Boston, 2012. \href{http://arxiv.org/abs/0910.3518}{arXiv:0910.3518}.

\bibitem{Joyc11} D. Joyce, {\it Algebraic Geometry over\/ $C^\iy$-rings}, \href{http://arxiv.org/abs/1001.0023}{arXiv:1001.0023}, 2010.

\bibitem{Joyc12} D. Joyce, {\it A new definition of Kuranishi space}, \href{http://arxiv.org/abs/1409.6908}{arXiv:1409.6908}, 2014.

\bibitem{Joyc13} D. Joyce, {\it A generalization of manifolds with corners}, \href{http://arxiv.org/abs/1501.00401}{arXiv:1501.00401}, 2015.

\bibitem{Joyc14} D. Joyce, {\it Kuranishi spaces as a\/ $2$-category}, \href{http://arxiv.org/abs/1510.07444}{arXiv:1510.07444}, 2015.

\bibitem{JoSa} D. Joyce and S. Salur, {\it Deformations of asymptotically cylindrical coassociative submanifolds with fixed boundary}, Geom. Topol. 9 (2005), 1115--1146. \href{http://arxiv.org/abs/math/0408137}{math.DG/0408137}.

\bibitem{Kova} A.G. Kovalev, {\it Twisted connected sums and special Riemannian holonomy}, J. Reine Angew. Math. 565 (2003), 125--160. \href{http://arxiv.org/abs/math/0012189}{math.DG/0012189}.

\bibitem{KoNo} A.G. Kovalev and J. Nordstr\"om, {\it Asymptotically cylindrical\/ $7$-manifolds of holonomy $G_2$ with applications to compact irreducible $G_2$-manifolds}, Ann. Global Anal. Geom. 38 (2010), 221--257. \href{http://arxiv.org/abs/0907.0497}{arXiv:0907.0497}.

\bibitem{Kron1} P.B. Kronheimer, {\it The construction of ALE spaces as hyperk\"ahler quotients}, J. Diff. Geom. 29 (1989), 665--683.

\bibitem{Kron2} P.B. Kronheimer, {\it A Torelli-type theorem for gravitational instantons}, J. Diff. Geom. 29 (1989), 685–697.

\bibitem{Lock1} R. Lockhart, {\it Fredholm properties of a class of elliptic operators on a class of noncompact manifolds}, Duke Math. J. 48 (1981), 289--312.

\bibitem{Lock2} R. Lockhart, {\it Fredholm, Hodge and Liouville Theorems
on noncompact manifolds}, Trans. A.M.S. 301 (1987), 1--35.

\bibitem{LoMc} R.B. Lockhart and R.C. McOwen, {\it Elliptic
Differential Operators on Noncompact Manifolds}, Ann. Sc.
norm. sup. Pisa, Classe di scienze 12 (1987), 409--447.

\bibitem{Lota1} J. Lotay, {\it Coassociative $4$-folds with conical singularities}, Comm. Anal. Geom. 15 (2007), 891--946. \href{http://arxiv.org/abs/math/0601762}{math.DG/0601762}.

\bibitem{Lota2} J. Lotay, {\it Deformation theory of asymptotically conical coassociative $4$-folds}, Proc. L.M.S. 99 (2009), 386--424. \href{http://arxiv.org/abs/math/0411116}{math.DG/0411116}.

\bibitem{Lota3} J. Lotay, {\it Desingularization of coassociative $4$-folds with conical singularities}, Geom. Funct. Anal. 18 (2009), 2055--2100. \href{http://arxiv.org/abs/math/0611183}{math.DG/0611183}.

\bibitem{Lota4} J. Lotay, {\it Desingularization of coassociative $4$-folds with conical singularities: obstructions and applications}, Trans. A.M.S. 366 (2014), 6051--6092. \href{http://arxiv.org/abs/1209.5116}{arXiv:1209.5116}.

\bibitem{Loya1} P. Loya, {\it Index theory of Dirac operators on manifolds with corners up to codimension two}, pages 131--166 in J. Gil et al., editors, {\it Aspects of boundary problems in analysis and geometry}, Birkh\"auser, Basel, 2004.

\bibitem{Loya2} P. Loya, {\it The index of b-pseudodifferential operators on manifolds with corners}, Ann. Global Anal. Geom. 27 (2005), 101--133.

\bibitem{LoMe} P. Loya and R.B. Melrose, {\it Fredholm perturbations of elliptic operators on manifolds with corners}, preprint, 2003.

\bibitem{Mars} S.P. Marshall, {\it Deformations of special Lagrangian
submanifolds}, Oxford PhD thesis, 2002.\! Available at \href{http://people.maths.ox.ac.uk/~joyce/theses/theses.html}{\tt http://people.maths.ox.ac.uk/$\sim$joyce/}.\!

\bibitem{Mazz} R. Mazzeo, {\it Resolution blowups, spectral convergence and quasi-asymp\-tot\-ic\-ally conical spaces}, Journ\'ees \'Equations aux D\'eriv\'ees Partielles (2006), Expos\'e VIII.

\bibitem{McSa} D. McDuff and D. Salamon, {\it $J$-holomorphic curves and Symplectic Topology}, second edition, A.M.S., 2012.

\bibitem{Melr1} R.B. Melrose, {\it Pseudodifferential operators, corners and singular limits}, pages 217--234 in Proc. Int. Cong. Math. Kyoto, 1990.

\bibitem{Melr2} R.B. Melrose, {\it Calculus of conormal
distributions on manifolds with corners}, IMRN 1992 (1992), 51--61.

\bibitem{Melr3} R.B. Melrose, {\it The Atiyah--Patodi--Singer Index
Theorem}, A.K. Peters, Wellesley, MA, 1993.

\bibitem{Melr4} R.B. Melrose, {\it Differential Analysis on
Manifolds with Corners}, unfinished book available at \href{http://math.mit.edu/~rbm}{\tt
http://math.mit.edu/$\sim$rbm}, 1996.

\bibitem{MeMe} R.B. Melrose and G. Mendoza, {\it Elliptic operators of totally characteristic type}, MSRI preprint, 1983.

\bibitem{MeNi} R.B. Melrose and V. Nistor, {\it K-Theory of\/ $C^*$-algebras of b-pseu\-do\-diff\-er\-ent\-ial operators}, Geom. Funct. Anal. 8 (1998), 88--122. 

\bibitem{MePi} R.B. Melrose and P. Piazza, {\it Analytic K-Theory on manifolds with corners}, Advances in Math. 92 (1992), 1--26.

\bibitem{Mont} B. Monthubert, {\it Groupoids and pseudodifferential
calculus on manifolds with corners}, J. Funct. Anal. 199 (2003),
243--286.

\bibitem{Nord} J. Nordstr\"om, {\it Deformations of asymptotically cylindrical $G_2$-manifolds}, Math. Proc. Camb. Phil. Soc. 145 (2008), 311--348. \href{http://arxiv.org/abs/0705.4444}{arXiv:0705.4444}.

\bibitem{Paci1} T. Pacini, {\it Deformations of Asymptotically
Conical Special Lagrangian Submanifolds}, Pacific J. Math. 215 (2004), 151--181. \href{http://arxiv.org/abs/math/0207144}{math.DG/0207144}.

\bibitem{Paci2} T. Pacini, {\it Special Lagrangian conifolds, I: moduli spaces}, Proc. L.M.S. 107 (2013), 198--224. \href{http://arxiv.org/abs/1002.1222}{arXiv:1002.1222} (extended version \href{http://arxiv.org/abs/1211.2800}{arXiv:1211.2800}).

\bibitem{Paci3} T. Pacini, {\it Special Lagrangian conifolds, II: gluing constructions in $\C^m$}, Proc. L.M.S. 107 (2013), 225--266. \href{http://arxiv.org/abs/1109.3339}{arXiv:1109.3339}.

\bibitem{Piaz} P. Piazza, {\it On the index of elliptic operators on manifolds with boundary}, J. Functional Anal. 117 (1993), 308--359.

\bibitem{Schw} M. Schwarz, {\it Morse homology}, Progr. Math. 111, Birkh\"auser, Basel, 1993.

\bibitem{Seid} P. Seidel, {\it Fukaya categories and
Picard--Lefschetz theory}, E.M.S., Z\"urich, 2008.

\bibitem{Wehr} K. Wehrheim, {\it Smooth structures on Morse trajectory spaces, featuring finite ends and associative gluing}, pages 369--450 in {\it Proceedings of the Freedman Fest}, Geom. Topol. Monogr. 18, Coventry, 2012. \href{http://arxiv.org/abs/1205.0713}{arXiv:1205.0713}.

\bibitem{Yang} D. Yang, {\it A choice-independent theory of Kuranishi structures
and the polyfold-Kuranishi correspondence}, PhD thesis, New York University, 2014. Available at 
\href{http://webusers.imj-prg.fr/~dingyu.yang/thesis.pdf}{\tt http://webusers.imj-prg.fr/$\sim$dingyu.yang/thesis.pdf}.

\end{thebibliography}
\end{document}